\documentclass[11pt]{article}%
\usepackage{amsfonts,amsmath,amssymb,bookmark,epsfig,graphicx,here,latexsym,mathtools,setspace}
\usepackage[all]{xy}
\usepackage[usenames]{color}

\providecommand{\U}[1]{\protect\rule{.1in}{.1in}}
\providecommand{\U}[1]{\protect\rule{.1in}{.1in}}

\def\<{\langle}
\def\>{\rangle}

\newcommand{\be}{\begin{equation}}
\newcommand{\ee}{\end{equation}}
\newcommand{\bea}{\begin{eqnarray}}
\newcommand{\eea}{\end{eqnarray}}
\newcommand{\beas}{\begin{eqnarray*}}
\newcommand{\eeas}{\end{eqnarray*}}

\newtheorem{theorem}{Theorem}[section]

\newtheorem{corollary}[theorem]{Corollary}

\newtheorem{example}[theorem]{Example}

\newtheorem{lemma}[theorem]{Lemma}

\newtheorem{proposition}[theorem]{Proposition}
\newtheorem{remark}[theorem]{Remark}

\newenvironment{proof}[1][Proof]{\noindent\textbf{#1.} }{\ \rule{0.5em}{0.5em}}
\newenvironment{Proof}[1][Proof]{\noindent\textbf{#1} }{\ \rule{0.5em}{0.5em}}

\DeclareMathOperator{\argmin}{arg\,min}
\DeclareMathOperator{\esssup}{ess\,sup}
\DeclareMathOperator{\essinf}{ess\,inf}

\newcommand{\E}[1]{{\mathbb{E}}\left[#1\right]}

\makeatletter
\newcommand*{\transpose}{%
  {\mathpalette\@transpose{}}%
}
\newcommand*{\@transpose}[2]{%
  \raisebox{\depth}{$\m@th#1\intercal$}%
}
\makeatother

\begin{document}

\title{Robust Multiple Stopping --- A Pathwise Duality Approach
%MOR:
%\thanks{We are very grateful to
%Christian Bender and seminar participants at the University of Waterloo for comments and suggestions.
%This research was funded in part by
%the Netherlands Organization for Scientific Research (NWO) under grants NWO-Vidi and NWO-Vici (Laeven)
%and by the DFG Excellence Cluster Math+ Berlin, project AA4-2 (Schoenmakers).}
}
\author{Roger J. A. Laeven\thanks{Corresponding author.}\\
{\footnotesize Dept. of Quantitative Economics}\\
{\footnotesize University of Amsterdam, EURANDOM}\\
{\footnotesize and CentER}\\
{\footnotesize {\tt R.J.A.Laeven@uva.nl}}\\[1mm]
\and John G. M. Schoenmakers\\
{\footnotesize Stochastic Algorithms and Nonparametric Statistics}\\
{\footnotesize Weierstrass Institute Berlin}\\
{\footnotesize {\tt schoenma@wias-berlin.de}}\\
\and Nikolaus Schweizer\\
{\footnotesize Dept. of Econometrics and Operations Research}\\
{\footnotesize Tilburg University}\\
%{\footnotesize and CentER}\\
{\footnotesize {\tt N.F.F.Schweizer@tilburguniversity.edu}}\\
\and Mitja Stadje\\
{\footnotesize Insurance Science and Financial Mathematics}\\
{\footnotesize Faculty of Mathematics and Economics}\\
{\footnotesize Ulm University}\\
{\footnotesize {\tt Mitja.Stadje@uni-ulm.de}}}
\date{This Version: \today}
\maketitle

\thispagestyle{empty} %MOR: to suppress page number at first page

\begin{abstract}
We develop a method to solve, theoretically and numerically, general optimal stopping problems.
Our general setting allows for multiple exercise rights, i.e., optimal multiple stopping,
for a robust evaluation that accounts for model uncertainty,
and for general reward processes driven by multi-dimensional jump-diffusions.
Our approach relies on first establishing robust martingale dual representation results
for the multiple stopping problem
that satisfy appealing pathwise optimality (i.e., almost sure) properties.
Next, we exploit these theoretical results to develop upper and lower bounds that,
as we formally show, not only converge to the true solution asymptotically,
but also constitute genuine pre-limiting upper and lower bounds.
%with only linear complexity.
We illustrate the applicability of our approach in a few examples
and analyze the impact of model uncertainty on optimal multiple stopping strategies.
\end{abstract}

\noindent\textbf{Keywords:} Optimal stopping; Multiple stopping;
Robustness; Model uncertainty; Ambiguity;
Pathwise duality; $g$-expectations; BSDEs; Regression.\\[3mm]
\noindent\textbf{AMS 2010 Subject Classification:} \textit{Primary}: 49L20; 60G40; 62L15;
\textit{Secondary}: 91B06; 91B16.\\[3mm]
\noindent\textbf{OR/MS Classification:} Dynamic programming/optimal control: Models;
Dynamic programming/optimal control: Applications;
Decision analysis: Risk.

\newpage

\clearpage\pagenumbering{arabic} %MOR: to restart page numbering

\setcounter{equation}{0}

\section{Introduction}

%\doublespacing

%{\footnotesize
%
%\noindent Related work:
%
%\begin{itemize}
%\item Brown, Smith, Sun, Operations Research, 2010.
%
%\item M. Kobylanski, M.-C. Quenez, and E. Rouy-Mironescu. Optimal multiple
%stopping time problem. Ann. Appl. Probab., 21(4):1365-1399, 2011.
%
%\item KLLSS, Mathematics of Operations Research, 2018.
%
%\item Bender, Schweizer, Mathematical Finance, 2017, Mathematics of Operations
%Research, 2018.
%\end{itemize}
%
%\noindent Our contribution:
%
%\begin{itemize}
%\item Elegant a.s. dual representation
%
%\item Pathwise optimality, surely optimal martingales
%
%\item Converging upper bounds (approximated and genuine), lower bounds,
%algorithm exploits stability result
%\end{itemize}
%
%\noindent Generality:
%
%\begin{itemize}
%\item Robustness (ambiguity)
%
%\item Jumps
%
%\item Multiple stop
%\end{itemize}
%
%\noindent Applications and implications:
%
%\begin{itemize}
%\item Examples and Economic analysis
%\end{itemize}
%
%\noindent Examples:
%
%\begin{itemize}
%\item Utility indifference pricing (seller's perspective)
%
%\item Risk measures
%
%\item No-arbitrage pricing
%\end{itemize}
%}
%
%\vskip 0.3cm

In this paper we analyze general optimal stopping problems of the following form:
\begin{equation}
Y^{\ast,L}_{t}:=\sup_{\substack{t\leq\tau_{1}<\cdots<\tau_{L} \\
(\tau_{1},\ldots,\tau_{L})\in\mathcal{T}_{t}(L)}}
\sup_{\mathbb{Q}\in %M\subset
\mathcal{Q}}
\mathbb{E}_{\mathbb{Q}}\left[\sum_{l=1}^{L}H_{\tau_{l}}\Big|\mathcal{F}_{t}\right],\qquad
0\leq t\leq T, \label{eq:problem}%
\end{equation}
where $\mathcal{T}_{t}(L)$ is a family of stopping time vectors,
$L$ is a number of exercise rights,
$T<\infty$ is a fixed time horizon,
%not cemetery time!
$\mathcal{Q}$ is a family of probabilistic models,
and $H$ is a general $\mathcal{F}_{\centerdot}$-adapted reward process.
(The operator $\sup$ is to be understood as $\esssup$ if it applies to an uncountable family of random variables.)
The optimal stopping problem \eqref{eq:problem} features generality along three dimensions:
(i) it allows for optimal multiple stopping (when $L>1$),
(ii) it allows for a robust evaluation that explicitly takes probabilistic model uncertainty (ambiguity) into account
(when $\mathcal{Q}$ is not a singleton), and
(iii) it allows for general reward processes that will be driven by
multi-dimensional jump-diffusion processes.
The process $Y^{\ast,L}_{t}$ is referred to as the upper Snell envelope of $H$ due to $L$ exercise rights after the seminal work of Snell \cite{S52}.
Problems of this type, or special cases thereof, occur naturally in a wide variety of applications
in probability, operations research, economics and finance.

Our aim is to develop upper and lower bounds on $Y^{\ast,L}_{t}$ that satisfy several desirable properties.
%Our approach %to obtaining upper [***ONLY UPPER***] bounds
%relies on
We achieve this by first establishing
suitable martingale dual representations for problem \eqref{eq:problem}
that can be viewed as significant generalizations of the
classical additive dual representations for \textit{standard}
(i.e., $\mathcal{Q}=\{\mathbb{Q}\}$ as opposed to \textit{robust}) optimal stopping problems,
developed independently by Rogers \cite{R02} and Haugh and Kogan \cite{HK04}
(see also the early Davis and Karatzas \cite{DK94}) and their extension to
standard multiple stopping problems in Schoenmakers \cite{S12}.
%%and see Schoenmakers \cite{S12} for the standard multiple stopping case
%%and \cite{KS10} and \cite{KLLSS18} for the robust single stopping case). IS NOT A.S.
(A multiplicative dual representation for standard optimal stopping problems
was proposed by Jamshidian \cite{J07}.)
Our %explicit
dual representations take the %(elegant)
form of an infimum over (robust) martingales,
%(only),
with no appearance of stopping times.

An appealing feature---both theoretically and for numerical stability---of the
dual representations we establish is their pathwise optimality,
i.e., their almost sure property.
Already when $L=1$ these results are new and of independent interest
for robust optimal single stopping.
They are developed here in the general setting of robust optimal multiple stopping
\eqref{eq:problem}.
The almost sure nature of the dual representations suggests that finding
a `good' %(robust)
martingale that is `close' to a `surely optimal' martingale will yield tight and nearly constant upper bounds.
%(Indeed, numerical approaches based on dual representations commonly rely on finding a suitable martingale.)
The target can be the unique (robust) Doob martingale,
to be constructed from an approximation to the upper Snell envelope or,
more generally, a martingale for which the dual representation's infimum is attained
and the almost sure property is satisfied.
While this phenomenon of tightness and constancy is known in the case of standard, non-robust single stopping problems
(i.e., when $\mathcal{Q}=\{\mathbb{Q}\}$ and $L=1$, see %e.g.,
Rogers \cite{Rog10} and
Schoenmakers, Zhang and Huang \cite{SZH13}), %and the references therein),
we will analyze %and exploit
it %in detail
in our general setting of robust %multiple
stopping.
%[***ADAPT.***]
We will show in particular that a low (vanishing in probability) robust variance implies a
tight (converging in $L^{1}$) approximation.
The mathematical details of these results are delicate.

These new theoretical results justify and enable us to next develop a numerically
implementable method to obtain upper and lower bounds to $Y^{\ast,L}_{t}$ with
desirable properties.
Our lower bound, derived from the proposed exercise strategy, will, as we formally show,
not only converge to the
optimal solution asymptotically but also be `biased low' at the pre-limiting level in a Brownian-Poisson filtration.
This is not the case for the initially proposed upper bound: it converges to
the true solution but is not in general `biased high'.
We therefore also develop a second upper bound
that as we prove both converges to the true solution asymptotically
and is biased high in a Brownian-Poisson filtration.
It is based on a Lipschitzian $L^{2}$-approximation,
with a Lipschitz constant that we explicitly derive,
and a suitable (reversed) application of Jensen's inequality.
We will refer to this second upper bound as our \textit{genuine} upper bound.
The computational complexity of both upper bounds is only linear in the number of exercise rights,
and our method does not require nested simulation.
%[***DISCUSS: algorithm exploits stability result.***]

We provide extensive numerical examples, including
single and multiple stopping problems, univariate and multivariate stochastic drivers,
increasing and decreasing reward functions,
and pure diffusion and jump-diffusion models,
to illustrate the applicability and generality of our approach.
They demonstrate that our %pathwise duality
approach yields %(very)
tight upper and lower bounds that, due to almost sure properties,
moreover have %(very)
low standard errors.
They also analyze the impact of multiple vs. single stopping rights
and reveal that employing a robust evaluation that takes ambiguity into account
is highly relevant for optimal stopping, especially in the presence of multiple exercise rights.
%... [***multiple? multi-dimensional? jump?***],
%[***ADD SOME DETAILS.***]

Embedded in a Brownian-Poisson filtration, the problems we analyze are naturally represented as stopping problems with respect to
$g$-expectations (Peng \cite{P97,P04}),
leading to backward stochastic differential equations (BSDEs).
Hence, we explicitly construct novel genuine upper and lower bounds to BSDE solutions
with positively homogeneous convex drivers
in a Brownian-Poisson filtration,
as a development of independent interest.
Bender, Schweizer and Zhuo \cite{BSZ17},
when analyzing solutions to discrete-time (reflected) BS$\Delta$Es
rather than the continuous-time BSDEs we consider,
develop upper and lower bounds applying techniques different from the techniques we employ.
Bender, G\"artner and Schweizer \cite{BGS18} construct Monte Carlo upper and lower bounds
for a class of discrete-time stochastic dynamic programs
which includes %suitable
discretizations of multiple stopping problems.
Our genuine upper and lower bounds apply directly to our original %underlying
continuous-time problem.
Our genuine lower bound takes advantage of an almost sure property of a `second kind'
that we formally establish in order to reduce its variance---`second kind' to distinguish it from the additive dual representation's almost sure property.
This almost sure property entails that the difference between the BSDEs terminal condition and the associated (robust) martingale
is constant almost surely.
Our genuine upper bound for the continuous-time problem is based on forward simulation of an approximate BSDE solution.
The construction is somewhat related to the a posteriori criterion
for error evaluation
introduced in Bender and Steiner \cite{BS13}
%for BS$\Delta$Es
in a Brownian filtration,
and developed here to obtain explicit genuine upper bounds for BSDEs
in a Brownian-Poisson filtration.

The development of numerically implementable methods to obtain approximations to problems
of the type \eqref{eq:problem} but with $\mathcal{Q}$ a singleton (no ambiguity),
$L=1$ (single stopping), and with $H$ multi-dimensional but satisfying strong conditions,
started with the regression-based Monte Carlo methods of Carriere \cite{C96} and Longstaff and Schwartz \cite{LS01};
see also Tsitsiklis and Van Roy \cite{TR00} and Cl\'ement, Lamberton and Protter \cite{CLP02}.
These methods yield lower bounds to $Y^{\ast,1}_{t}$ by approximating the optimal stopping time using regression
and are commonly referred to as ``primal'' approaches.
An important example of a non-regression based primal approach
is the stochastic mesh method of Broadie and Glasserman \cite{BG04}
(see, for further details, Glasserman \cite{G04} and also Belomestny, Kaledin and Schoenmakers \cite{BKS20}).
%They are related to the stochastic mesh method of
%Broadie and Glasserman \cite{BG04} (see, for further details, also Glasserman \cite{G04}).
``Dual'' algorithms that exploit additive dual representations
to numerically compute upper bounds were first proposed by Andersen and Broadie \cite{AB04}
in the standard single stopping problem and were further developed by e.g., Belomestny, Bender and Schoenmakers \cite{BBS09}
to allow for non-nested simulation.
While primal methods rely in a sense on constructing an
appropriate stopping time, dual methods rely on constructing an appropriate martingale.
Brown, Smith and Sun \cite{BSS10} in an innovative paper enlarge
the information on which an exercise decision may depend in dual optimization,
yielding tight upper bounds.
%via (non-nested) simulation.

Model uncertainty, and the distinction between risk
%(probabilistic model given)
and ambiguity,
%(probabilistic model unknown),
has received much
attention in recent years.
Under the Bayesian paradigm, as adopted in Savage's
\cite{S54} subjective expected utility model, this distinction is, in a sense,
nullified, through subjective probabilities resulting from a subjective prior
probability over probabilistic models that quantifies model uncertainty.
A popular approach beyond the Bayesian paradigm is provided by the multiple
priors model of Gilboa and Schmeidler \cite{GS89}, which is a
decision-theoretic formalization of the classical Waldian maxmin decision rule
(Wald \cite{W50}; see also Huber \cite{H81}) and experimentally motivated by
the Ellsberg \cite{E61} paradox.
These models are intimately related to
coherent, convex and entropy convex measures of risk in financial risk
measurement (F\"ollmer and Schied \cite{FS02,FS11}, Frittelli and Rosazza
Gianin \cite{FG02}, Ruszczy\'nski and Shapiro \cite{RS06a,RS06b}, and Laeven
and Stadje \cite{LS13}).
They explicitly recognize that probabilistic models
may be misspecified and are often referred to as \textit{robust} approaches
(Hansen and Sargent \cite{HS07}).
The literature on robust single
stopping theory is rapidly growing; it includes Riedel \cite{R09},
Kr\"atschmer and Schoenmakers \cite{KS10}, Bayraktar, Karatzas and Yao
\cite{BKY10}, Bayraktar and Yao \cite{BY11}, Cheng and Riedel \cite{CR13},
{\O }ksendal, Sulem and Zhang \cite{OSZ14}, Belomestny and Kr\"atschmer
\cite{BK16a,BK16b},
%(For ambiguity aversion with non-dominated families of measures, see
Bayraktar and Yao \cite{BY14,BY15a,BY15b}, Ekren, Touzi, and Zhang
\cite{ETZ14}, Matoussi, Piozin, and Possama\"{\i} \cite{MPP14}, Matoussi,
Possama\"{\i}, and Zhou \cite{MPZ13}, and Nutz and Zhang \cite{NZ15}.
However, numerically implementable methods to solve general optimal stopping problems
of the form \eqref{eq:problem} have not been well-developed as yet.
Kr\"atschmer \textit{et al.} \cite{KLLSS18} propose a numerically
implementable method for single stopping problems under uncertainty in drift
and jump intensity.
Their approach is dual but not path-wise, i.e., it does not rely on a dual representation with the appealing almost sure property,
and cannot handle multiple stopping problems.

The multiple stopping problem can be viewed as $L$ nested single stopping problems,
where the decision-maker first chooses between stopping at time $\tau_{1}$ on the one hand,
thus collecting the reward and entering into a new contract with $L-1$ exercise rights,
and retaining $L$ exercise rights on the other hand, and so on.
Multiple exercise rights occur naturally in many applications across various fields.
For example,
in environmental economics,
a swing option gives the investor the right to change his purchased
energy quantity a number of times per time period;
in finance, a flexible interest rate cap gives the investor
the right to exercise at each interest rate reset date a number of times over the life of the contract;
and in insurance, a partial surrender option provides a payoff to the policyholder
each time he partially surrenders his life insurance contract;
see e.g., Carmona and Dayanik \cite{CD08} and Carmona and Touzi \cite{CT08} and the references therein.
%[***ADD EXAMPLE: ADAPTIVE INVENTORY CONTROL. SEE Brown, Smith and Sun \cite{BSS10}. NOT REALLY MULTIPLE STOP.***]
Kobylanski, Quenez and Rouy-Mironescu \cite{KQR11} analyze the standard multiple stopping
problem (without ambiguity) allowing the payoff to be
a general functional of an ordered sequence of stopping times.
Bender, Schoenmakers and Zhang \cite{BSZ15} develop a dual approach to
generalized multiple stopping problems
with respect to standard conditional expectations
that is intimately related to the information
relaxation approach of Brown, Smith and Sun \cite{BSS10}.
%These methods do not explicitly allow for model uncertainty.
A primal-dual algorithm for standard multiple stopping with respect to standard conditional expectations
in the context of flexible interest rate caps has been proposed in Balder, Mahayni and Schoenmakers \cite{BMS13}.

%General Rewards and drivers. Reward = discounted cash-flow. We allow for both
%simple rewards and general rewards. [***ADD FURTHER DISCUSSION.***]

As an important application, our approach may be used for robust no-arbitrage
pricing (Hansen and Jagannathan \cite{HJ91},
Cochrane and Sa\'{a}-Requejo \cite{CS00})
of American-style derivatives with possibly multiple exercise rights,
via superhedging.
This entails a significant advancement of the standard
approach, where in a usually incomplete market the corresponding stopping
problem is solved with respect to an arbitrarily chosen (local equivalent
martingale) measure.
Our results can also be applied to indifference
valuation (seller's perspective; Carmona \cite{C09}, Laeven and Stadje \cite{LS14}) of general optimally stopped reward processes
under the multiple priors model.
Another application is that
of robust risk measurement (Ben-Tal and Nemirovski \cite{BN98}, Bertsimas and Brown \cite{BB09}, F\"ollmer and Schied \cite{FS11})
to determine e.g., the risk capital required to cover
optimally stopped reward processes.

The remainder of this paper is organized as follows. %[***TO FOLLOW.***]
In Section \ref{sec:robstop} we recall some basic notions, establish some general properties,
introduce the robust optimal multiple stopping problem,
and provide some examples.
In Section \ref{sec:pathdual}, we present our pathwise dual representations
and establish our results on surely optimal (robust) martingales.
In Section \ref{sec:MSAL}, we outline a general primal-dual algorithm and prove its convergence.
Section \ref{sec:BPcon} presents explicit upper and lower bounds in a Brownian-Poisson filtration.
%Section \ref{sec:MS} generalizes our results to multiple stopping problems. [***]
Section \ref{sec:num} provides extensive numerical results.
%Section \ref{sec:con} contains a short conclusion.
All proofs and several auxiliary results are in the Online Appendix.

\setcounter{equation}{0}

\section{Robust Optimal Multiple Stopping}

\label{sec:robstop}

\subsection{Basic Notions and General Properties}

We start by considering a general stochastic setup.
We let $\big(\Omega,\left(\mathcal{F}_{t}\right)_{t\in\left\{0,\ldots,T\right\}},\mathbb{P}\big)$
be a filtered probability space and
let $\mathfrak{X}$ be a linear subspace of $L^{0}(\Omega,\mathcal{F},\mathbb{P})$
with $\mathcal{F}:=\mathcal{F}_{T}$.
We further assume that $\mathfrak{X}$ has a lattice structure,
i.e., $\mathfrak{X}$ is closed under the operations
$\wedge$ (min) and $\vee$ (max), and that $\mathfrak{X}$ contains all
indicator functions $1_{A},$ $A\in\mathcal{F}$.
((In)equalities
between random variables are understood in the $\mathbb{P}$-almost sure sense,
often without explicit mention.)

To represent preferences, we consider
a family of mappings $\rho:=\left(\rho_{t}\right)_{t=0,\ldots,T}$,
\begin{equation*}
\rho_{t}:\mathfrak{X}\rightarrow\mathfrak{X}\cap L^{0}(\Omega,\mathcal{F}_{t},\mathbb{P}).
\end{equation*}
It is referred to as a monotone, regular, recursive, conditional translation invariant
\emph{dynamic monetary utility functional},
henceforth \emph{DMU} for short,
if it satisfies the following conditions:

\begin{enumerate}
\item[(C1)] $\rho_{t}\left(  X\right)  \leq\rho_{t}\left(  Y\right)  $ for all
$X,Y\in\mathfrak{X}$ with $X\leq Y$ %a.s.
and $t\in\left\{  0,\ldots,T\right\}
$ \emph{(monotonicity)}.

\item[(C2)] $\rho_{t}\left(  1_{A}X\right)  =1_{A}\rho_{t}\left(  X\right)  $
for all $X\in\mathfrak{X}$, $A$ $\in\mathcal{F}_{t}$ and $t\in\left\{
0,\ldots,T\right\}  $ \emph{(regularity).}

\item[(C3)] $\rho_{t}=\rho_{t}\circ\rho_{t+1}$ for all $t\in\left\{
0,\ldots,T-1\right\}  $ \emph{(recursiveness).}

\item[(C4)] $\rho_{t}\left(  X+Y\right)  =\rho_{t}\left(  X\right)  +Y$ for
all $X,Y\in\mathfrak{X}$ with $Y$ $\in\mathcal{F}_{t}$ and $t\in\left\{
0,\ldots,T\right\}  $ \emph{(conditional translation invariance).}
\end{enumerate}

As additional %`optional'
properties we consider:

\begin{enumerate}
\item[(P1)] $\rho_{t}(X+Y)\leq\rho_{t}(X)+\rho_{t}(Y)$ for all $X,Y\in\mathfrak{X}$
and $t\in\left\{0,\ldots,T\right\}$ (\emph{subadditivity}).

\item[(P2)] %The functional $\rho$ is \emph{sensitive} if we have
%\begin{equation*}
$\left[  X\leq0\text{ \ \ and \ \ }\rho_{t}(X)\geq0\right]  \Longrightarrow
X=0,\text{ \ \ for all }X\in\mathfrak{X}\text{ and }t\in\left\{
0,\ldots,T\right\}$ (\emph{sensitivity}).
%\end{equation*}

\item[(P3)] %If
$\rho_{t}\left(\lambda X\right) =\lambda\rho_{t}\left(X\right)$ for all $X\in\mathfrak{X},$
$\lambda\geq0$ and $t\in\left\{0,\ldots,T\right\}$
%,
%$\rho$ is said to be
(\emph{positive homogeneity}).
\end{enumerate}
Conditions (C1)--(C4) will always be assumed.
In the sequel, we will mention explicitly which of the properties (P1)--(P3) is required.
Properties (P1) and (P2) also entail the implication
$\left[  X\geq0\text{ \ \ and \ \ }\rho_{t}(X)\leq0\right]  \Longrightarrow X=0$;
see Lemma \ref{lem:P2a} in Appendix \ref{sec:A1}.
DMUs that satisfy (P1)--(P3), in addition to (C1)--(C4), take the form of robust, or worst case, expectations
and have been widely used in applied probability, operations research, economics and finance;
see the references in the Introduction and Section \ref{sec:examples} below.

In this paper, we will frequently use the following implications of (C2) and (C4):
\begin{enumerate}
\item[(C5)] $\rho_{t}\left(0\right)=0$ for all $t\in\left\{0,\ldots,T\right\}$
(\emph{normalization}).

\item[(C6)] $\rho_{t}\left(X\right)=X$ for all $X\in\mathfrak{X}$ with $X$
$\in\mathcal{F}_{t}$ and $t\in\left\{0,\ldots,T\right\}$
(\emph{$\mathcal{F}_{t}$-invariance}).
\end{enumerate}

%For subadditive, sensitive $\rho$ we have the following alternative characterization.
%\begin{enumerate}
%\item[(P2a)] If (P1) applies then sensitivity (P2) of $\rho_{t}$ implies%
%\begin{equation}
%\left[  X\geq0\text{ \ \ and \ \ }\rho_{t}(X)\leq0\right]  \Longrightarrow
%X=0,\text{ \ \ for all }X\in\mathfrak{X,}\text{ and }t\in\left\{
%0,\ldots,T\right\}  . \label{se1}%
%\end{equation}
%\end{enumerate}
%This can be seen as follows. Let $\rho$ be sensitive (P2) and subadditive.
%Suppose $\rho_{i}\left(  Y\right)  \leq0$ and $Y\geq0.$ Then $-Y\leq0$ and so
%by subadditivity,
%\[
%0=\rho_{t}\left(  Y-Y\right)  \leq\rho_{t}\left(  Y\right)  +\rho_{t}\left(
%-Y\right)  \leq\rho_{t}\left(  -Y\right)
%\]
%Hence, by (P2), $-Y=0,$ i.e., $Y=0$ a.s.

%\begin{definition}
Let $\mathcal{H}$ be the set of adapted processes $\left(U_{t}\right)_{t\in\{0,\ldots,T\}}$
such that $U_{t}\in\mathfrak{X}\cap L^{0}(\Omega,\mathcal{F}_{t},\mathbb{P})$.
A process $M=\left(M_{t}\right)_{t\in\{0,\ldots,T\}}\in\mathcal{H}$ is said to be a $\rho$-martingale if
\begin{equation}
M_{t}=\rho_{t}(M_{t+1}),\qquad 0\leq t<T.
\label{eq:rhomartingale}
\end{equation}
%\end{definition}
We present two auxiliary lemmas.
The first lemma provides a generalization of Doob's optional sampling theorem
towards our setup:
\begin{lemma}
\label{lem:Doob} (Doob) Suppose $\rho$ satisfies (C1)--(C4).
Then, for any $\rho$-martingale $M$ %satisfying (C2)--(C4),
and any stopping time $\tau_{i}$, $i\leq\tau_{i}\leq T$, it holds that $\rho_{i}\left(
M_{\tau_{i}}\right)  =M_{i},$ $0\leq i\leq T.$
\end{lemma}
Due to the next lemma, the properties of recursiveness (C3)
and conditional translation invariance (C4) carry over to stopping times, as we will exploit later:
\begin{lemma}
\label{L3} Let $\rho$ satisfy (C1)--(C4), and let $t\in\left\{  0,\ldots,T\right\}$ be fixed.
Consider, for any stopping time $\tau$, $t\leq\tau\leq T$, the functional%
\begin{equation*}
\rho_{\tau}(X):=\sum_{j=t}^{T}1_{\tau=j}\rho_{j}(X).
\end{equation*}
Then, $\rho_{\tau}$ acts from $\mathcal{F}_{T}\rightarrow\mathcal{F}%
_{\tau}\supset\mathcal{F}_{t}$, and
\begin{itemize}
\item[(i)] $\rho_{\tau}$ satisfies $\rho_{t} = \rho_{t}\circ\rho_{\tau}$;
\item[(ii)] $\rho_{\tau}(X+Y) = X+\rho_{\tau}(Y),
\text{ \ for \ }X\in\mathcal{F}_{\tau},\text{ \ \ }Y\in\mathcal{F}_{T}$.
\end{itemize}
\end{lemma}

\subsection{The Stopping Problem}

Consider a fixed adapted reward, or (discounted) cash-flow, process
$H=\left(H_{t}\right)_{t\in\{0,\ldots,T\}}\in\mathcal{H}$
and a DMU decision-maker with $L$ exercise rights
that have to be exercised at different exercise dates.
For each fixed $t$ and $L$, $0\leq t\leq T$,
let $\mathcal{T}_{t}(L)$ be the family of stopping vectors
$\left(\tau_{1},\ldots,\tau_{L}\right)$
such that $\tau_{1}\geq t$ and $\tau_{l}\geq\tau_{l-1}+1$ for all $l$, $1<l\leq L$.
The decision-maker faces the following robust optimal multiple stopping problem:

\begin{equation}
Y_{t}^{\ast, L}:=
%\rho_{t}\left(  \sum_{j=1}^{L}H_{\tau_{j}^{\ast, l}}\right)
%=
\underset{t\leq\tau_{1}<\tau_{2}<\cdots<\tau_{L}}{\esssup }\text{ }\rho_{t}
\left(\sum_{l=1}^{L}H_{\tau_{l}}\right),\qquad t\in\{0,\ldots,T\},
\label{eq:multiplestop}
\end{equation}
for a DMU functional $\rho$ that satisfies (C1)--(C4).
(We note that problem \eqref{eq:multiplestop} is even slightly more general than problem \eqref{eq:problem},
which arises when additionally (P1)--(P3) are satisfied.)
Henceforth, we write $\sup$ (and $\inf$) instead of $\esssup$ (and $\essinf$) for convenience,
understanding that they apply to an uncountable family of random variables.
For a clean formulation of the multiple stopping problem (\ref{eq:multiplestop}),
we extend the cash-flow process by setting $H_{j}\equiv 0$ and $\mathcal{F}_{j}\equiv\mathcal{F}_{T}$, for $j=T+1,T+2,\ldots$.
That is, the subset of rights $l$, $l=2,\ldots,L$, not exercised by time $T$ become valueless.
Hence, for any $\rho$-martingale $M$, $M_{j}=M_{T}$, $j>T$.
%(see also the proof of Theorem \ref{th:Th14} below).
%we introduce a cemetery time $\partial:>T,$ where all rights will be exercised,
%including those that are not exercised up to time $T.$
%Furthermore, by assumption, $H_{\partial}:=0$ and $\mathcal{F}_{\partial}:=\mathcal{F}_{T}$,
%hence $M_{\partial}=M_{T}$ for a $\rho$-martingale $M$

When $L\equiv 1$, the single stopping problem
\begin{equation}
Y_{t}^{\ast}\equiv Y_{t}^{\ast,1}=\sup_{\tau\in\mathcal{T}_{t}}\rho_{t}(H_{\tau}),\qquad
t\in\{0,\ldots,T\},
\label{eq:OST}
\end{equation}
occurs as a special case,
where the family of stopping times $\mathcal{T}_{t}\equiv\mathcal{T}_{t}(1)$
takes values in the set $\left\{t,\ldots,T\right\}$.

The multiple stopping problem can be viewed as $L$ nested single stopping problems with only a single exercise right.
Indeed, setting $Y^{\ast,0}\equiv 0$, $Y^{\ast,1}\equiv Y^{\ast}$
is the upper Snell envelope of $H$ due to a single exercise right.
Then, for multiple exercise rights $L\geq 1$,
$Y^{\ast,L}$ can be viewed as the upper Snell envelope
of the process
\begin{equation*}
H_{t}+\rho_{t}\left(Y_{t+1}^{\ast,L-1}\right),\qquad t\in\{0,\ldots,T-1\},
\end{equation*}
due to only a single exercise right.

Let us denote the set of $\rho$-martingales $M$ with $M_{0}=0$ by $\mathcal{M}_{0}^{\rho}$.
There exists a unique $\rho$-martingale $M^{\ast\rho}\in\mathcal{M}_{0}^{\rho}$
and a non-decreasing predictable $A^{\ast\rho}\in\mathcal{H}$ such that
\begin{equation}
Y_{t}^{\ast}=Y_{0}^{\ast}+M_{t}^{\ast\rho}-A_{t}^{\ast\rho},\qquad
t\in\{0,\ldots,T\},
\label{eq:DoobMeyerdisc}
\end{equation}
which represents the $\rho$-Doob decomposition of
$Y^{\ast}=\left(Y_{t}^{\ast}\right)_{t\in\{0,\ldots,T\}}$.
It is easy to verify that, for $t\in\{0,\ldots,T-1\}$,
\begin{equation}
M_{t+1}^{\ast\rho}-M_{t}^{\ast\rho}=Y_{t+1}^{\ast}-\rho_{t}\left(Y_{t+1}^{\ast}\right),
\text{ \ and \ }
A_{t+1}^{\ast\rho}-A_{t}^{\ast\rho}=Y_{t}^{\ast}-\rho_{t}\left(  Y_{t+1}^{\ast}\right).
\label{eq:DD}
\end{equation}
Henceforth, the $\rho$-martingale $M^{\ast\rho}$ will often be referred to as the $\rho$-Doob martingale
and we often suppress its superscript $\rho$ to simplify notation.

In Appendix \ref{sec:singlestop}, we %recall and
establish some auxiliary results for problem \eqref{eq:OST}
that will be exploited in the proofs of the results that follow.

\subsection{Examples}

\label{sec:examples}

We provide the following examples in which specific versions of the robust optimal multiple stopping problem
of the general form \eqref{eq:multiplestop}
%\eqref{eq:problem}
occur naturally:

\begin{enumerate}
\item[(A.)] No-arbitrage pricing: Let $\mathcal{Q}$ be the set
of local equivalent martingale measures.
(Only if markets are complete $\mathcal{Q}$ is a singleton,
i.e., $\mathcal{Q}=\{\mathbb{Q}\}$.)
Then, the \textit{superhedging} price $\pi^{L}$ of a contract with $L\geq 1$ exercise rights
and associated payoff $\sum_{l=1}^L H_{\tau_l}$ is given by
\begin{equation*}
\pi^{L}=\sup_{\tau_1 < \tau_2 < \cdots < \tau_L}\sup_{\mathbb{Q}\in\mathcal{Q}}
\mathbb{E}_{\mathbb{Q}}\left[\sum_{l=1}^L H_{\tau_l}\right].
\end{equation*}
Many different approaches to no-arbitrage pricing have been proposed in the
literature; see, e.g., the good-deal bounds of Cochrane and Sa\'{a}-Requejo \cite{CS00},
Hansen and Jagannathan \cite{HJ91} and Bj\"ork and Slinko \cite{BS06},
or the acceptable opportunities of Carr, Geman and Madan \cite{CGM01}.
All these approaches yield prices of the form
\begin{equation*}
\tilde{\pi}^{L}=\sup_{\tau_1 < \tau_2 < \cdots < \tau_L}\sup_{\mathbb{Q}\in\mathcal{Q}_{\mathrm{restricted}}}
\mathbb{E}_{\mathbb{Q}}\left[\sum_{l=1}^L H_{\tau_l}\right],
\end{equation*}
where $\mathcal{Q}_{\mathrm{restricted}}\subset\mathcal{Q}$.

Prototypical situations leading to single and multiple stopping problems in economics and finance
are the pricing and exercising of American-style, Bermudan-style, and swing options.
American options give the holder the %(single)
right to exercise the option (once) %(e.g., call or put)
on any preferred trading day before expiration.
Different from American options, Bermudan options prescribe a set of trading days %before expiration
on which the option can be exercised (once).
Swing options, more generally, give the holder the right to exercise the option multiple times,
at a pre-specified set of exercise dates.
With $L\geq 1$ exercise rights, exercised at $\tau_{1}<\tau_{2}<\cdots<\tau_{L}$,
the payoff equals $\sum_{l=1}^{L}H_{\tau_{l}}$ for a cash-flow process $H\in\mathcal{H}$.
Swing options are particularly popular in energy markets
to manage the risk of fluctuations in oil, gas, or electricity prices.

\item[(B.)] Indifference valuation---the seller's perspective:
Suppose that the seller of a contract has a max-min utility functional
of the form
\begin{equation*}
U(H)=\inf_{\mathbb{Q} \in \mathcal{Q}} \mathbb{E}_{\mathbb{Q}}[H],
\end{equation*}
for a family of probabilistic models (i.e., priors) $\mathcal{Q}$
and adopts a utility indifference valuation approach (Carmona \cite{C09}, Laeven and Stadje \cite{LS14}).
Then, the value $V^{L}$ of a contract with $L\geq 1$ exercise rights
and associated payoff $\sum_{l=1}^L H_{\tau_l}$
is determined from the indifference relation
\begin{equation*}
U(0)=\inf_{\tau_1 < \tau_2 < \cdots < \tau_L} U\left(-\sum_{l=1}^L H_{\tau_l}+ V^{L}\right)
=\inf_{\tau_1 < \tau_2 < \cdots < \tau_L}\inf_{\mathbb{Q} \in \mathcal{Q}}
\mathbb{E}_{\mathbb{Q}}\left[-\sum_{l=1}^L H_{\tau_l}+V^{L}\right].
\end{equation*}
Hence,
\begin{equation*}
V^{L}=\sup_{\tau_1 < \tau_2 < \cdots < \tau_L}\sup_{\mathbb{Q} \in \mathcal{Q}}
\mathbb{E}_{\mathbb{Q}}\left[\sum_{l=1}^L H_{\tau_l}\right].
\end{equation*}

\item[(C.)] Robust risk measurement:
Suppose that $\rho$ is a robust, or worst case, expectation, that is,
\begin{equation}
\rho(H)=\sup_{\mathbb{Q} \in \mathcal{Q}} \mathbb{E}_{\mathbb{Q}}[H],
\label{eq:coherentrm}
\end{equation}
for a family of probabilistic models $\mathcal{Q}$.
In financial risk measurement, \eqref{eq:coherentrm} is referred to as a coherent risk measure
and $\mathcal{Q}$ as a set of generalized scenarios
(Artzner \textit{et al.} \cite{ADEH99} and F\"ollmer and Schied \cite{FS11});
see also Ben-Tal and Nemirovski \cite{BN98} for the intimately connected
robust optimization paradigm.
It determines the minimal amount of risk capital required to be added to the financial position $H$
to make it `safe' from the viewpoint of the regulatory authority.
Applications of coherent risk measures and generalized scenarios to
decision and optimization include
Lesnevski, Nelson and Staum \cite{LNS07},
Bertsimas and Brown \cite{BB09},
%Ben-Tal, Bertsimas and Brown \cite{BBB10},
Choi, Ruszczy\'nski and Zhao \cite{CRZ11},
Philpott, de Matos and Finardi \cite{PMF13}
and Tekaya, Shapiro, Soares and da Costa \cite{TSPP13}.
Assume now that $H_{\tau_{l}}$ is a payout obligation (i.e., liability)
at time $\tau_{l}$,
where $\tau_{l}$ is a stopping time,
due to e.g., a flexible interest rate cap in interest rate markets
%a swing option in energy markets,
or a partial surrender option in life insurance
to be paid to a policyholder who decides to partially surrender his insurance contract.
Then, the required amount of risk capital due to $L\geq 1$ stopping rights is given by
\begin{equation*}
\sup_{\tau_1 < \tau_2 < \cdots < \tau_L}\rho\left(\sum_{l=1}^L H_{\tau_l}\right)
=\sup_{\tau_1 < \tau_2 < \cdots < \tau_L}\sup_{\mathbb{Q} \in \mathcal{Q}}
\mathbb{E}_{\mathbb{Q}}\left[\sum_{l=1}^L H_{\tau_l}\right].
\end{equation*}
\end{enumerate}

\setcounter{equation}{0}

\section{Pathwise Duality}

\label{sec:pathdual}

\subsection{Pathwise Dual Representation}

\label{sec:pathdualrep}

The following theorem establishes our pathwise (i.e., almost sure)
additive dual representation for general multiple stopping problems of the form \eqref{eq:multiplestop}.
\begin{theorem}
\label{th:mulsub}
Suppose %the DMU functional
$\rho$ satisfies (C1)--(C4) and is subadditive (P1). %and sensitive (P2).
Then,
for any adapted process $H=\left(H_{t}\right)_{t\in\{0,\ldots,T\}}\in\mathcal{H}$
and each fixed $t\in\left\{ 0,\ldots,T\right\}$,
\begin{itemize}
\item[(i)] we have the dual representation
\begin{equation}
Y_{t}^{\ast, L}=\underset{M^{(1)},\ldots,M^{(L)}\in\mathcal{M}_{0}^{\rho}}{\inf}\rho_{t}\left(\max_{t\leq j_{1}<j_{2}<\cdots< j_{L}}\sum_{k=1}^{L}
\left(H_{j_{k}}+M_{j_{k-1}}^{(k)}-M_{j_{k}}^{(k)}\right)\right);
\end{equation}
\item[(ii)] the dual representation's infimum is attained:
\begin{equation}
Y_{t}^{\ast, L}=\rho_{t}\left(\max_{t\leq j_{1}<j_{2}<\cdots< j_{L}}
\sum_{k=1}^{L}\left(H_{j_{k}}+M_{j_{k-1}}^{\ast, L-k+1}-M_{j_{k}}^{\ast, L-k+1}\right)\right);
\end{equation}
\item[(iii)] if in addition $\rho$ is sensitive (P2), we have the pathwise dual representation
\begin{equation}
Y_{t}^{\ast, L}=\max_{t\leq j_{1}<j_{2}<\cdots< j_{L}}
\sum_{k=1}^{L}\left(H_{j_{k}}+M_{j_{k-1}}^{\ast, L-k+1}-M_{j_{k}}^{\ast, L-k+1}\right),\qquad\mathrm{almost\ surely};
\end{equation}
\end{itemize}
where the $\rho$-martingales $M^{\ast, L-k+1}$ satisfy
\begin{equation}
M_{r+1}^{\ast, L-k+1}-M_{r}^{\ast, L-k+1}=Y_{r+1}^{\ast, L-k+1}-\rho_{r}
\left(Y_{r+1}^{\ast, L-k+1}\right),
\label{eq:DoobM}
\end{equation}
and $Y^{\ast, L-k+1}$ is the upper Snell envelope due to $L-k+1$ exercise rights,
satisfying the Bellman principle,
\begin{equation}
Y_{r}^{\ast, L-k+1}=\max\left[  H_{r}+\rho_{r}\left(Y_{r+1}^{\ast,L-k}\right),
\rho_{r}\left(Y_{r+1}^{\ast, L-k+1}\right)\right].
\label{eq:Bellman}%
\end{equation}
\end{theorem}

Already when $L\equiv 1$, Theorem \ref{th:mulsub} is new and of significant independent interest.
In this case it simplifies to:
\begin{corollary}
\label{P1} Suppose %the DMU functional
$\rho$ satisfies (C1)--(C4) and
(P1). %(subadditivity).
Then, for any adapted process $H\in\mathcal{H}$ and each
fixed $t\in\left\{ 0,\ldots,T\right\}$, we have the %alternative
dual representation
\begin{align}
Y_{t}^{\ast}  &  =\inf_{M\in\mathcal{M}_{0}^{\rho}}\rho_{t}\left(  \max_{t\leq
j\leq T}\left(  H_{j}+M_{t}-M_{j}\right)  \right) \label{as0}\\
&  =\rho_{t}\left(  \max_{t\leq j\leq T}\left(  H_{j}+M_{t}^{\ast}-M_{j}%
^{\ast}\right)  \right), \label{as00}%
\end{align}
where $M^{\ast}$ %(and $A^{\ast}$) are the ingredients of
is the $\rho$-Doob martingale in Eqn.~\eqref{eq:DoobMeyerdisc}.
If in addition $\rho$ is sensitive (P2), we have
the almost sure property:
\begin{equation}
Y_{t}^{\ast}=\max_{t\leq j\leq T}\left(  H_{j}+M_{t}^{\ast}-M_{j}^{\ast
}\right),\qquad\mathrm{almost\ surely.} \label{as}
\end{equation}
\end{corollary}

\begin{remark}
We note that the %robust
single stopping problem also admits
an alternative but non-pathwise additive dual representation
for functionals $\rho$ satisfying (C1)--(C4);
%without requiring the functional $\rho$ to be subadditive (P1); %and sensitive (P2);
see Proposition~\ref{prop:ADP} in the Appendix.
In Theorem~\ref{th:mulsub} and Corollary~\ref{P1}, the subadditivity property (P1) of $\rho$ is required,
and exploited through application of Lemma~\ref{L2} in the proof of Theorem~\ref{th:mulsub}.
%In particular, the almost sure representation (ii) is closely connected with (P1) and (P2).
A dual representation theorem in the spirit of Theorem~\ref{th:mulsub} and Corollary~\ref{P1}
without assuming (P1) seems not possible to us.
\end{remark}

\subsection{Surely Optimal \texorpdfstring{$\rho$-}{rho-}Martingales}

\label{sec:suroptmart}

The $\rho$-Doob martingale in Eqn.~\eqref{eq:DoobMeyerdisc} plays a special role
in the set of $\rho$-martingales $\mathcal{M}_{0}^{\rho}$
as its appearance in Corollary~\ref{P1}, Eqn.~\eqref{as00}
(and indirect appearance in Theorem~\ref{th:mulsub}, (ii))
confirms.
In our numerically implementable method developed and applied in Sections~\ref{sec:MSAL}--\ref{sec:num}
we rely on the $\rho$-Doob martingale.
From a theoretical perspective, however,
and as a general justification of our pathwise dual, martingale-based approach,
we develop in this section several results
on so-called \textit{surely optimal $\rho$-martingales}.
To achieve this, we generalize the concept of standard surely optimal martingales
(see Schoenmakers, Zhang and Huang \cite{SZH13}
in the context of standard conditional expectations and optimal single stopping problems)
to subadditive DMU functionals.
The results in this section show formally that
if a general $\rho$-martingale---not necessarily the $\rho$-Doob martingale---induces `small' (robust) variance,
then the associated bounds obtained from the dual representation can be expected to be `tight'
and nearly constant.

Our results on surely optimal $\rho$-martingales can also serve as a diagnostic device
to assess the quality of the estimated $\rho$-Doob martingale,
derived from an (input) approximation to the upper Snell envelope.
If the (robust) variance the estimate induces fails to be small, then it must be far from the $\rho$-Doob martingale.
If, on the other hand, this variance is small, then the estimate will be close to an optimal $\rho$-martingale
(attaining the dual representation's infimum),
even though not necessarily close to the $\rho$-Doob martingale.

For ease of exposition, we focus attention first on optimal single stopping problems.
The next theorem generalizes the analogous key measurability result for standard conditional expectations and optimal single stopping problems
to DMU functionals satisfying (C1)--(C4) and (P1).
\begin{theorem}
\label{th:sure1}
Let $Y_{i}^{\ast}$ be the upper Snell envelope of the cash-flow process $H$
with respect to a subadditive DMU functional $\rho$ satisfying (C1)--(C4) and (P1) as in Corollary~\ref{P1}
and let $M$ be a $\rho$-martingale.
Then, for any $i\in{\{0,\ldots,T\}}$,
\begin{equation*}
\max_{i\leq j\leq T}(H_{j}-M_{j}+M_{i})\in\mathcal{F}_{i}\Rightarrow
\max_{i\leq j\leq T}(H_{j}-M_{j}+M_{i})=Y_{i}^{\ast}.
\end{equation*}
\end{theorem}

The following lemma will later allow for a generalization of the results in this section to multiple stopping.
\begin{lemma}
\label{th:sure1_cor}
Let $Y^{\ast}$, $H$, $M$ and $\rho$ be as in Theorem~\ref{th:sure1}.
Then, for any fixed $0\leq i<T$,
\begin{align*}
\theta_{i+} &  :=\max_{i<j\leq T}(H_{j}-M_{j}+M_{i})\in\mathcal{F}%
_{i}\Rightarrow\\
(i)\text{ \ }\theta_{i+} &  =\rho_{i}\left(  Y_{i+1}^{\ast}\right)  \text{
\ \ and \ \ }(ii)\text{ \ }M_{i+1}-M_{i}=Y_{i+1}^{\ast}-\rho_{i}\left(
Y_{i+1}^{\ast}\right),\text{ }%
\end{align*}
hence $M_{i+1}-M_{i}$ is a $\rho$-Doob martingale increment.
(Note the strict first inequality under the max operator.)
Thus, in particular, if
$\theta_{i+}\in\mathcal{F}_{i}$ for every $0\leq i<T,$ then $M$ is the $\rho$-Doob martingale.
\end{lemma}

Let us define the conditional $\rho$-variance as follows:
\begin{equation}
\mathrm{Var}_{\rho_{i}}\left(  X\right)  :=\rho_{i}\left(  \left(  X-\rho
_{i}\left(  X\right)  \right)  ^{2}\right). \label{cv}%
\end{equation}
It admits a conditional Chebyshev inequality,
exploited in the proof of Theorem \ref{th:Th14} below,
as follows:
\begin{proposition}
\label{prop:ch}
Assume (C1)--(C4).
If $\rho$ is positively homogeneous (P3),
then
\begin{equation}
\rho_{i}\left(  1_{\left\vert X-\rho_{i}\left(  X\right)  \right\vert
\geq\epsilon}\right)  \leq\frac{\mathrm{Var}_{\rho_{i}}\left(  X\right)
}{\epsilon^{2}}. \label{ch}
\end{equation}
\end{proposition}
Next, we state the following lemma:
\begin{lemma}
\label{sensv}
Assume (C1)--(C4).
Let $\rho$ be subadditive (P1) and sensitive (P2).
Then,
\begin{equation*}
\mathrm{Var}_{\rho_{i}}\left(  X\right)  =0\Longleftrightarrow\text{ }%
X\in\mathcal{F}_{i}.
\end{equation*}
\end{lemma}

By virtue of Lemma~\ref{sensv}, Theorem~\ref{th:sure1} implies that if a $\rho$-martingale
$M$ is such that, for some $i\leq j\leq T$, the conditional $\rho$-variance
\begin{equation*}
\mathrm{Var}_{\rho_{i}}\left(  \theta_{i}(M)\right)  :=\mathrm{Var}_{\rho_{i}%
}\left(  \max_{i\leq j\leq T}(H_{j}-M_{j}+M_{i})\right)
\end{equation*}
is zero a.s., then $\theta_{i}(M)=Y_{i}^{\ast}$ a.s.
In that case, we say that the $\rho$-martingale $M$ is \textit{surely optimal} at $i$.
(Note that, in particular, the $\rho$-Doob martingale in (\ref{eq:DoobMeyerdisc}) is surely optimal.)

We then present a stability result for $\rho$-martingales $M$ that are,
in loose terms, `close' to be surely optimal, %at some $i$,
in the sense that the conditional $\rho$-variance $\mathrm{Var}_{\rho_{i}}(\theta_{i}(M))$ is `small'.
In particular, for a sequence of $\rho$-martingales $(M^{(n)})_{n\geq1}$
that induces vanishing conditional $\rho$-variance,
we establish weak conditions
guaranteeing that the corresponding upper bounds converge to the upper Snell
envelope (in $L^{1}$),
even though the sequence of $\rho$-martingales $(M^{(n)})_{n\geq1}$ itself does not necessarily converge.

\begin{theorem}
\label{th:Th14}
Assume (C1)--(C4).
Let $\rho$ be subadditive (P1) and positively homogeneous (P3).
Suppose that
\begin{equation*}
\mathrm{Var}_{\rho_{i}}\left(  \theta_{i}^{(n)}\right)  \overset{\text{P}}{\rightarrow}0,\quad
\mathrm{with}\quad
\theta_{i}^{(n)}=\max_{i\leq j\leq T}\left(  H_{j}-M_{j}^{(n)}+M_{i}^{(n)}\right).
\end{equation*}
If, in addition, for every $i$ and every $\epsilon>0$ there exists
$K_{\epsilon}>0$ such that
\begin{equation}
\sup_{n\geq1}\mathbb{E}\left[\rho_{i}\left(  \left\vert M_{i}^{(n)}\right\vert
1_{\left\vert M_{i}^{(n)}\right\vert >K_{\epsilon}}\right)\right]  <\epsilon,
\label{gUI}%
\end{equation}
then
\begin{equation*}
\rho_{i}\left(  \theta_{i}^{(n)}\right)  \overset{L_{1}}{\longrightarrow}Y_{i}^{\ast}.
\end{equation*}
\end{theorem}

Note that, if $\rho_{i}\equiv\mathbb{E}_{i}$, \eqref{gUI} boils down to a standard
uniform integrability condition.
More generally, we have the following:

\begin{proposition}
\label{prop:EM}
Assume (C1)--(C4).
Let $\rho$ be subadditive (P1) and positively homogeneous (P3).
If, for some $\eta>0$,
\begin{equation*}
\sup_{n\geq1}\mathbb{E}\left[\rho_{i}\left(\left\vert M_{i}^{(n)}\right\vert
^{1+\eta}\right)\right]  <\infty,
\end{equation*}
then $\left(M_{i}^{(n)}\right)_{n\geq1}$ satisfies (\ref{gUI}).
\end{proposition}

%\paragraph{$\rho_{t}$-variance}
Under an additional Lipschitz continuity condition,
Theorem~\ref{th:Th14} may be readily applied as follows.
Let us assume that, for some number $p$,
\begin{equation}
\mathbb{E}\left[\left\vert \rho_{i}\left(  Z\right)  \right\vert ^{p}\right]\leq
C_{p}\mathbb{E}\left[  \left\vert Z\right\vert ^{p}\right],  \label{l1}%
\end{equation}
with $C_{p}>0.$
In particular, if (\ref{l1}) holds for $p=1$, one obviously has
\begin{equation}
\mathbb{E}\left[\mathrm{Var}_{\rho_{i}}\left(X\right)\right]:=\mathbb{E}\left[
\rho_{i}\left(\left(X-\rho_{i}\left(X\right)\right)^{2}\right)
\right]  \leq C_{1}\widetilde{\mathrm{Var}}_{i}\left(X\right),
\label{vc1}
\end{equation}
where
%\begin{equation}
$\widetilde{\mathrm{Var}}_{i}(X):=\mathbb{E}_{\mathcal{F}_{i}}\left[  \left(
X-\rho_{i}\left(  X\right)  \right)  ^{2}\right]$.
%\label{evar}
%\end{equation}
That is, if we achieve in an algorithm that
%\begin{equation*}
$\widetilde{\mathrm{Var}}_{i}\left(\theta_{i}^{(n)}\right)
\overset{\text{P}}{\rightarrow}0$,
%\end{equation*}
and the $\left(  M_{i}^{(n)}\right)_{n\geq1}$ are standard uniformly
integrable, i.e.,
\begin{equation*}
\sup_{n\geq1}\mathbb{E}\left[  \left\vert M_{i}^{(n)}\right\vert 1_{\left\vert
M_{i}^{(n)}\right\vert >K_{\epsilon}}\right]  <\epsilon,
\end{equation*}
then on the one hand,
\begin{align*}
\sup_{n\geq1}\mathbb{E}\left[\rho_{i}\left(  \left\vert M_{i}^{(n)}\right\vert
1_{\left\vert M_{i}^{(n)}\right\vert >K_{\epsilon}}\right)\right]   &  \leq C_{1}%
\sup_{n\geq1}\mathbb{E}\left[\mathbb{E}_{i}\left\vert M_{i}^{(n)}\right\vert 1_{\left\vert
M_{i}^{(n)}\right\vert >K_{\epsilon}}\right]\\
&  =C_{1}\sup_{n\geq1}\mathbb{E}\left[\left\vert M_{i}^{(n)}\right\vert
1_{\left\vert M_{i}^{(n)}\right\vert >K_{\epsilon}}\right]<\epsilon,
\end{align*}
i.e., the $\left(  M_{i}^{(n)}\right)_{n\geq1}$ satisfy the notion of
$\rho_{i}$-uniform integrability,
and on the other hand we have due to (\ref{vc1}) that
\begin{equation*}
\mathrm{Var}_{\rho_{i}}\left(  \theta_{i}^{(n)}\right)  \overset{\text{P}%
}{\rightarrow}0.
\end{equation*}
Then, Theorem~\ref{th:Th14} implies that
\begin{equation*}
\rho_{i}\left(  \theta_{i}^{(n)}\right)  \overset{L_{1}}{\rightarrow}%
Y_{i}^{\ast}.
\end{equation*}

The next theorem generalizes Theorem~\ref{th:sure1} and Lemma~\ref{th:sure1_cor} to multiple stopping:
\begin{theorem}
\label{th:sure1_cor_mult}
Assume (C1)--(C4) and (P1).
Let us define for a set of $\rho$-martingales $M^{(k)}$, $k=1,\ldots,L$,
\begin{equation*}
\Theta_{i+}^{q}:=\max_{i<j_{1}<j_{2}<\cdots<j_{q}}\sum_{k=1}%
^{q}\left(H_{j_{k}}+M_{j_{k-1}}^{(q-k+1)}-M_{j_{k}}^{(q-k+1)}\right)\text{ \ \ for
\ \ }q=1,\ldots,L,
\end{equation*}
with $j_{0}:=i.$
(Note the strict first inequality under the max operator.)
Then it holds that
\begin{align*}
\Theta_{i+}^{q}  &  \in\mathcal{F}_{i}\text{ \ \ for \ \ }q=1,\ldots,L,\text{
\ \ }0\leq i<T\text{ \ \ \ }\Longrightarrow\\
&\left\{
\begin{array}
[c]{l}%
(i)\text{ \ \ }\Theta_{i+}^{q}=\rho_{i}\left(  Y_{i+1}^{\ast, q}\right) \\
(ii)\text{ \ }M_{i+1}^{(q)}-M_{i}^{(q)}=Y_{i+1}^{\ast, q}-\rho_{i}\left(
Y_{i+1}^{\ast, q}\right)
\end{array}
\right.  \text{for \ \ }q    =1,\ldots,L,\text{ \ \ }0\leq i<T.
\end{align*}
\end{theorem}

\begin{remark}
Without doubt it is also possible to derive a version of Theorem~\ref{th:Th14} for the multiple stopping setting.
However, as our algorithm in Section~\ref{sec:MSAL} below aims at approximative construction
of $\rho$-Doob martingale increments associated with
the upper Snell envelopes $Y^{\ast, l}$ of the
generalized cash-flows
\begin{equation}
U_{j}^{\ast, l}:=H_{j}+\rho_{j}\left(  Y_{j+1}^{\ast, l-1}\right)  ,\text{
\ \ }l=1,\ldots,L, \label{genc0}%
\end{equation}
respectively, rather than approximative construction of merely surely optimal
$\rho$-martingales, we refrain from such an analysis.
\end{remark}

\setcounter{equation}{0}

\section{A General Primal-Dual Pseudo Algorithm}% for Robust Multiple Stopping}

\label{sec:MSAL}

In this section, we develop a primal-dual pseudo algorithm for robust multiple stopping (henceforth called \textit{algorithm} for short).
Our treatment in this section applies to DMUs satisfying (C1)--(C4), (P1) and weak continuity conditions,
and to general reward processes in a Markovian environment;
in particular, our treatment in this section is not restricted to $g$-expectations.
The following lemma will serve as a cornerstone in our construction.

\begin{lemma}
\label{var+}
Let $\rho$ satisfy (C1)--(C4), (P1)
and be Lipschitz continuous in the sense of (\ref{l1}) for $p=2$.
Furthermore, let $\mathcal{C}^{N},\mathcal{C},\mathcal{U}\in\mathcal{F}_{j}$,
$\mathcal{Y}\in\mathcal{F}_{j+1}$,
and let $\mathfrak{m}^{N}\in\mathcal{F}_{j+1}$ be a $\rho$-martingale increment,
that is, $\rho_{j}\left(\mathfrak{m}^{N}\right)=0,$ for $j=0,\ldots,T$, %[***OR: $j=0,\ldots,T-1$***],
$N\in\mathbb{N}$,
such that
\begin{align}
\mathbb{E}\left[  \left(  \mathcal{Y}-\mathfrak{m}^{N}-\mathcal{C}^{N}\right)
^{2}\right]   &  \rightarrow0,\text{ \ \ and}\label{coy}\\
\mathbb{E}\left[  \left(  \mathcal{C}^{N}-\mathcal{C}\right)  ^{2}\right]   &
\rightarrow0,\text{ \ \ for }N\rightarrow\infty.\nonumber
\end{align}
Then,
\begin{equation*}
\mathcal{C}=\rho_{j}\left(  \mathcal{Y}\right)  \text{ \ \ and \ \ }%
\mathfrak{m}^{N}\overset{L_{2}}{\rightarrow}\mathcal{Y}-\rho_{j}\left(
\mathcal{Y}\right).
\end{equation*}
\end{lemma}

\begin{corollary}
\label{loc}
Let $Y_{j+1}^{l}$ be an already constructed approximation to a random variable $Y_{j+1}^{\ast, l}$.
Furthermore, let $\rho$,
the $\rho$-martingale increment $\mathfrak{m}_{j+1}^{l,N}\in\mathcal{F}_{j+1}$,
and $\mathcal{C}_{j}^{l,N},\mathcal{C}_{j}^{l}\in\mathcal{F}_{j}$ be as in Lemma~\ref{var+},
such that
\begin{align*}
\mathbb{E}\left[  \left(  Y_{j+1}^{l}\,-\mathfrak{m}_{j+1}^{l,N}%
-\mathcal{C}_{j}^{l,N}\right)  ^{2}\right]\rightarrow0,
\quad\text{and}\quad
\mathbb{E}\left[  \left(  \mathcal{C}_{j}^{l,N}-\mathcal{C}_{j}^{l}\right)
^{2}\right]\rightarrow0,
\quad\text{for}\ N\rightarrow\infty.
\end{align*}
Then,
\begin{equation*}
Y_{j+1}^{l}\,-\mathfrak{m}_{j+1}^{l,N}\overset{L_{2}}{\rightarrow}\rho
_{j}\left(  Y_{j+1}^{l}\right).
\end{equation*}
\end{corollary}

Guided by Lemma~\ref{var+} and Corollary~\ref{loc}, we now develop
a primal-dual algorithm in the context of a Markovian underlying process $X$
with state space $\mathbb{R}^{d}$, possibly in continuous time, that is
monitored at the exercise dates as $X_{j},$ $j=0,\ldots,T$.
As usual, we assume that $\mathcal{F}_{j}$ is the $\sigma$-field generated
(directly or, as in the next section, indirectly)
by the process
$X$ up to exercise date $j$.
Furthermore, we assume that the cash-flows are of the form
\begin{equation*}
H_{j}=f_{j}(X_{j}),\text{ \ \ for \ \ }j=0,\ldots,T,
\end{equation*}
where $f_{j}:\mathbb{R}^{d}\rightarrow\mathbb{R}_{\geq0}$, $j=0,\ldots,T$, are
given nonnegative payoff functions
such that $H\in\mathcal{H}$.
Note that, due to the Bellman principle (\ref{eq:Bellman}),
$Y^{\ast, l}$ can be seen as the upper Snell envelope corresponding
to the generalized cash-flow
\begin{equation}
U_{j}^{\ast, l}:=H_{j}+\rho_{j}\left(Y_{j+1}^{\ast, l-1}\right)
=:\ H_{j}+c_{j}^{\ast, l}(X_{j}),\text{\ \ }l=1,\ldots,L,\quad j=0,\ldots,T-1,
\label{genc}
\end{equation}
due to a single exercise right,
where the so-called continuation functions $c_{j}^{\ast, l}$ exist by
Markovianity. %Let us further introduce for $0\leq l\leq L,$ approximate
%continuation functions $c_{r}^{l}(\cdot)$ defined by%
%\begin{equation}
%c_{r}^{l}(X_{r}):=\rho_{r}\left(  Y_{r+1}^{l}\right)  ,\text{ \ \ }1\leq r\leq
%T. \label{continuation}
%\end{equation}
We also assume to have a set of Monte Carlo simulated training
trajectories
\begin{equation*}
X^{(n)}\equiv X^{n},\qquad n=1,\ldots,N.
\end{equation*}

We proceed in the following steps:
\begin{itemize}
\item[1.] Initialize $\overline{M}^{0}$ $=$ $\overline{Y}^{0}$
$=\overline{c}^{0}$ $=$ $0$, for $l=0$.
\item[2.] Suppose that, for a particular $l$ with $0\leq l<L$:\newline
(i) we have constructed a set of (approximate) continuation functions
$\overline{c}_{j}^{l}:\mathbb{R}^{d}\rightarrow\mathbb{R}_{\geq0},$ $1\leq j\leq T$,
hence an (approximate) continuation value process (for up to $l$ exercise times)
of the form
\begin{equation*}
\overline{C}_{j}^{l}=\overline{c}_{j}^{l}(X_{j});
\end{equation*}
(ii) we have constructed a (true) $\rho$-martingale $\overline{M}_{j}^{l}$;
and \newline
(iii) we have constructed, on each trajectory $n$, a path
\begin{align}
%\overline{Y}_{r}^{l,n}\,  & :=\overline{Y}_{r+1}^{l,n}-\overline{M}_{r+1}^{l,n} +\overline{M}_{r}^{l,n} \mbox{ if } r \notin \{0,1,\ldots,T\}, \label{barzero}\\
\overline{Y}_{j}^{l,n}\,  &  :=\max\left[ \overline{U}_{j}^{l,n}, \overline{c}_{j}^{l}(X_j^n) \right]
 \mbox{ if } j \in \{0,1,\ldots,T\}, \qquad\text{\ where }  \label{bar}\\
\overline{U}_{j}^{l,n}  &  :=\left\{
\begin{array}
[c]{l}%
f_{j}(X_{j}^{n})+\overline{c}_{j}^{l-1}(X_{j}^{n})\text{ \ \ if
\ \ \ }l>0\\
0\text{ \ \ \ if \ \ }l=0
\end{array}
\right.  ,\qquad 0\leq j\leq T,  %\\
%\text{and } \overline{c}_r^l & := \overline{Y}_{r+1}^{l,n} -\overline{M}_{r+1}^{l,n} + \overline{M}_{r}^{l,n} \nonumber
\end{align}
as an approximation to $Y^{\ast, l,n}$.
\item[3.] Now construct, using these trajectories,
a subsequent (true) $\rho$-martingale $\overline{M}^{l+1}$,
a subsequent set of continuation functions $\overline{c}_{j}^{l+1}$, $j=0,\ldots,T$,
and subsequent trajectories $\overline{Y}^{l+1,n}$, $n=1,\ldots,N$
(as approximations to $M^{\ast, l+1}$ and $Y^{\ast, l+1}$, respectively)
such that (\ref{bar}) holds for $l+1$.
To this end, we carry out the
following backward procedure, or \textquotedblleft backward
subroutine\textquotedblright, at level $l+1$:
\begin{itemize}
\item As initialization, set
$\overline{Y}_{T}^{l+1}=H_{T},$ $\overline{c}_{T}^{l+1}=0$.
(We also set $\overline{c}^0_T = 0$.) %[***REDUNDANT??? MENTIONED ALREADY IN STEP 1???***]
\item Suppose that, for $0<j+1\leq T$,
the values $\overline{Y}_{j+1}^{l+1,n}$, $n=1,\ldots,N$,
the set of $\rho$-martingale increments $\left(\overline{M}_{r}^{l+1}-\overline{M}_{j+1}^{l+1}\right)_{j+1<r\leq T}$
(which is empty if $j+1=T$), and the continuation function $\overline{c}_{j+1}^{l+1}$ have been
constructed.
\item Then construct, according to the regression subroutine in Section~\ref{sec:PD} below,
a continuation function $\overline{c}_{j}^{l+1}$,
a $\rho$-martingale increment $\overline{\mathfrak{m}}_{j+1}^{l+1}$ $\in\mathcal{F}_{j+1}$
with $\rho_{j}(\overline{\mathfrak{m}}_{j+1}^{l+1})=0$,
and set $\left(\overline{M}_{r}^{l+1}-\overline{M}_{j}^{l+1}\right)
%_{j<r\leq T}
=\left(\overline{M}_{r}^{l+1}-\overline{M}_{j+1}^{l+1}
+\overline{\mathfrak{m}}_{j+1}^{l+1}\right)$,
%_{j<r\leq T}
for
\begin{align}
\overline{Y}_{j}^{l+1,n}  &  =\max\left[  \overline{U}_{j}^{l+1,n},
\overline{c}_{j}^{l+1}(X_j^n) \right], \label{algorecursion} \text{ \ \ with}\\
\overline{U}_{j}^{l+1,n}  &  =f_{j}(X_{j}^{n})
+\overline{c}_{j}^{l}(X_{j}^{n}),\text{ \ \ and \ \ }n=1,\ldots,N. \nonumber
\end{align}
\end{itemize}
\end{itemize}
%Note that (\ref{algorecursion}) is equivalent to $ \overline{Y}_{r}^{l+1,n} = \max_{j\leq r \leq t} \left[ \overline{U}_{r}^{l+1,n} - \bar{M}_{r}^{l+1}\left( \omega_{n}\right) \right]$.

Proceeding this way:
\begin{itemize}
\item[(a.)] Working forward from $l=0,\ldots,L$ thus yields
a family of continuation functions $\overline{c}^{l}$ and
a family of (true) $\rho$-martingales $\overline{M}^{l}$,
respectively:
\begin{equation*}
\overline{c}_{j}^{l}(\cdot),\text{ \ \ and \ \ }
\overline{M}_{j}^{l}:=\sum_{r=1}^{j}\overline{\mathfrak{m}}_{r}^{l},\text{ \ \ \ }l=1,\ldots
,L,\text{ }j=0,\ldots,T.
\end{equation*}

\item[(b.)] An upper bound for the solution to the robust multiple stopping problem
at $t=0$ due to $L$ exercise rights
is now given by (cf. Theorem~\ref{th:mulsub}, (ii)):
\begin{equation}
Y_{0}^{\text{upp},L}:=\rho_{0}\left(  \max_{0\leq j_{1}<j_{2}<\cdots< j_{L}}%
\sum_{l=1}^{L}\left(  f_{j_{l}}(X_{j_{l}})-\overline{M}_{j_{l}}^{L-l+1}%
+\overline{M}_{j_{l-1}}^{L-l+1}\right)  \right),
\label{eq:upp}
\end{equation}
which needs to be estimated by a separate (Monte Carlo) procedure.

\item[(c.)] A lower bound for the solution to the robust multiple stopping problem
at $t=0$ due to $L$ exercise rights
may next be obtained from the family of stopping times
\begin{equation}
\label{taul}\tau^{l}:=\min\left\{  j:\tau^{l-1}<j\leq T,\text{ }f_{j}%
(X_{j})+\overline{c}_{j}^{l-1}(X_{j})\geq\overline{c}_{j}^{l}(X_{j})\right\},
\end{equation}
via a (Monte Carlo) estimation of:
\begin{equation}
\label{eq:low}
Y_{0}^{\text{low},L}:=\rho_{0}\left(  \sum_{l=1}^{L}\text{ }f_{\tau^{l}}%
(X_{\tau^{l}})\right).
\end{equation}
\end{itemize}

\subsection{Regression Subroutine}

\label{sec:PD}

Let there be given a collection of
`elementary'
$\rho$-martingale increments
$\mathcal{E}_{j}^\beta(X)$,
with $\beta=(\beta_1,\ldots,\beta_{K^{\prime}}) \in \mathbb{R}^{K^{\prime}}$ and ${K^{\prime}}\in \mathbb{N}$,
that is,
\begin{equation*}
\rho_{j}\left(  \mathcal{E}_{j}^\beta(X) \right)=0,\qquad j=0,\ldots,T-1,
\end{equation*}
and a collection of basis functions
$\psi_{1},\ldots,\psi_{K^{\prime\prime}}:\mathbb{R}^{d}\rightarrow\mathbb{R}$.
We assume that the set of $\rho$-martingale increments
\begin{equation*}
\bigg\{ \mathcal{E}^{(\beta_1,\ldots,\beta_{K^{\prime}})}_j(X)  \bigg\}, %\qquad K^{\prime} \in \mathbb{N},
\end{equation*}
is $L^2$-dense among the $\mathcal{F}_{j+1}$-measurable square-integrable random variables $\mathcal{E}_j$
such that $\rho_{j}(\mathcal{E}_j)=0$.
We then solve, in view of Lemma~\ref{var+} and Corollary~\ref{loc},
for fixed $N$ and $K^{\prime}$ and $K^{\prime\prime}$ the least squares problem
\begin{align}
\text{MSE} &:= \sum_{n=1}^{N}\left(  Y_{j+1}^{l+1,n}-\mathcal{E}_{j}^\beta(X)-\sum_{k=1}^{K^{\prime\prime}}\gamma_{k}\psi_{k}(X_{j}^{n})\right)
^{2} = \sum_{n=1}^{N} \left(  Y_{j+1}^{l+1,n}-\mathcal{E}%
_{j}^\beta(X)-\gamma\psi(X_{j}^{n})\right)
^{2}\nonumber\\
&\longrightarrow \argmin_{\beta\in\mathbb{R}^{K^{\prime}},\gamma\in\mathbb{R}^{K^{\prime\prime}}}=:\left[  \beta
^{l+1,j,K^{\prime},N},\gamma^{l+1,j,K^{\prime\prime},N}\right], \label{argmin}
\end{align}
where we used vector notation $\gamma=(\gamma_1,\ldots,\gamma_{K^{\prime\prime}})$
and $\psi=(\psi_1,\ldots,\psi_{K^{\prime\prime}})^{\transpose}$.
For the algorithm to converge, it is actually sufficient that the MSE above
converges to zero as $K^{\prime},K^{\prime\prime} \rightarrow \infty$ for our choice of $\beta$ and $\gamma$.
We will suppress the superscripts $K^{\prime},K^{\prime\prime}$ and $N$ whenever there is no ambiguity.
We set
\begin{equation*}
\overline{\mathfrak{m}}_{j+1}^{l+1}:=\overline{\mathfrak{m}}_{j+1}^{l+1}(X):=%\sum_{k=1}^{K}\beta_{k}^{l+1,r,K^{\prime},N}%
\mathcal{E}_{j}^{\beta^{l+1,j}}(X),\text{ \ \ \ }\overline{c}_{j}^{l+1}%
(\cdot):=\sum_{k=1}^{K^{\prime\prime}}\gamma_{k}^{l+1,j}\psi_{k}(\cdot).
\end{equation*}

%\subsection*{Discussion}

%The trust of Corollary \ref{loc} is that the pseudo algorithm sketched above
%is locally consistent in the following sense. If $U_{r}^{l+1}=U_{r}^{\ast
%l+1}$ (i.e., $\overline{c}_{r}^{l}=c_{r}^{\ast l}$), and $\overline{Y}%
%_{r+1}^{l+1}=Y_{r+1}^{\ast l+1},$ then $\overline{Y}_{r+1}^{l+1}%
%-\mathfrak{m}_{r+1}^{l,N}$ converges to $\rho_{r}\left(  Y_{r+1}^{\ast
%l}\right)  $ in $L_{2},$ if $K\rightarrow\infty$ and $N_{K}\rightarrow\infty.$

%Concerning the estimations $\overline{c}_{r}^{l+1},$ well-known distribution
%free results connected with (nonparametric) regression, see for example
%\cite{GKKW02}, imply the following result in our context. . . %. . . .

%\bigskip

%UNDER\ CONSTRUCTION

%Set $Y_t^{l,K,N} := Y_{t_i}^{l,K,N}$ for $t_i \leq t \leq t_{i+1}$ and define $M_t^{l,K,N}$ and $c_t^{l,K,N}$ similarly. Set $\mathfrak{m}^{l,K,N}_{t_i} = \beta^{l,K,N}\mathcal{E}_{i+1}(X_{t_{i+1}})+\gamma^{l,K,N}\psi(X_{t_i})$ where $X$ should be sampled independent of $X^1,\ldots,X^N$.

\subsection{Convergence Theorem}

We state the following theorem:

\begin{theorem}
\label{theoConv}
Let $\rho$ be subadditive (P1)
and Lipschitz continuous in the sense of (\ref{l1}) for $p=2$.
We set $K=\min{\{K^{\prime},K^{\prime\prime}\}}$ and denote by
$\overline{M}_{j}^{l,K,N}:=\overline{M}_{j}^{l,K,N}(X_j):=\overline{M}_{j}^{l}$,
$\overline{c}_{j}^{l,K,N}:=\overline{c}_{j}^{l,K,N}(X_j):=\overline{c}_{j}^{l}$ and
$\overline{Y}_{j}^{l,K,N}:=\overline{Y}_{j}^{l,K,N}(X_j):=\overline{Y}_{j}^{l,N}$
the functions constructed in the algorithm above. %[***WHERE: BE MORE EXPLICIT.***]
Then,
\begin{equation}
\lim_{K\rightarrow \infty} \lim_{N\rightarrow \infty} \overline{M}_j^{l,K,N} = M_j^{\ast, l} \mbox{ in } L^2, \label{Mconv}
\end{equation}
\begin{equation}
\lim_{K\rightarrow \infty} \lim_{N\rightarrow \infty} \overline{c}_j^{l,K,N} = c_j^{\ast, l} \mbox{ in } L^2, \label{Cconv}
\end{equation}
\begin{equation}
\lim_{K\rightarrow \infty} \lim_{N\rightarrow \infty} \overline{Y}_j^{l,K,N} = Y_j^{\ast, l} \mbox{ in } L^2, \label{Yconv}
\end{equation}
for all $j= T, T-1, \ldots, 0$ and $l= 1,\ldots,L$.
Furthermore,
\begin{align*}
\left\vert Y_{j}^{\ast, l} \right. \left. -\rho_{j}\left(  \max_{j\leq r\leq T}\left(
\overline{U}_{r}^{l}-\overline{M}_{r}^{l}\right)  \right)  \right\vert
&\leq %C_{\mathrm{Lipschitz}}
\rho_{j}\left(  \max_{j\leq r\leq T}\left\vert M_{r}^{\ast, l}-\overline{M}_{r}^{l}\right\vert \right)
+ %C_{\mathrm{Lipschitz}}
\rho_{j}\left(  \max_{j\leq r\leq T}\,\left\vert \overline{c}_{r}^{l-1}-c_{r}^{\ast, l-1}\right\vert \right) \\
&\longrightarrow_{K\rightarrow \infty, N\rightarrow \infty} 0.
\end{align*}
\end{theorem}

%\begin{remark}
%We note that, now that Theorem \ref{theoConv} has been established,
%the convergence in probability of optimal stopping times
%follows from applying similar arguments as in Kr\"atschmer \textit{et al.} \cite{KLLSS18},
%Theorem 18.
%\end{remark}

\subsection{Complexity}

At first sight, the path-wise maximum in \eqref{eq:upp} would require the evaluation
of $T!/(L!(T-L)!)$ terms.
Fortunately, due to following proposition, it only requires $O(LT)$ evaluations.

\begin{proposition}
Define, for $1\leq q\leq L$ and $0\leq i\leq T,$%
\begin{equation*}
\Theta_{i}^{q}:=\max_{i\leq j_{1}<j_{2}<\cdots< j_{q}}\sum_{l=1}^{q}\left(
f_{j_{l}}(X_{j_{l}})-\overline{M}_{j_{l}}^{q-l+1}+\overline{M}_{j_{l-1}%
}^{q-l+1}\right)  , \text{ \ \ with \ }j_{0}=i,
\end{equation*}
and naturally $\Theta_{i}^{q}=0$, $i>T$.
Then,
\begin{equation}
\Theta_{i}^{q}=\max\left[  f_{i}(X_{i})+\overline{M}_{i}^{q-1}-\overline
{M}_{i+1}^{q-1}+\Theta_{i+1}^{q-1},\overline{M}_{i}^{q}-\overline{M}_{i+1}%
^{q}+\Theta_{i+1}^{q}\right]. \label{rec}
\end{equation}
\label{prop:com}
\end{proposition}

Thus, the evaluation of \eqref{eq:upp} may be described as follows:

\paragraph{Recursive evaluation of \eqref{eq:upp}}

\begin{itemize}
\item[1.] Initialize $\Theta_{i}^{0}=0,$ for $i=0,\ldots,T;$

\item[2.] Suppose that, for $0\le q-1<L$ and for $i=0,\ldots,T$, the construction
of $\Theta_{i}^{q-1}$ has been conducted;

\item[3.] Backward subroutine: Initialize $\Theta_{T}^{q}=f_{T}(X_{T})$.
When $\Theta_{i+1}^{q}$ has been constructed for $i+1\leq T$,
compute $\Theta_{i}^{q}$ via (\ref{rec}).
\end{itemize}

\setcounter{equation}{0}

\section{Explicit %Numerical
Construction %of the Upper Bound and Lower Bound
in a Brownian-Poisson Filtration}

\label{sec:BPcon}

In the sequel, we assume that we have a completed continuous-time filtration $\mathbb{F}=(\mathcal{F}_t)_{t\in[0,T]}$
on a filtered probability space $(\Omega, \mathcal{F}, \mathbb{F}, \mathbb{P})$
generated by a $d_1$-dimensional standard (i.e., zero mean and unit variance)
Brownian motion $W=(W_1,\ldots,W_{d_1})^{\transpose}$
and a $d_2$-dimensional Poisson process $N=(N_1,\ldots,N_{d_2})^{\transpose}$
with arrival intensity $\lambda_{\mathbb{P}}=(\lambda_{\mathbb{P}}^{1},\ldots,\lambda_{\mathbb{P}}^{d_{2}})^{\transpose}$.
As usual, we define the compensated counterpart of $N$ as $\tilde{N}_{t}=N_{t}-\lambda_{\mathbb{P}}t$.
%The processes $W$ and $N$ are assumed to mutually independent,
%as are all their components.
The components of the processes $W$ and $N$ are assumed to be independent.
The stochastic drivers $W$ and $N$ generate
the underlying Markovian adapted reward process $(X_{t})_{t\in[0,T]}$ with state space $\mathbb{R}^{d}$
of Section~\ref{sec:MSAL}.

Furthermore, we assume that $\rho$ satisfies (the continuous-time analogs of) (C1)--(C4) and (P1)--(P3).
This means, in particular, that $\rho$ is a coherent risk measure.
By classical duality results (e.g., F\"ollmer and Schied \cite{FS11}),
the robust multiple stopping problem at time $t$ is then given by
\begin{align}
\label{eq:prob}
Y_{t}^{\ast,L}= \sup_{\substack{t\leq\tau_{1}<\cdots<\tau_{L} \\
(\tau_{1},\ldots,\tau_{L})\in\mathcal{T}_{t}(L)}}
\rho_{t}\left(  \sum_{l=1}^{L}H_{\tau_{l}}\right)
=\sup_{\substack{t\leq\tau_{1}<\cdots<\tau_{L} \\
(\tau_{1},\ldots,\tau_{L})\in\mathcal{T}_{t}(L)}} \sup_{\mathbb{Q}\in\mathcal{Q}}
\mathbb{E}_{\mathbb{Q}}\left[\sum_{l=1}^{L}H_{\tau_{l}}\Big|\mathcal{F}_{t}\right],
\quad 0\leq t\leq T,\end{align}
with $\mathcal{T}_{t}(L)$ our family of stopping vectors,
and $\mathcal{Q}$ a closed convex set of probability measures absolutely continuous with respect to $\mathbb{P}$
and satisfying a stability assumption.
In such a continuous-time setting, it is known that every recursive coherent risk measure
can be identified with a solution to a backwards stochastic differential equation (BSDE) also called a $g$-expectation,
modulo a compactness assumption;
see Section \ref{sec:tcdp} for the precise definitions and results.
Exploiting our algorithm presented in Section \ref{sec:MSAL},
this section constructs explicit upper and lower bounds to $Y^{\ast,L}$
with desirable properties.

\subsection{Bellman's Principle, the Set of Priors, and BSDE drivers}

\label{sec:tcdp}

A probability measure change from $\mathbb{P}$ to an absolutely continuous measure $\mathbb{Q}\in\mathcal{Q}$
admits an explicit representation in our Brownian-Poisson setting.
Consider the Radon-Nikodym derivative
\begin{equation*}
D_{t}:=\mathbb{E}\left[\frac{d\mathbb{Q}}{d\mathbb{P}}|\mathcal{F}_{t}\right],\qquad t\in[0,T].
\end{equation*}
As is well-known, $D_{t}$ has the Dol\'eans-Dade exponential form
\begin{equation}
\label{eq:Radon}
D_{t}=\exp\left(\int_{0}^{t}q_{s}dW_{s}
+ \int_{0}^{t}\log\left(\frac{\lambda_{s}}{\lambda_{\mathbb{P}}}\right)dN_{s}
-\int_{0}^{t}\left(\frac{|q_{s}|^{2}}{2} + \lambda_{s}-\lambda_{\mathbb{P}}\right)ds\right),
\end{equation}
where $q$ is a predictable, $\mathbb{R}^{d_{1}}$-valued, stochastic drift,
and $\lambda$ is a positive, predictable, $\mathbb{R}^{d_{2}}$-valued process
with $\frac{\lambda_{s}}{\lambda_{\mathbb{P}}}:=(\frac{\lambda^{1}_{s}}{\lambda^{1}_{\mathbb{P}}},
\ldots,\frac{\lambda^{d_{2}}_{s}}{\lambda^{d_{2}}_{\mathbb{P}}})^{\transpose}$,
which jointly uniquely characterize $\mathbb{Q}$.
From Girsanov's theorem, we know that
$W^{\mathbb{Q}}_{t}:=W_{t}-\int_{0}^{t} q_{s} ds$ is a standard Brownian motion under $\mathbb{Q}$
while the process ${N}_{t}$ has arrival intensity $\lambda_{t}$.
In particular, the reference model $\mathbb{P}$ corresponds to $q\equiv 0$ and $\lambda\equiv \lambda_{\mathbb{P}}$.
The stochastic drift $q$ may be given the interpretation of a drift in the diffusive component
that is misspecified to be absent by the reference model $\mathbb{P}$.
Similarly, $\lambda_{s}-\lambda_{\mathbb{P}}$ represents a deviation
from the misspecified arrival intensity $\lambda_{\mathbb{P}}$ under $\mathbb{P}$.

In our dynamic setting,
with DMU evaluations satisfying the continuous-time analogs of (C1)--(C4) and (P1)--(P3),
time-consistency of choice under uncertainty is satisfied
as it is equivalent to recursiveness
or Bellman's dynamic programming principle.
Time-consistency of a dynamic evaluation $(\rho_{t}(H))_{t\in[0,T]}$ requires---according to its usual definition,
also referred to as `strong' time-consistency---that $\rho_{s}(H_{1}) \ge\rho_{s}(H_{2})$
whenever $\rho_{t}(H_{1}) \ge\rho_{t}(H_{2})$, $t\geq s$.
That is, if $H_{2}$ is preferred over $H_{1}$,
in each state of nature at time $t$,
then the same preference necessarily applies prior to time $t$;
see e.g., %Duffie and Epstein \cite{DE92}, Chen and Epstein \cite{CE02},
Riedel \cite{R04}, Ruszczy\'nski and Shapiro \cite{RS06b},
Shapiro, Dentcheva and Ruszczy\'nski \cite{SDR09}, Chapter 6,
Ruszczy\'nski \cite{R10} and Shapiro \cite{S16}.
Indeed, requiring recursiveness or Bellman's dynamic programming principle
is equivalent to requiring time-consistency for
$\rho_{t}(H)=\sup_{\{\mathbb{Q}\sim \mathbb{P}|\mathbb{Q}=\mathbb{P}\mbox{ on }\mathcal{F}_{t}\}}
\mathbb{E}_{\mathbb{Q}}[H|\mathcal{F}_{t}]$, $t\in[0,T]$,
which is, in turn, equivalent to the set of priors $\mathcal{Q}$ being $m$-stable; see Delbaen \cite{D06}.
More formally, the following statements are equivalent
(see Lemma 11.11 of F\"ollmer and Schied \cite{FS11} for the equivalence (i)--(ii),
Delbaen \cite{D06} and Delbaen, Peng and Rosazza Gianin \cite{DPR10} for (ii)--(iii) in a Brownian setting,
and Tang and Wei \cite{TW12} and Laeven and Stadje \cite{LS14} for (ii)--(iii) in a general semi-martingale setting):
%[***RESTRICTION TO BOUNDED REWARDS***]
\begin{itemize}
\item[(i)] $\rho$ is recursive, i.e., $\rho$ satisfies Bellman's dynamic
programming principle $\rho_{0}(\rho_{t}(H)I_{A})=\rho_{0}(H I_{A})$
for every $t\in[0,T]$, $A\in\mathcal{F}_{t}$, and bounded $H$.

\item[(ii)] $\rho$ is time-consistent over bounded rewards.

\item[(iii)] There exists a closed, convex, set-valued predictable mapping $C$
taking values in $\mathbb{R}^{d_{1}}\times(-\lambda_{\mathbb{P}}^{1},\infty
)\times\cdots\times(-\lambda_{\mathbb{P}}^{d_{2}},\infty)$
such that %[***VERIFY***]
\begin{equation*}
\rho_{t}(H)=\sup_{(q,\lambda)\in C}\mathbb{E}_{\mathbb{Q}}\left[H|\mathcal{F}_{t}\right],\qquad t\in[0,T].
\end{equation*}%$r(s,q,v)=I_{C_{s}}(q,v).$
\end{itemize}

\vskip -0.2cm As in our continuous-time Brownian-Poisson setting recursiveness (C3) is equivalent to (iii),
we assume henceforth:
\begin{itemize}
\item[(A1)] $C=(C_{t})_{t\in[0,T]}\subset[0,T]\times\mathbb{R}^{d_{1}}
\times(-\lambda_{\mathbb{P}}^{1}+\varepsilon,\infty)\times\cdots\times(-\lambda_{\mathbb{P}}^{d_{2}}+\varepsilon,\infty)$
with $\varepsilon> 0$ is compact.
%$C=(C_s)_{s\in[0,T]}\subset [0,T]\times
%\mathbb{R}^d\times (-\lambda_{P}^{1},\infty)\times\ldots\times(-\lambda_{P}^{k},\infty)$.
\end{itemize}
We note that (A1) also implies that $\rho$ is recursive and time-consistent over square-integrable rewards.

%convex/concave $\mathbb{P}\otimes dt-$a.s. (cf. Jiang \cite{J08}).
For $t\in[0,T]$, $z\in\mathbb{R}^{1\times d_{1}}$ and $\tilde{z}\in\mathbb{R}^{1\times d_{2}}$,
and $C$ satisfying Assumption (A1),
let us define a function $g$ via Fenchel's duality:
\begin{equation}
\label{g1}
g(t,z,\tilde{z}):=\sup_{(q,\lambda-{\lambda}_{\mathbb{P}})\in C_{t}}
\{zq+\tilde{z} (\lambda-\lambda_{\mathbb{P}})\}.
\end{equation}
One easily verifies that $g$ is convex, positively homogeneous and Lipschitz continuous.
%Suppose now that, for every exercise date $r$, $r=0,\ldots,L-1$, we have a
%finer time grid $\pi_{r}=\{s_{r0}=r<s_{r1}<\ldots<s_{rP}=r+1\}$.
%Denote the corresponding overall time grid by $\pi=\{s_{00},s_{01}%
%,\ldots,s_{LP}\}$. The following theorem provides a way to practically approximate
%$M^{\ast}$ in Theorem \ref{mulsub} by connecting it to specific
%semi-martingale dynamics that can be dealt with numerically in an efficient
%way.
%and is also known in the literature
%under the term backward stochastic differential equation
%(BSDE).
%which will turn out to be very useful when constructing numerically dual Doob-martingales for the stopping problem (\ref{OST}).
Then, from Kr\"atschmer \textit{et al.} \cite{KLLSS18}, we have the following statement,
which is essentially (with (A1)) equivalent to (i)--(iii) above:
%\begin{theorem}
%\label{theobsde}
\begin{itemize}
\item[(iv)] For every $H\in L^{2}(\mathcal{F}_{j+1})$,
there exists a unique square-integrable predictable $(Z,\tilde{Z})$ such that
\begin{align*}
d\rho_{t}(H_{j+1})=-g(t,Z_{t},\tilde{Z}_{t})dt
+Z_{t}dW_{t} + \tilde{Z}_{t}d\tilde{N}_{t},\ \mbox{for
		}\,t\in[j,j+1],\ j\in\{0,\ldots,T-1\}.
\end{align*}
\end{itemize}
In particular, there exists a unique square-integrable predictable
$(Z^{\ast},\tilde{Z}^{\ast})$ such that
\begin{equation}
\label{eq:BSDE2}
d\rho_{t}(Y^{\ast}_{j+1})=-g(t,Z^{\ast}_{t},\tilde{Z}^{\ast}_{t})dt
+Z^{\ast}_{t}dW_{t} + \tilde{Z}^{\ast}_{t}d\tilde{N}_{t},\ \mbox{for
        }\,t\in[j,j+1],\ j\in\{0,\ldots,T-1\}.
\end{equation}
Furthermore, for $t\in[0,T]$, the $(Z^{\ast},\tilde{Z}^{\ast})$ in \eqref{eq:BSDE2}
recover---and later allow to practically compute---the $\rho$-Doob martingale as follows
(cf. Eqns.~\eqref{eq:DoobMeyerdisc}--\eqref{eq:DD}):
\begin{equation}
\label{sttar}M^{\ast}_{t} =\rho_{t}(M^{\ast}_{T})=-\int_{0}^{t} g(s,Z^{\ast}_{s},
\tilde{Z}^{\ast}_{s})ds + \int_{0}^{t} Z^{\ast}_{s}\,dW_{s}+ \int_{0}^{t}\tilde
{Z}^{\ast}_{s}\,d\tilde{N}_{s}.
\end{equation}
Here, $j=0,\ldots,L-1$ should be interpreted as exercise dates
and $t\in[j,j+1]$ as the continuous embedding.
%\end{theorem}

%\begin{remark}
Because $\rho_{j+1}(Y^{\ast}_{j+1})=Y^{\ast}_{j+1}$, by (iv), for $t\in[j,j+1]$,
\begin{equation}
\label{ruck}
\rho_{t}(Y^{\ast}_{j+1})=Y^{\ast}_{j+1}+\int_{t}^{j+1}g(s,Z^{\ast}_{s},\tilde{Z}^{\ast}_{s})ds
-\int_{t}^{j+1}Z^{\ast}_{s}dW_{s}
-\int_{t}^{j+1}\tilde{Z}^{\ast}_{s}d\tilde{N}_{s}.
\end{equation}
%Similarly,
%	it follows that, for $t\in(j,j+1]$, \be \label{ruck2}
%	\rho^h_{t}=\rho^h_{j+1}+\int_t^{j+1}
%	g(s,Z^h_{s},\tilde{Z}^h_s)ds-\int_t^{j+1}Z^h _{s}dW_{s} -
%	\int_t^{j+1}\tilde{Z}^h_{s}d\tilde{N}_{s}.\ee
%\end{remark}
Eqn. (\ref{eq:BSDE2}) is referred to as a backward stochastic differential equation (BSDE).
Formally, given a terminal payoff $H\in L^{2}$
and a function $g:[0,T]\times\mathbb{R}^{d_{1}}\times\mathbb{R}^{d_{2}}\to\mathbb{R}$,
referred to as a driver, the solution to the corresponding BSDE is a
triple of square-integrable and suitably measurable processes $(Y,Z,\tilde{Z})$
that satisfies
\begin{equation*}
dY_{t}=-g(t,Z_{t},\tilde{Z}_{t})dt+Z_{t} dW_{t}+\tilde{Z}_{t}d\tilde{N}%
_{t},\quad\mbox{and}\quad Y_{T}=H.
\end{equation*}
The solution is often referred to as a (conditional) $g$-expectation;
see, e.g., Peng \cite{P04}.

%With $Z^*$ and $\tilde{Z}^*$ at hand we can compute $M^{*,U^h}$
%and with (\ref{kez2}) also $U^h_{j}(V^{\ast}_{t_{j + 1}}).$

As a means of illustrating the generality of our setup given by (\ref{eq:prob}) with (A1)
we provide a few examples.

\begin{example}
\label{g-ex}
\begin{itemize}
\item[(1.)] Ball scenarios:
Consider a decision-maker endowed with a set of priors
constituting a small ball environment surrounding $\mathbb{P}$
all deemed equally plausible.
Then,
\begin{align*}
\mathcal{Q}=\Big\{\mathbb{Q}^{(q,\lambda)} \ll \mathbb{P} \Big||q_{t}|\leq\delta_{1},\quad|\lambda_{t}-{\lambda}_{\mathbb{P}}|\leq\delta_{2},
\mbox{ for Lebesgue-a.s. all }t\Big\},\qquad \delta_{1},\delta_{2}>0,
\end{align*}
and
%\begin{align*}
$C_{t}= \big\{(q,\lambda) \big{|} |q|\leq\delta_{1}, |\lambda-{\lambda}_{\mathbb{P}}|\leq\delta_{2}\big\}$.
%\end{align*}
Suppose without losing generality that $|{\lambda}_{\mathbb{P}}|\geq\delta_{2}$.
Then, from \eqref{g1}, in explicit form,
%\begin{equation*}
$g(t,z,\tilde{z})=\delta_{1}|z| + \delta_{2}|\tilde{z}|$.
%\end{equation*}
\item[(2.)] Discrete scenarios:
Imagine a decision-maker who considers,
at each time $t>0$,
finite-dimensional families $\{q_{1,t},\ldots,q_{m,t}\}$ and $\{\lambda_{1,t},\ldots,\lambda_{m,t}\}$,
$m\in\mathbb{N}$, with all elements deemed equally plausible.
Then,
\begin{align*}
\mathcal{Q}=\Big\{\mathbb{Q}^{(q,\lambda)} \ll \mathbb{P} \Big|(q_{t},\lambda_{t})\in\{(q_{i,t},\lambda_{j,t}), i,j\in\{1,\ldots,m\}\},
\mbox{ for Lebesgue-a.s. all }t\Big\},
\end{align*}
and
$C_{t}= \big\{(q,\lambda) \big{|} (q,\lambda)\in\mathrm{conv}\left(\{(q_{i,t},\lambda_{j,t}), i,j\in\{1,\ldots,m\}\}\right)\big\}$,
with $\mathrm{conv}(\cdot)$ the convex hull.
We can assume that $0\in\mathrm{conv}\left(\{(q_{i,t},\lambda_{j,t}), i,j\in\{1,\ldots,m\}\}\right)$
without losing generality, upon redefining the reference measure.
Furthermore,
%\begin{equation*}
$g(t,z,\tilde{z})=\max_{i=1,\ldots,m}q_{i,t}z+\max_{j=1,\ldots,m}\lambda_{j,t}\tilde{z}$.
%\end{equation*}
\end{itemize}
\end{example}

To obtain genuine upper and lower bounds to the optimal solution of the stopping problem \eqref{eq:prob},
we henceforth impose the following additional assumption:

\begin{itemize}
\item[(A2)] $H_{t} = f_{t}(X_{t})$ and we can simulate i.i.d. copies of
$(X_{t})_{t\in[0,T]}$.

%and that either (i) $X$ is a
%submartingale, or (ii) $X$ has independent
%increments, or (iii) $\log(X)$ has independent
%increments.
%\item[(b)] $f=0$ and $r$ is independent of $x.$
%\end{itemize}

\end{itemize}

%Assumption (G2) is needed (only) to obtaine genuine upper and lower bounds.
%(An Euler scheme, for example, may induce a bias.) Assumption (G2) is
%satisfied in particular if $X$ follows a linear SDE, which holds e.g., in the
%case of a Brownian motion with drift, a Poisson process with drift, an
%Ornstein-Uhlenbeck process, or a geometric Brownian motion with drift (with or
%without Poisson type jumps). But note there are by now also very general
%results available on exact sampling of more general diffusions and
%jump-diffusions; see, e.g., Beskos and Roberts \cite{BR05}, Broadie and Kaya
%\cite{BK06}, Chen and Huang \cite{CH13}, Giesecke and Smelov \cite{GS13}, and
%Henry-Labord\`{e}re, Tan and Touzi \cite{HTT15}.

%Proposition \ref{ADH} gives a way how to obtain an upper bound

%Let $g$ be any driver of a BSDE satisfying the standard assumptions (A1)--(A2),
%and let $\rho=(\rho)_{t\in[0,T]}$ denote the continuous-time coherent risk measure
%according to the associated family of conditional $g$-expectations.
%Furthermore, consider the $\rho$-Snell envelope $Y^{\ast}$
%as defined before.
%Our aim is to find upper and lower bounds for $Y_{0}^{\ast}$.

Once, we have constructed a `good' family of $\rho$-martingales $\overline{M}$,
we are faced with the computation of $Y_{0}^{\text{upp},L}$ in \eqref{eq:upp}.
%\begin{equation*}
%\label{gen1}Y_{0}^{\text{upper}}:=\rho_{0}\left(  \max_{0\leq j_{1}%
%<j_{2}<\cdots< j_{L}}\sum_{l=1}^{L}\left(  f_{j_{l}}(X_{j_{l}})-\overline
%{M}_{j_{l}}^{L-l+1}+\overline{M}_{j_{l-1}}^{L-l+1}\right)  \right)  ,
%\end{equation*}
%according to Proposition~\ref{prop:ADP}. [***ADAPT: (Theorem \ref{th:mulsub}, (ii))***]
An important advantage of this dual approach for numerical stability is that
%for the optimal martingale
we have, in fact, the pathwise dual representation
\begin{equation*}
Y_{0}^{\ast,L}:= \max_{0\leq j_{1}<j_{2}<\cdots< j_{L}}\sum_{l=1}^{L}\left(
f_{j_{l}}(X_{j_{l}})-M_{j_{l}}^{\ast,L-l+1}+M_{j_{l-1}}^{\ast,L-l+1}\right)
,\qquad\text{almost surely}.
\end{equation*}
%For example, in the case where $\rho_{0}$ stands for the standard expectation
Indeed, we will obtain an estimate with `low' variance
provided the $\rho$-martingale $\overline{M}$ is `good'.
Furthermore, to obtain a lower bound we employ $Y_{0}^{\text{low},L}$ in \eqref{eq:low}.
%\begin{align}
%\label{primal}Y_{0}^{\text{low}}:=\rho_{0}\left(  \sum_{l=1}^{L}\text{
%}f_{\tau^{l}}(X_{\tau^{l}}) \right)  .
%\end{align}
By the results discussed in this subsection, for a square-integrable payoff $U$, $\rho$ has a representation of the form
\begin{align}
\rho_{t}(U)  &  = \sup_{(q,\lambda) \in C} \mathbb{E}_{\mathbb{Q}^{(q,\lambda)}}\left[
U \big|\mathcal{F}_{t}\right] \label{bs2}\\
&  =U+\int_{t}^{T}g(s,\mathcal{Z}_{s},\mathcal{\tilde{Z}}_{s})ds-\int_{t}%
^{T}\mathcal{Z}_{s}dW_{s} - \int_{t}^{T}\mathcal{\tilde{Z}}_{s}d\tilde{N}_{s}.
\label{bs}%
\end{align}
Let us now first consider the question of how to explicitly obtain
a `good' family of $\rho$-martingales $\overline{M}$
in our Brownian-Poisson filtration with (C1)--(C4), (P1)--(P3) and (A1)--(A2).
Next, we develop in Sections~\ref{sec:LB} and~\ref{sec:UB} explicit genuine %two approaches for the computation of
lower and upper bounds. % to \eqref{eq:prob}.

\subsection{Parameterization of the \texorpdfstring{$\rho$-}{rho-}Martingale Increments} %in a Brownian-Poisson Filtration}

\label{sec:algorithm}

As in Section \ref{sec:MSAL},
we assume that $(X_{j})_{0\leq j\leq T}$ is an $\mathcal{F}_{j}$-adapted
$d$-dimensional underlying Markovian process, now in a Brownian-Poisson filtration,
and that our cash-flow process has a structure of the form
$H_{j}=f_{j}(X_{j})$, $j=0,\ldots,T$ such that $H\in\mathcal{H}$.
We further assume that the resulting random variables $H_{j}$ are square integrable.

%In view of Theorem \ref{th:mulsub},
We are going to construct a $\rho$-martingale backwardly.
To this end, we consider, between two exercise dates $j$ and $j+1$,
the (fine) grid $\pi_{j}=\{s_{(j-1)n_{0}}=j,s_{j1},\ldots,s_{jn_{0}}=j+1\}$,
where $s_{jp}=j+p\Delta$ with $\Delta=n_{0}^{-1}$.
We also define
\begin{equation*}
\Pi:=\{s_{00}=0,s_{01},\ldots,s_{0n_{0}}=1,s_{11},\ldots,%s_{1n_{0}}=2,s_{21},\ldots,
s_{(T-2)n_{0}}=T-1,s_{(T-1)1},\ldots,s_{(T-1)n_{0}}=T\},
\end{equation*}
and sometimes use the notation $\Pi=\{t_{0}=0,t_{1},t_{2},\ldots,t_{n_{1}}=T\}$,
where the $t_{i}$ are simply the enumerated $s_{jp}$.
For our numerical schemes we always assume that Assumption (A2) is in place,
next to (A1).
In particular, we can also simulate i.i.d. copies
of the $(Z_{s_{jp}},\tilde{Z}_{s_{jp}})_{s_{jp}\in\Pi}$.

We formally initialize
$\left(\overline{M}_{j}^{l+1}- \overline{M}_{T}^{l+1}\right)_{j\geq T}=(0)_{T}$.
Suppose that for some (fixed) $j<T$, an approximation $Y_{j+1}^{l+1}$ to the upper
Snell envelope $Y_{j+1}^{\ast,l+1}$ and the set of $\rho$-martingale increments
$\left(\overline{M}_{q}^{l+1}-\overline{M}_{j+1}^{l+1}\right)_{j+1\leq q\leq T}$ have been constructed.
Then, we carry out the following loop.
For $p=n_{0}$, we initialize $0\leq U_{n_{0}}=Y_{j+1}^{l+1}\approx\rho_{j+1}\left(Y_{j+1}^{\ast,l+1}\right)$.
Now, if $0\leq U_{p+1}$, $p<n_{0}$, has been constructed,
we solve the piecewise linear minimization problem
\begin{align}
\label{defag}
%&
\left[  \gamma^{N_1}_{s_{jp}},\beta^{N_1}_{s_{jp}},\tilde{\beta}^{N_1}_{s_{jp}}\right]  %\nonumber\\
=&\underset
{\gamma,\beta,\tilde{\beta}\in\mathbb{R}^{K_1}}{\arg\min}\,\,\frac{1}{N_1}\sum_{n=1}^{N_1}  \bigg(
U^n_{p+1}- \sum_{k=1}^{K_1}\gamma_{k} \psi_{k}(s_{jp},X^n_{s_{jp}})\nonumber\\
\hspace{0.2cm}&  + g\Big(s_{jp},\sum_{k=1}^{K_1}%
\beta_{k}\varphi_{k}(s_{jp},X^n_{s_{jp}}),\sum_{k=1}^{K_1}%
\tilde{\beta}_{k}\tilde{\varphi}_{k}(s_{jp},X^n_{s_{jp}})\Big)(s_{j(p+1)}-s_{jp}) \nonumber\\
\hspace{0.2cm}& - \sum_{k=1}^{K_1}%
\beta_{k}\varphi_{k}(s_{jp},X^n_{s_{jp}})\Delta W^n_{s_{jp}}- \sum_{k=1}^{K_1}%
\tilde{\beta}_{k}\tilde{\varphi}_{k}(s_{jp},X^n_{s_{jp}})\Delta\tilde{N}^n_{s_{jp}}  \bigg)^{2},
\end{align}
%\end{document}
for certain basis functions $(\psi_{k})$, $(\varphi_{k})$, and $(\tilde{\varphi}_{k})$,
$K_{1}\in\mathbb{N}$,
and $N_{1}$ trajectories.
Alternatively to (\ref{defag}), we can solve
\begin{align}
\label{defag2}
%&
\left[  \gamma^{N_1}_{s_{jp}},\beta^{N_1}_{s_{jp}},\tilde{\beta}^{N_1}_{s_{jp}}\right]  %\nonumber\\
=&\underset
{\gamma,\beta,\tilde{\beta}\in\mathbb{R}^{K_1}}{\arg\min}\,\,\frac{1}{N_1}\sum_{n=1}^{N_1}  \bigg(
U^n_{p+1}- \sum_{k=1}^{K_1}\gamma_{k} \psi_{k}(s_{jp},X^n_{s_{jp}}) \nonumber\\
\hspace{0.2cm}&- \sum_{k=1}^{K_1}\beta_{k}\varphi_{k}(s_{jp},X^n_{s_{jp}})\Delta W^n_{s_{jp}}
- \sum_{k=1}^{K_1}\tilde{\beta}_{k}\tilde{\varphi}_{k}(s_{jp},X^n_{s_{jp}})
\Delta \tilde{N}^n_{s_{jp}}\bigg)^{2},
\end{align}
which has a closed-form solution.

\begin{remark}
We note that in the case the filtration is generated by a one-dimensional process
(either a Brownian motion or a Poisson process),
the minimization problem (\ref{defag}) corresponds to a linear programming problem.
This is seen as follows.
As $\rho$ is a coherent risk measure, it follows that $g$ is positively homogeneous.
Thus, there exist $f_{s_{jp}}^+ , f_{s_{jp}}^- \geq 0$
such that $g(s_{jp},z) = f_{s_{jp}}^+ z^+ + f_{s_{jp}}^- z^-$.
Hence, the function
\begin{equation*}
z\mapsto h(s_{j(p+1)},z)^2 := \left( U_{p+1}^n + g(s_{jp},z) (s_{j(p+1)}-s_{jp}) - (W_{s_{j(p+1)}}^n-W_{s_{jp}}^n)z\right)^2
\end{equation*}
is convex as $U_{p+1}^n \geq 0$.
(The reason is that $h(s_{jp},\cdot)$ is linear on its negative part.)
Thus, the minimization problem (\ref{defag}) is convex.
Because any piecewise linear function that is convex can be written as
a supremum of finitely many linear functions,
the minimization problem
can be expressed as a linear programming problem.
%Now it is well-known that a piecewise linear convex optimization problem
%can be transformed into a linear programming problem.
\end{remark}

Next, define
\begin{align*}
C_{p}(X^n_{s_{jp}})  &  :=\sum_{k=1}^{K_1}\gamma^{N_1}_{s_{jp}k}\psi_{k}(s_{jp},X_{s_{jp}}^n),
\quad\mathcal{Z}_{s_{jp}}^{l+1,N_1}(X^n_{s_{jp}})
:=\sum_{k=1}^{K_1}\beta^{N_1}_{s_{jp}k}\varphi_{k}(s_{jp},X_{s_{jp}}^n),
\end{align*}
and $\tilde{\mathcal{Z}}_{s_{jp}}^{l+1,N_1}(X^n_{s_{jp}})$ similarly.
We set
\begin{equation*}
U_{p}(X^n_{s_{jp}}):=\max \Big( C_{p}(X^n_{s_{jp}})+g(s_{jp},\mathcal{Z}^{l+1,N_1}_{s_{jp}}(X^n_{s_{jp}}),\tilde{\mathcal{Z}}^{l+1,N_1}_{s_{jp}}(X^n_{s_{jp}}))(s_{j(p+1)}-s_{jp}),0\Big).
\end{equation*}

%For $ j\leq s_{jp}\leq t<  s_{j(p+1)}$
We then obtain the desired $\rho$-martingale increments
$\overline{M}_{s_{jp}}^{l+1,K_1,N_1}-\overline{M}_{j}^{l+1,K_1,N_1}$
by defining
\begin{align}
&\overline{M}_{s_{jp}}^{l+1,\Delta,K_1,N_1}(X_{s_{jp}}^n)-\overline{M}_{j}^{l+1,\Delta,K_1,N_1}(X_{s_{jp}}^n)\equiv
\overline{M}_{s_{jp}}^{l+1,K_1,N_1}(X_{s_{jp}}^n)-\overline{M}_{j}^{l+1,K_1,N_1}(X_{s_{jp}}^n)  \nonumber\\
&:= -\sum_{u=0}^{p-1}\int_{s_{ju}}^{s_{j(u+1)}}
g(u,\mathcal{Z}_{s_{ju}}^{l+1,N_1}(X_{s_{ju}}^n),\tilde{\mathcal{Z}}_{s_{ju}}^{l+1,N_1}(X_{s_{ju}}^n) )du\nonumber\\
&\qquad+\sum_{u=0}^{p-1}\mathcal{Z}_{s_{ju}}^{l+1,N_1}%
(X_{s_{ju}}^n) \Delta W_{s_{ju}}^n  +\sum_{u=0}^{p-1}\tilde{\mathcal{Z}}_{s_{ju}}^{l+1,N_1}%
(X_{s_{ju}}^n) \Delta \tilde{N}_{s_{ju}}^n . \label{m}%
\end{align}
In the end, when we have arrived at $p=0$,
we define $\overline{c}_{j}^{l+1,K_{1},N_{1}}(\cdot):= C_0(\cdot)$ and
$\overline{Y}_{j}^{l+1,K_{1},N_{1}}(\cdot)$ according to Eqn. (\ref{algorecursion}).
This way, we have recursively constructed the parameterized space of $\rho$-martingale increments.
The following proposition establishes convergence of our construction.
\begin{proposition}
The $\rho$-martingale increments constructed in (\ref{m}) are dense in $L^2$ and, in particular,
\begin{equation*}
\lim_{\Delta\to 0}\lim_{K_1\to\infty}\lim_{N_1\to\infty} \overline{M}_t^{l+1,\Delta,K_1,N_1}
= M_t^{\ast,l+1}.
\label{PropConv}
\end{equation*}
\end{proposition}

We finally note that \eqref{m} gives rise to a true discrete-time $\rho$-martingale
$(\overline{M}^{l+1}_{j})_{j\in \{0,1,2,\ldots,T\}}$.
The thus constructed
$(\overline{M}^{l+1}_{j})_{j\in \{0,1,2,\ldots,T\}}$ will be exploited to
establish an upper bound to the upper Snell envelope via \eqref{eq:upp},
while $(\overline{M}_{s_{jp}}^{l+1})_{j,p}$
(living on the finer grid $\Pi$)
is needed for the numerical approximation.

\subsection{Converging Genuine Lower Bound}

\label{sec:LB}

This subsection develops an explicit lower bound
that converges to the upper Snell envelope asymptotically
and constitutes a genuine (biased low) lower bound at the pre-limiting level.
Consider \eqref{bs} with $U =\sum_{l=1}^{L} f_{\tau^{l}}(X_{\tau^{l}})$
and $\tau^{l}$ constructed by \eqref{taul}.
From e.g., Barrieu and El Karoui \cite{BK09}, we know that the supremum in \eqref{bs2}
is attained in
\begin{equation}
\label{BK}
\frac{d\mathbb{Q}^{g}}{d\mathbb{P}}=\mathcal{D}
\left(\int_{0}^{T} H_{s} dW_{s}+\int_{0}^{T} \tilde{H}_{s}d\tilde{N}_{s} \right),
\qquad\mbox{with } (H_{s},\tilde{H}_{s})
\in\partial g(s,\mathcal{Z}_{s},\mathcal{\tilde{Z}}_{s}),
\end{equation}
where $\mathcal{D}$ denotes the Doleans-Dade exponential,
and $\partial g(s,\cdot)$ denotes the mapping of subdifferentials
of the convex driver $g(s,\cdot)$.
We shall exploit this to compute the lower bound numerically.
For simplicity, assume that $g(s,\cdot)$ is differentiable.
(If that is not satisfied, then our approach may still be applied
by considering elements in the subgradient.)

Let $N_{2}\in\mathbb{N}$, simulate paths $(W^{n}_{s_{jp}})$, $(N^{n}_{s_{jp}})$
and $(X^{n}_{s_{jp}})$ for $n=1,\ldots,N_{2}$,
and also consider (the true, non-simulated)
$(W_{s_{jp}})$, $(N_{s_{jp}})$ and $(X_{s_{jp}})$.
%In the sequel we assume that $\mathcal{D}_{j+1}(X_{t_j}^n) = \varphi(t_j,X_{t_j}^n)\Delta W_{j+1}$ for a suitable vector of function $\varphi$.
Then, define and construct
\begin{align*}
\mathcal{Z}^{N_{2}}_{t}:=  &  z^{N_{2}}_{t}(s_{jp},X_{s_{jp}})
:= {\beta}^{N_{2}}_{s_{jp}}\varphi(s_{jp},X_{s_{jp}})\mbox{ for }s_{jp}\leq t<s_{j(p+1)},\\
\mathcal{\tilde{Z}}^{N_{2}}_{t}:=  &  \tilde{z}^{N_{2}}_{t}(s_{jp},X_{s_{jp}})
: =\tilde{\beta}^{N_{2}}_{s_{jp}}\tilde{\varphi}(s_{jp},X_{s_{jp}})\mbox{ for }s_{jp}\leq t<s_{j(p+1)},
\end{align*}
using least squares Monte Carlo regression as described in Section \ref{sec:algorithm} with
$K_{2}$ basis functions
and terminal condition given by $U =\sum_{l=1}^{L} f_{\tau^{l}}(X_{\tau^{l}})$.
(Henceforth, we suppress $K_{2}$ in the notation.)
Furthermore, define and construct,
using $N_{3}\in\mathbb{N}$ new i.i.d. simulations,
\begin{align*}
\mathcal{Z}^{N_{2},n}_{t}:=  &  z^{N_{2}}_{t}(s_{jp},X^n_{s_{jp}})
: ={\beta}^{N_{2}}_{s_{jp}}\varphi(s_{jp},X^{n}_{s_{jp}})\mbox{ for }s_{jp}\leq t<s_{j(p+1)},\qquad n=1,\ldots,N_{3},\\
\mathcal{\tilde{Z}}^{N_{2},n}_{t}:=  &  \tilde{z}^{N_{2}}_{t}(s_{jp},X^n_{s_{jp}})
:= \tilde{\beta}^{N_{2}}_{s_{jp}}\tilde{\varphi}(s_{jp},X^n_{s_{jp}})\mbox{ for }s_{jp}\leq t<s_{j(p+1)},\qquad n=1,\ldots,N_{3},
\end{align*}
and moreover the partial derivatives
\begin{align*}
q^{n}_{s_{jp}}  &  :=g_{\mathcal{Z}}(s_{jp},\mathcal{Z}^{N_{2},n}_{s_{jp}} , \mathcal{\tilde{Z}}^{N_{2},n}_{s_{jp}}),\qquad n=1,\ldots,N_{3},\\
\lambda^{n}_{s_{jp}}-{\lambda}_{\mathbb{P}}
&  :=g_{\mathcal{\tilde{Z}}}(s_{jp},\mathcal{Z}^{N_{2},n}_{s_{jp}},
\mathcal{\tilde{Z}}^{N_{2},n}_{s_{jp}}),\qquad n=1,\ldots,N_{3}.
\end{align*}
Next, define $N_{3}$ i.i.d. simulations of the measure $\frac{d\mathbb{Q}^{\mathrm{approx}}}{d\mathbb{P}}$
via the Radon-Nikodym derivative
\begin{align*}
D^{n}  &  := \exp\bigg( \sum_{0\leq s_{jp}} q^{n}_{s_{jp}}\Delta W^{n}_{jp}
+\sum_{0\leq s_{jp}}\log\bigg( \frac{\lambda^{n}_{s_{jp}}}{\lambda_{\mathbb{P}}}\bigg)\Delta N^{n}_{jp}\\
&  \hspace{1.5cm}-\sum_{0\leq s_{jp}}\Big( \frac{1}{2} | q^{n}_{s_{jp}}|^{2}
+ \lambda^{n}_{s_{jp}}-{\lambda}_{\mathbb{P}}\Big)\Delta_{jp}\bigg),\qquad n=1,\ldots,N_{3}.
\end{align*}
Finally, set
\begin{align}
\label{eq:lowerbound}
\widetilde{Y}_{0}^{\mathrm{low},L}
:=\frac{1}{N_{3}}\sum_{n=1}^{N_{3}}D^{n} \sum_{l=1}^{L} f_{\tau^{l,n}}(X^{n}_{\tau^{l,n}}),
\end{align}
with $f_{\tau^{l,n}}(X^{n}_{\tau^{l,n}})$, $n=1,\ldots,N_{3}$,
%Monte Carlo
simulated copies of $f_{\tau^{l}}(X_{\tau^{l}})$
constructed by applying the numerical scheme of Section \ref{sec:algorithm}.
%[***HOW OBTAINED???***]
%where $M_{j}%
%^{\pi,M_{1},N_{1},n}$ should be simulated using %${\alpha}^{\pi,M_{1},N_{1}}$
%and ${\gamma}^{\pi,M_{1},N_{1}}$ estimated previously (see Section
%\ref{algorithm}).

Recall from \eqref{eq:low} that $\rho_{0}(U)$ gives a lower bound to the upper Snell envelope.
Thus, it follows from \eqref{bs2}, \eqref{BK},
and the definition of $U$
%and Theorem~\ref{th:mulsub} [***NEEDED?***]
that
\begin{align}
\mathbb{E}\left[\widetilde{Y}_{0}^{\mathrm{low},L}\right]
&=\mathbb{E}\left[\frac{1}{N_{3}}\sum_{n=1}^{N_{3}}D^{n} \sum_{l=1}^{L} f_{\tau^{l,n}}(X^{n}_{\tau^{l,n}})\right]\nonumber\\
&=\mathbb{E}_{\mathbb{Q}^{\mathrm{approx}}}\left[\sum_{l=1}^{L} f_{\tau^{l}}(X_{\tau^{l}})\right]
\leq\rho_{0}\left(\sum_{l=1}^{L} f_{\tau^{l}}(X_{\tau^{l}})\right)
\leq Y^{\ast,L}_{0}.
\label{eq:genuinelowerbound}
\end{align}
That is, our estimator \eqref{eq:lowerbound} constitutes a genuine lower bound.
%[***ADD SOME MORE FORMAL DETAILS???***]
This means that, on average, we indeed obtain a lower bound to the optimal solution
given by the upper Snell envelope.
Furthermore, as a consequence of Proposition \ref{PropConv} and Theorem \ref{theoConv},
the lower bound converges to the optimal solution.

Summarizing succinctly: given a time grid $\Pi$,
our explicit numerical construction of the lower bound consists of the following steps:
\begin{itemize}
\item[(1.)] Select $K_{1}$ basis functions and run $N_{1}$ Monte Carlo simulations to
determine $\overline{M}^{K_{1},N_{1}}$
and $\overline{c}^{K_{1},N_{1}}$.
To describe their evolution, it is sufficient to store the corresponding
$(\gamma^{N_{1}}_{s_{jp}})_{j,p}$,
$(\beta^{N_{1}}_{s_{jp}})_{j,p}$
and $(\tilde{\beta}^{N_{1}}_{s_{jp}})_{j,p}$.
[(BSDE 1).]

\item[(2.)] Use $(\gamma^{N_{1}}_{s_{jp}})_{j,p}$,
$(\beta^{N_{1}}_{s_{jp}})_{j,p}$ and
$(\tilde{\beta}^{N_{1}}_{s_{jp}})_{j,p}$
to estimate $\sum_{l=1}^{L} f_{\tau^{l}}(X_{\tau^{l}})$.
Select $K_{2}$ basis functions and run $N_{2}$ Monte Carlo simulations to
determine $(\mathcal{Z}^{N_{2}},\tilde{\mathcal{Z}}^{N_{2}})$
using the estimated $\sum_{l=1}^{L} f_{\tau^{l}}(X_{\tau^{l}})$ as terminal condition.
To describe the evolution of this process it is sufficient to store the corresponding
%$({\gamma}^{N_{1}}_{s_{jp}})_{j,p}$,
$(\beta^{N_{2}}_{s_{jp}})_{j,p}$ and
$(\tilde{\beta}^{N_{2}}_{s_{jp}})_{j,p}$.
[(BSDE 2).]

%	\item[(2.)] $({\alpha}^{\pi,M_{1},N_{1}}_{j})_{j}$ and $({\gamma}^{\pi
%		,M_{1},N_{1}}_{j})_{j}$ give rise to a terminal condition $U^{\pi,M_{1},N_{1}%
%	}$ and a Markov process $\mathcal{X}^{\pi,M_{1},N_{1}}$ defined above. Choose
%	$M_{2}$ basis functions, and use $N_{2}$ MC simulations to calculate
%	$(Y^{\pi,M_{2},N_{2}},\mathcal{Z}^{\pi,M_{2},N_{2}}) $, the solution of
%	corresponding BS$\Delta$Es with the Markov process $\mathcal{X}^{\pi
%		,M_{1},N_{1}}$ and terminal condition $U^{\pi,M_{1},N_{1}}$. Store the
%	corresponding $({\gamma}^{\pi,M_{2},N_{2}}_{j})_{j}.$

\item[(3.)] With $(\beta_{s_{jp}}^{N_{2}})_{j,p}$
and $(\tilde{\beta}_{s_{jp}}^{N_{2}})_{j,p}$ at hand from Step (2.),
simulate $N_{3}$ copies of $\frac{d\mathbb{Q}^{\mathrm{approx}}}{d\mathbb{P}}$.
Furthermore, with $(\gamma^{N_{1}}_{s_{jp}})_{j,p}$,
$(\beta^{N_{1}}_{s_{jp}})_{j,p}$ and
$(\tilde{\beta}^{N_{1}}_{s_{jp}})_{j,p}$ at hand
from Step (1.),
simulate $N_{3}$ copies of $\sum_{l=1}^{L} f_{\tau^{l,n}}(X^{n}_{\tau^{l,n}})$.
%	Moreover, with
%	$({\alpha}^{\pi,M_{1},N_{1}}_{j})_{j}$ and $({\gamma}^{\pi,M_{1},N_{1}}%
%	_{j})_{j}$ at hand from Step 1, $N_{3}$ copies of $U^{\pi,M_{1},N_{1}}$ can be
%	simulated.
Using \eqref{eq:lowerbound}, a genuine lower bound to the upper Snell envelope is then obtained.
\end{itemize}

Note that the Monte Carlo simulations are done consecutively and are not
nested, so that the total computation time depends only on $N_{1}+N_{2}+N_{3}$.
We summarize the results of this subsection in the following theorem:

\begin{theorem}
The estimator $\widetilde{Y}_{0}^{\mathrm{low},L}$ defined in \eqref{eq:lowerbound} is a
genuine lower bound to the upper Snell envelope, i.e.,
${\mathbb{E}}\left[\widetilde{Y}_{0}^{\mathrm{low},L}\right]  \leq Y^{\ast,L}_{0}$.
Furthermore, $\widetilde{Y}_{0}^{\mathrm{low},L}$ converges to $Y^{\ast,L}_{0}$
as $N_{1}$, $N_{2}$, $N_{3}$, $K_{1}$ and $K_{2}$ %and $K_{3}$
tend to infinity and $\Delta$ tends to zero, i.e.,
\begin{equation*}
\lim_{\Delta\rightarrow 0}
\lim_{K_{i}\rightarrow \infty, i=1,2}%,3}
\lim_{N_{i}\rightarrow \infty,i=1,2,3} \widetilde{Y}_{0}^{\mathrm{low},L}
=Y^{\ast,L}_{0}.
\end{equation*}
\label{th:LB}
\end{theorem}

\subsubsection{Subtracting the associated $\rho$-martingale in Eqn. \eqref{eq:lowerbound}:
Reducing the variance while not inducing a bias}

\label{sec:martsubt}

Let $(\mathcal{Z}^N,\mathcal{\tilde{Z}}^N)$
be approximations of the `true' solution to the BSDE with terminal condition $U$.
%These approximations have been obtained before and can be simulated through our basis functions and the parameters obtained previously by linear regression.
Denote by $\mathbb{Q}$ a probability measure defined by
\begin{equation}
\frac{d\mathbb{Q}}{d\mathbb{P}}=\mathcal{D}\left(  \int_{0}^{T} H_{s} dW_{s}+
\int_{0}^{T} \tilde{H}_{s} d\tilde{N}_{s} \right),\qquad\mbox{
	with } (H_{s},\tilde{H}_{s})\in\partial g(s,\mathcal{Z}^{N}_{s},\tilde
	{\mathcal{Z}}^N_{s} ).
\end{equation}
By \eqref{bs2},
$\mathbb{E}_{\mathbb{Q}}\left[U\right]$ yields a lower bound to $\rho_{0}(U)$. %, see the main paper [***ADAPT***].
(If the $\mathcal{Z}^N$ and $\mathcal{\tilde{Z}}^N$ were exact,
then $\mathbb{Q}=\mathbb{Q}^g$ and
$\mathbb{E}_{\mathbb{Q}}\left[U\right]= \mathbb{E}_{\mathbb{Q}^g}\left[U\right]=\rho_{0}(U)$.)

The following proposition will help to reduce the variation in our numerical scheme.
\begin{proposition}
\label{prop:subtracting}
Let $\overline{M}^{N}$ be the $\rho$-martingale defined by
\begin{equation*}
\overline{M}^{N}_t:=-\int_0^t g(s,\mathcal{Z}^{N}_s,\mathcal{\tilde{Z}}^{N}_s)ds
+\int_0^t \mathcal{Z}^{N}_s dW_s
+\int_0^t \mathcal{\tilde{Z}}^{N}_s d\tilde{N}_s.
\end{equation*}
Then, subtracting the $\rho$-martingale from the terminal condition does not induce a bias, i.e.,
\begin{equation*}
\mathbb{E}_{\mathbb{Q}}\left[U-\overline{M}^{N}_T\right]=\mathbb{E}_{\mathbb{Q}}\left[U\right].
\end{equation*}
\end{proposition}

We finally note that if we would have that $\mathbb{Q}=\mathbb{Q}^g$ from Eqn. \eqref{BK},
then, by the definition of $\overline{M}^{N}_t$,
$U-\overline{M}^{N}_T=\mathbb{E}_{\mathbb{Q}^g}\left[U\right]=\rho_{0}(U)$ were constant a.s.
Hence, if $\mathbb{Q}$ is approximately $\mathbb{Q}^g$,
then $U-\overline{M}^{N}_T$ is approximately constant.
More formally, if $\overline{M}^{N}_T$ converges to $\overline{M}_T$ in $L^2$,
we have that
\begin{align*}
\lim_{N \to \infty} U-\overline{M}^{N}_T = U- \overline{M}_T = \text{constant},
\end{align*}
where the convergence should be understood in $L^2$.
In particular, if $\overline{M}^{N}_T \rightarrow \overline{M}_T$ in $L^2$,
then $\mathrm{Var}(U-\overline{M}^{N}_T) \rightarrow 0 $ as $N \rightarrow \infty$.

Inspired by this theoretical result, which may be referred to as an almost sure property of a second kind
to distinguish it from the additive dual representation's almost sure property,
we will in our numerical analysis subtract the associated $\rho$-martingale
from the right-hand side of Eqn. \eqref{eq:lowerbound} when computing the genuine lower bound:
it will reduce the variance without inducing a bias.

\subsection{Converging Approximate and Genuine Upper Bounds}

\label{sec:UB}

This subsection develops an explicit \textit{approximate} upper bound
that converges to the upper Snell envelope asymptotically,
and an explicit \textit{genuine} upper bound that not only converges to the upper Snell envelope
but is also biased high at the pre-limiting level.
To obtain an upper bound with Theorem~\ref{th:mulsub},
we set
%in view of \eqref{eq:upp},
\begin{equation*}
U:=\max_{0\leq j_{1}<j_{2}<\cdots< j_{L}}\sum_{l=1}^{L}
\left(f_{j_{l}}(X_{j_{l}})-M_{j_{l}}^{\ast,L-l+1}+M_{j_{l-1}}^{\ast,L-l+1}\right).
\end{equation*}
Since we cannot compute $(M^{\ast})$ exactly,
we approximate, in view of \eqref{eq:upp},
the terminal condition $U$ by the $\rho$-martingale $(\overline{M}^{K_{1},N_{1}})$
constructed in Section \ref{sec:algorithm}, i.e.,
we set
\begin{equation*}
U^{N_{1}}:= \max_{0\leq j_{1}<j_{2}<\cdots< j_{L}}\sum_{l=1}^{L}
\left(f_{j_{l}}(X_{j_{l}})-\overline{M}_{j_{l}}^{L-l+1,K_{1},N_{1}}
+\overline{M}_{j_{l-1}}^{L-l+1,K_{1},N_{1}}\right).
\end{equation*}
Next, we define
\begin{align*}
\mathcal{X}_{t_{i}}^{N_{1}}
:= \left(X_{t_{i}},\overline{M}_{t_{i}}^{K_{1},N_{1}},
\max_{k \leq L: 0\leq j_{1}<j_{2}<\cdots<j_{k} \leq\lfloor{t_{i}}\rfloor}\sum_{l=1}^{k}
\left(f_{j_{l}}(X_{j_{l}})-\overline{M}_{j_{l}}^{L-l+1,K_{1},N_{1}}
+\overline{M}_{j_{l-1}}^{L-l+1,K_{1},N_{1}}\right)\right).
\end{align*}
Clearly, the terminal condition $U^{N_{1}}$ depends only on $\mathcal{X}_{T}^{N_{1}}$.
Next, note that, for every $t_{i}$, we have that $H_{t_{i}}$ is a function of $X_{t_{i}}$,
and, for every $r>i$,
$\overline{M}_{t_{r}}^{K_{1},N_{1}}$ only depends on $\overline{M}_{t_{i}}^{K_{1},N_{1}}$,
$(X_{t_{l}})_{i\leq l\leq r}$,
$(W_{t_{l}}-W_{t_{i}})_{i\leq l\leq r}$ and
$(N_{t_{l}}-N_{t_{i}})_{i\leq l\leq r}$ where both of the latter are
independent of $\mathcal{F}_{t_{r}}$.
From this we may conclude that $\mathcal{X}^{N_{1}}$ is a Markov process on the time grid $\Pi$.

Next, we solve numerically the BSDE (\ref{bs}) with $U^{N_{1}}$.
To do so we will consider an approximation scheme.
We can simulate paths of the adapted process%
\begin{align*}
&  (\mathcal{X}^{N_{1},n}_{t_{i}})_{i}=\bigg(X^{n}_{t_{i}},\overline{M}_{t_{i}}^{K_{1},N_{1},n},\\
&  \max_{k \leq L: 0\leq j_{1}<j_{2}<\cdots< j_{k} \leq\lfloor{t_{i}}\rfloor}\sum_{l=1}^{k}
\left(  f_{j_{l}}(X^{n}_{j_{l}})-\overline{M}_{j_{l}}^{L-l+1,K_{1},N_{1},n}+\overline{M}_{j_{l-1}}^{L-l+1,K_{1},N_{1},n}\right)\bigg)_{i},
\end{align*}
for $n=1,\ldots,N_{4}$.
To compute the BSDE with terminal condition $U^{N_{1}}$,
let $K_{4}$ be the number of basis functions in the least squares Monte Carlo regression.
Employing the algorithm described in Section \ref{sec:algorithm},
we can construct the coefficients
$\gamma_{s_{jp}}^{N_{4}}$,
$\beta_{s_{jp}}^{N_{4}}$ and
$\tilde{\beta}_{s_{jp}}^{N_{4}}$,
and processes
\begin{equation}
(\overline{Y}^{N_{4},L},\mathcal{Z}^{N_{4}},\tilde{\mathcal{Z}}^{N_{4}}).
\label{eq:conupp}
\end{equation}
Note that by applying Proposition \ref{PropConv} and Theorem \ref{theoConv} twice we may conclude that,
in $L^{2}$,
\begin{equation}
\label{convhaty}
%\lim_{M_{2},M_{1}\to\infty}
\lim_{\Delta\to 0}\lim_{K_{i}\to\infty,i=1,4}\lim_{N_{i}\to\infty,i=1,4}(\overline{Y}^{N_{4},L},\mathcal{Z}^{N_{4}},\tilde{\mathcal{Z}}^{N_{4}})
= (Y^{\ast,L},\mathcal{Z}^{\ast},\tilde{\mathcal{Z}}^{\ast}).
\end{equation}
In particular, $\overline{Y}^{N_{4},L}$ constitutes a converging approximate upper bound
to the upper Snell envelope.
%as the mesh ratio of the grid, $\pi,$ tends to zero, and the number of basis
%functions and the number of simulations tend to infinity.
%So in principle our
%algorithm will converge to the \emph{true} value of the $\rho$-Snell envelope.
%However, for finite computation time our estimates for the upper bound may be
%biased meaning that in the average our upper bound may not provide enough
%protection. For this reason we will look at algorithms yielding estimates for
%an upper bound which are unbiased so that with finite computation time in the
%average the company applying the upper bound will be `safe'. Therefore, we
%will \emph{not} take $Y^{\pi,M_{2},N_{2}}_{0}$ as the estimator for the upper
%bound. However, we will use the $({\gamma}^{\pi,M_{2},N_{2}}_{j})_{j}$ to
%obtain unbiased upper bounds.
From now on we assume that we have already estimated $\overline{Y}^{N_{4}}_{0}$,
$(\beta_{s_{jp}}^{N_{4}})_{j,p}$ and
$(\tilde{\beta}_{s_{jp}}^{N_{4}})_{j,p}$.

Next, let us develop a genuine upper bound.
It is well-known that under assumption (A1), the functional $\rho$ is
Lipschitz continuous (cf. Peng \cite{P97}).
%Throughout this section we assume that $\rho$ is Lipschitz continuous in
%$L^{p}$ for a $p\geq1,$ with a Lipschitz constant, say $K$. We do need the
%concavity assumption anymore. Note that it is well known that $\rho$ is always
%Lipschitz continuous in $L^{2}$ if the corresponding driver $g$ is Lipschitz continuous.
Define
\begin{align}
\label{fhat}\hat{U}^{N_{4}}:  &  =\overline{Y}^{N_{4},L}_{0} -\int
_{0}^{T}g(s,\mathcal{Z}^{N_{4}}_{s},\mathcal{\tilde{Z}}^{N_{4}}_{s})ds +\int_{0}^{T}\mathcal{Z}^{N_{4}}_{s}dW_{s} +
\int_{0}^{T}\mathcal{\tilde{Z}}^{N_{4}}_{s}d\tilde{N}_{s}\nonumber\\
&  =\overline{Y}^{N_{4},L}_{0}-\sum_{j,p} g(s_{jp},\beta_{s_{jp}}^{N_{4}} \psi(s_{jp},\mathcal{X}^{N_{1}}_{s_{jp}}),
\tilde{\beta}_{s_{jp}}^{N_{4}} \tilde{\psi}(s_{jp},\mathcal{X}%
^{N_{1}}_{s_{jp}}))\Delta_{s_{jp}}\nonumber\\
&  \hspace{1cm}+\sum_{j,p} \beta_{s_{jp}}^{N_{4}} \psi%
(s_{jp},\mathcal{X}^{N_{1}}_{s_{jp}})\Delta W_{s_{jp}} +\sum_{j,p} \tilde{\beta}_{s_{jp}}^{N_{4}} \tilde{\psi}(s_{jp},\mathcal{X}^{N_{1}}_{s_{jp}})\Delta\tilde{N}_{s_{jp}}.
\end{align}
Then, by Theorem \ref{th:mulsub} and the Lipschitz continuity of $\rho$,
we have that
\begin{align}
\label{inequality}
Y^{\ast,L}_{0}\leq\rho_{0}(U^{N_{1}})  &  \leq\rho_{0}(\hat{U}^{N_{4}})
+\mathcal{K}||\hat{U}^{N_{4}}-U^{N_{1}}||_{2}\nonumber\\
&  =\overline{Y}^{N_{4},L}_{0}
+\mathcal{K}||\hat{U}^{N_{4}}-U^{N_{1}}||_{2} ,
\end{align}
for a Lipschitz constant $\mathcal{K}$ analyzed later.
We will exploit inequality \eqref{inequality} to develop our genuine upper bound.
By Proposition \ref{PropConv} and Theorem \ref{theoConv},
$\overline{Y}^{N_{4},L}_{0}$ converges to $Y^{\ast,L}_{0}$,
$(\mathcal{Z}^{N_{4}}, \mathcal{\tilde{Z}}^{N_{4}})$
converges to $(\mathcal{Z}^{*},\mathcal{\tilde{Z}}^{*})$ in $L^{2}(d\mathbb{P}\times ds)$ and
$||\hat{U}^{N_{4}}-U^{N_{1}}||_{2}$ converges to zero,
as $N_{1}$, $N_{4}$, $K_{1}$ and $K_{4}$ tend to infinity
and $\Delta$ tends to zero.

%Through (\ref{inequality}) we can now derive an estimator for the upper bound
%$\rho_{0}(U^{N_{1}}).$
%Specifically,
For $N_{5}\in\mathbb{N}$,
simulate i.i.d. copies of $U^{N_{1}}$ through
\begin{equation*}
U^{N_{1}}_{n}=\max_{0\leq j_{1}<j_{2}<\cdots< j_{L}} \sum_{l=1}%
^{L}\left(  f_{j_{l}}(X^{n}_{j_{l}})-\overline{M}_{j_{l}}^{L-l+1,K_{1},N_{1},n}%
+\overline{M}_{j_{l-1}}^{L-l+1,K_{1},N_{1},n}\right),\quad n=1,\ldots,N_{5}.
\end{equation*}
Next, using \eqref{fhat}, simulate i.i.d. copies
$\hat{U}^{N_{4}}_{1},\hat{U}^{N_{4}}_{2},\ldots,\hat{U}^{N_{4}}_{N_{5}}$
of $\hat{U}^{N_{4}}$.
(Recall that $\overline{Y}^{N_{4},L}_{0}$,
$(\beta_{s_{jp}}^{N_{4}})_{j,p}$ and
$(\tilde{\beta}_{s_{jp}}^{N_{4}})_{j,p}$
are already available.)
Then, (\ref{inequality}) suggests to estimate the upper bound $\rho_{0}(U^{N_{1}})$ by
\begin{equation*}
\overline{Y}^{N_{4},L}_{0}+\mathcal{K} \sqrt{\frac{1}{N_{5}}\sum_{n=1}^{N_{5}} |\hat
{U}^{N_{4}}_{n}-U^{N_{1}}_{n}|^{2}}.
\end{equation*}
Note that $\frac{1}{N_{5}}\sum_{i=1}^{N_{5}} |\hat{U}^{N_{4}}%
_{i}-U^{N_{1}}_{i}|^{2}$ is an unbiased estimator of ${\mathbb{E}%
}\left[  |U^{N_{1}}-\hat{U}^{N_{4}}|^{2}\right]  $.
However, taking the square root of an estimator gives rise to a possible downward bias.
If we wish to eliminate the downward bias, we need to elaborate a little further, as follows.
We simulate independent $(\hat{U}^{N_{4}}_{i})_{i=1,\ldots, 2N_{5}}$.
Then we set
\begin{equation}
\label{eq:genupp}
\widetilde{Y}_{0}^{\mathrm{upp},L}: = \overline{Y}^{N_{4},L}_{0}+\mathcal{K} \frac{\frac{1}{N_{5}}
\sum_{n=1}^{N_{5}} |\hat{U}^{N_{4}}_{n}-U^{N_{1}}_{n}
|^{2}}{\sqrt{\frac{1}{N_{5}}\sum_{n=N_{5}+1}^{2N_{5}} |\hat{U}^{N_{4}}_{n}-U^{N_{1}}_{n}|^{2}}}.
\end{equation}
Since
\begin{align*}
 {\mathbb{E}}&\left[  \frac{\frac{1}{N_{5}}\sum_{n=1}^{N_{5}} |\hat{U}
^{N_{4}}_{n}-U^{N_{1}}_{n}|^{2}}{\sqrt{\frac{1}{N_{5}}
\sum_{n=N_{5}+1}^{2N_{5}} |\hat{U}^{N_{4}}_{n}-U^{N_{1}
}_{n}|^{2}}}\right] \\
&  = {\mathbb{E}}\left[  \frac{1}{N_{5}}\sum_{n=1}^{N_{5}}
|\hat{U}^{N_{4}}_{n}-U^{N_{1}}_{n}|^{2}\right]
{\mathbb{E}}\left[  \frac{1}{\sqrt{\frac{1}{N_{5}}\sum_{n=N_{5}+1}^{2N_{5}}
|\hat{U}^{N_{4}}_{n}-U^{N_{1}}_{n}|^{2}}}\right] \\
&  \geq||\hat{U}^{N_{4}}_{n}-U^{N_{1}}%
_{n}||^{2} \frac{1}{\sqrt{\frac{1}{N_{5}}\sum_{n=N_{5}+1}^{2N_{5}}{\mathbb{E}
}\left[  |\hat{U}^{N_{4}}_{n}-U^{N_{1}}_{n}|^{2}\right]
}}\\
&  = ||\hat{U}^{N_{4}}_{n}-U^{N_{1}}%
_{n}||^{2} \frac{1}{\sqrt{{\mathbb{E}}\left[  |\hat{U}^{N_{4}}
_{n}-U^{N_{1}}_{n}|^{2}\right]  }}
= ||\hat{U}^{N_{4}}_{n}-U^{N_{1}}_{n}||,
\end{align*}
$\widetilde{Y}_{0}^{\mathrm{upp},L}$ thus defined is biased high.
Note that the first equality follows from independence and the inequality is due to
a suitable (reversed) application of
Jensen's inequality.
As before, as a consequence of Proposition \ref{PropConv} and Theorem \ref{theoConv},
we have that $\widetilde{Y}_{0}^{\mathrm{upp},L}$ converges to $Y^{\ast,L}_{0}$ as
%$\pi$ tends to zero and $M_{1}$, $M_{2}$,
$N_{1}$, $N_{4}$, $N_{5}$,
$K_{1}$ and $K_{4}$ %$K_{5}$,
tend to infinity
and $\Delta$ tends to zero.

Summarizing succinctly:
given a time grid $\Pi$,
our explicit numerical construction of the upper bound
consists of the following steps:

\begin{itemize}
\item[(1.)] Select $K_1$ basis functions and run $N_{1}$ Monte Carlo simulations to
determine $\overline{M}^{K_{1},N_{1}}$
and $\overline{c}^{K_{1},N_{1}}$.
To describe their evolution, it is sufficient to store the corresponding
$(\gamma^{N_{1}}_{s_{jp}})_{j,p}$,
$(\beta^{N_{1}}_{s_{jp}})_{j,p}$ and
$(\tilde{\beta}^{N_{1}}_{s_{jp}})_{j,p}$.
[(BSDE 1).]

\item[(2.)] $(\gamma^{N_{1}}_{s_{jp}})_{j,p}$,
$(\beta^{N_{1}}_{s_{jp}})_{j,p}$ and
$(\tilde{\beta}^{N_{1}}_{s_{jp}})_{j,p}$
give rise to a terminal condition $U^{N_{1}}$ and a Markov process $\mathcal{X}^{N_{1}}$ defined above.
Select $K_{2}$ basis functions and run $N_{4}$ Monte Carlo simulations to calculate
$(\overline{Y}^{N_{4},L},\mathcal{Z}^{N_{4}}, \mathcal{\tilde{Z}}^{N_{4}})$
as the solution to the corresponding BS$\Delta$Es with the
Markov process $\mathcal{X}^{N_{1}}$ and terminal condition
$U^{N_{1}}$.
Store the corresponding $\overline{Y}^{N_{4},L}_{0}$,
$(\beta^{N_{4}}_{s_{jp}})_{j,p}$ and
$(\tilde{\beta}^{N_{4}}_{s_{jp}})_{j,p}$.
[(BSDE 2).]

\item[(3.)] With $\overline{Y}^{N_{4},L}_{0}$,
$(\beta^{N_{4}}_{s_{jp}})_{j,p}$ and
$(\tilde{\beta}^{N_{4}}_{s_{jp}})_{j,p}$
at hand from Step (2.), simulate $N_{5}$ copies of $\hat{U}^{N_{4}}$
defined by (\ref{fhat}).
Furthermore, with
$(\gamma^{N_{1}}_{s_{jp}})_{j,p}$,
$(\beta^{N_{1}}_{s_{jp}})_{j,p}$ and
$(\tilde{\beta}^{N_{1}}_{s_{jp}})_{j,p}$ at hand
from Step (1.),
simulate $N_{5}$ copies of $U^{N_{1}}$.
Using \eqref{eq:genupp}, a genuine upper bound to the upper Snell envelope is then obtained.
\end{itemize}

Note again that the Monte Carlo simulations are done consecutively
and are not nested,
so that the total computation time depends only on $N_{1}+N_{4}+N_{5}$.
We summarize the results of this subsection in the following theorem:

\begin{theorem}
The estimator $\widetilde{Y}_{0}^{\mathrm{upp},L}$ defined in
\eqref{eq:genupp} is a genuine upper bound to the upper Snell envelope,
i.e., ${\mathbb{E}}\left[\widetilde{Y}_{0}^{\mathrm{upp},L}\right]\geq Y^{\ast,L}_{0}$.
Furthermore, $\widetilde{Y}_{0}^{\mathrm{upp},L}$ converges to
$Y^{\ast,L}_{0}$ as $N_{1}$, $N_{4}$, $N_{5}$, $K_{1}$ and $K_{4}$ %and $K_{5}$
tend to infinity and $\Delta$ tends to zero,
i.e.,
\begin{equation*}
\lim_{\Delta\rightarrow 0}
\lim_{K_{i}\rightarrow \infty, i=1,4}%,5}
\lim_{N_{i}\rightarrow \infty,i=1,4,5} \widetilde{Y}_{0}^{\mathrm{upp},L}
=Y^{\ast,L}_{0}.
\end{equation*}
Moreover, the estimator $\overline{Y}_{0}^{N_{4},L}$ defined in \eqref{eq:conupp} also converges to
$Y^{\ast,L}_{0}$ as $N_{1}$, $N_{4}$, $K_{1}$ and $K_{4}$ tend to infinity and $\Delta$ tends to zero.
\label{th:UB}
\end{theorem}

%\subsubsection{The precise Lipschitz constant}
We finally state the following proposition on the precise Lipschitz constant $\mathcal{K}$
appearing in \eqref{inequality}.
\begin{proposition}
\label{prop:Lipschitz}
Let $\xi$ and $\xi'$ be square-integrable terminal conditions and denote by $(Y,Z,\tilde{Z})$
and $(Y',Z',\tilde{Z}')$ the associated BSDE solutions.
Then,
\begin{equation}
	|Y_0-Y'_0|^2 \leq \exp(\mathcal{L}^2 T)\E{|\delta \xi|^2},
\end{equation}
where $\delta \xi:= \xi-\xi'$
and with $\mathcal{L}$ the Lipschitz constant of the driver $g$.
\end{proposition}

\setcounter{equation}{0}

\section{Numerical Examples}

\label{sec:num}

%{\footnotesize
%\begin{itemize}
%\item Note that the algorithm proposed in this paper is linear and requires no
%nested simulation.
%
%\item Report CPU times.
%
%\item Compare our gaps (not the bounds themselves) to benchmark gaps.
%
%\item In addition to tables with the bounds and their standard errors for the
%various examples (univariate, multivariate, single, multiple, 4 examples of
%$\rho$) prepare plots of:
%
%\begin{itemize}
%\item The exercise region (thresholds).
%
%\item Densities of the stopping times.
%\end{itemize}
%
%and compare to the non-ambiguity averse case.
%\end{itemize}
%}

In this section we analyze our approach in numerical examples,
including single and multiple stopping, univariate and multivariate stochastic drivers,
increasing and decreasing reward functions,
and pure diffusion and jump-diffusion models.
As a general observation, we recall that the computational complexity of the
numerically implementable method proposed in this paper
is linear in the number of exercise rights
and that it does not require nested simulation.

\subsection{Single Stopping: Bermudan Option in a Diffusion Model}

The first example studied in this subsection is the pricing problem for a Bermudan-style option
in a single risky asset Black-Scholes model with dividends in the presence of ambiguity.
This example goes back to Andersen and Broadie \cite{AB04} in a setting without ambiguity,
%,
%and has since been considered many times.
%, e.g., in Kr\"atschmer \textit{et al.} \cite{KLLSS18}.
which can serve as a benchmark case.
Throughout this subsection, we consider the dynamics
$\tfrac{dX_{t}^{i}}{X_{t}^{i}}=\mu^{i} dt+\sigma^{i} dW_{t}^{i}$, %+J^{i}d\tilde{N}_{t}^{i}$,
$i=1,\ldots,d$, with $\mu^{i}\in\mathbb{R}$ and $\sigma^{i}\in\mathbb{R}_{>0}$,
%and $J^{i}\in(-1,\infty)$; cf.~\eqref{eq:markov}.
where $X_{t}^{i}$ is the price of asset $i$ at time $t$.
Following this literature, we assume that
there is a risk-free interest rate of $\rho=0.05$,
and that the option's underlying is a single dividend-paying stock $X$ with constant volatility $\sigma=0.2$ and dividend rate $\delta=0.1$,
resulting in a risk-neutral drift of $\mu=\rho-\delta=-0.05$.
%,
%and suppose that $J\equiv \lambda\equiv 0$.
The first product we study is a call option with strike price $K=100$ and maturity $T=3$.
The stock price at time 0 is varied between $x_0\in \{90,100,110\}$.
Exercise dates are specified as $t_{j}=\tfrac{jT}{10}$, $j=0,1,\ldots,10$,
i.e., there are 9 intermediate exercise dates, and the trivial ones at time $t=0$ and at maturity.
We allow for ambiguity in the drift and consider $g(t,z)=\delta_{1}|z|$, $\delta_{1}>0$;
cf. Example~\ref{g-ex}.

The three BSDEs---$\rho$-martingale, lower, and upper bounds---are solved with two sets of simulations.
We consider one set of simulations with $100\mathord{,}000$ trajectories and $1\mathord{,}000$ time steps
for the initial $\rho$-martingale BSDE,
and a second set of simulations with $100\mathord{,}000$ trajectories and $1\mathord{,}000$ time steps
for the BSDEs associated with the lower and upper bounds.
The number of basis functions is always 52 for $Y$ and 52 for $Z$.
These are $1$, $x$, $(x-q_{it})^+$ where the $q_{it}$ are $1,3,5,\ldots,99$ percent quantiles of $X_t$
estimated from the trajectories in the initial run of least squares Monte Carlo.
Thus, we approximate all unknown functions by linear splines.
Notice that this function basis is not problem-dependent and that its precision can be controlled by the chosen grid of quantiles.
Numerically, this basis works much better in our experiments than a locally linear approximation,
which would also include the discontinuous terms $1_{\{x>q_{it}\}}$.
In all our numerical experiments, the implementation of the regression is based on \eqref{defag2}.

For the evaluation of the lower bound, we draw a new sample with a larger number of $400\mathord{,}000$ trajectories.
For the evaluation of the upper bound, we need a very fine time discretization
to obtain a small `tracking error', defined as
$\sqrt{\frac{1}{N_{5}}\sum_{i=1}^{N_{5}} |\hat{U}^{N_{4}}_{i}-U^{N_{1}}_{i}|^{2}}$; cf.~\eqref{eq:genupp}.
We need, however, fewer trajectories,
because the pathwise dual representation leads to a very small variance,
as expected from the theoretical results in Section~\ref{sec:suroptmart}.
We thus choose the number of trajectories to be $1\mathord{,}000$
and increase the number of time steps by a factor 100 to $100\mathord{,}000$.
Here, we do not run new regressions but simply repeat each set of coefficients 100 times.
This device of increasing the time discretization and extrapolating regression coefficients is due to Belomestny, Bender, Schoenmakers \cite{BBS09}, %Section 6,
in a standard stopping setting without ambiguity.
However, due to the additional non-linearity from the BSDE,
it is not sufficient in our setting to just run the regressions at exercise times only;
we need them at our fine grid $\Pi$.

Tables~\ref{tab:singleunicall-1}--\ref{tab:singleunicall-3} summarize lower and upper bounds
for different degrees of ambiguity ($\delta_{1}$) and different values of the initial stock price ($x_0$).
The first four columns display lower bounds along with their standard errors.
Here, ``LB without M'' is the lower bound in Eqn.~\eqref{eq:lowerbound}
without subtraction of the $\rho$-martingale,
while ``LB'' is the definitive lower bound that subtracts the $\rho$-martingale,
as discussed in Section~\ref{sec:martsubt}.
As the results confirm, subtraction of the $\rho$-martingale leads to a substantial variance reduction.
The next three columns correspond to the upper bound.
Here, $\overline{Y}_0^{N_{4}}$ is the solution to the upper bound BSDE,
i.e., the approximate upper bound,
``TE'' is the tracking error defined above,
and ``UB'' combines the two to obtain the genuine upper bound in Eqn.~\eqref{eq:genupp}.
The final two columns display the mean of $\overline{Y}_T^{N_{4}}=U^{N_{1}}$ together with its standard error.
The point here is to illustrate that having only $1\mathord{,}000$ trajectories
in the upper bound simulation is already sufficient,
since the terminal condition has a (very) small variance.
In general, we observe that the gaps between LB on the one hand
and $\overline{Y}_0^{N_{4}}$ and UB on the other hand
are (very) small.
When $1/\delta_{1}\equiv\infty$ (i.e., in the no-ambiguity case),
we can compare our results to benchmark values, e.g., from Andersen and Broadie \cite{AB04}.
With $x_{0}=100$, the true value is 7.98.
This should be compared to our genuine lower bound of 7.98,
our approximate upper bound of 8.00, and our genuine upper bound of 8.07,
corresponding to gaps of 0.2\% and 1.0\% of the option value, respectively.

\begin{table}[h!]
\begin{center}
{\small
\begin{tabular}{c||cc|cc|ccc|cc}
$1/\delta_{1}$ & LB without M & s.e. & LB & s.e. & $\overline{Y}_0^{N_{4}}$ & TE & UB & $\mathbb{E}[\overline{Y}_T^{N_{4}}]$ & s.e. \\
\hline\hline
%1 &  38.1997 &  0.8383  & 38.5052  &  0.2502 &  38.3716  &  0.0735  & 38.7010  & 38.3092   & 0.0023\\
%3 &   9.4425  &  0.0558  &  9.4011 &   0.0154&    9.3836 &   0.0690 &   9.4651 &   9.3728  &  0.0022\\
10&    5.4901 &   0.0197  &  5.4635 &   0.0107&    5.4833 &   0.0843 &   5.5689 &   5.4679  &  0.0026\\
30&   4.7384 &   0.0158  &  4.7161 &   0.0096 &   4.7326 &   0.0836 &   4.8163 &   4.7184  &  0.0025\\
100 & 4.5005 &   0.0147  &  4.4795 &   0.0093 &   4.4962 &   0.0812 &   4.5775 &   4.4830  &  0.0025\\
300 &   4.4349 &   0.0144  &  4.4143 &   0.0092 &   4.4310 &   0.0806 &   4.5116 &   4.4180  &  0.0024\\
$1\mathord{,}000$&   4.4121 &   0.0143  &  4.3917 &   0.0092 &   4.4084 &   0.0804 &   4.4888  &  4.3955  &  0.0024\\
$10\mathord{,}000$ &   4.4030 &   0.0142  &  4.3827 &   0.0092 &   4.3997 &   0.0803 &   4.4800  &  4.3868  &  0.0024\\
$\infty$ &   4.4024 &   0.0142  &  4.3821 &   0.0092 &   4.3987 &   0.0803 &   4.4790  &  4.3859  &  0.0024\\
\end{tabular}
}
\caption{Bounds for $x_0=90$}
\label{tab:singleunicall-1}
%\end{table}
%\begin{table}[h!]
{\small
\vskip 0.1cm
\begin{tabular}{c||cc|cc|ccc|cc}
$1/\delta_{1}$ & LB without M & s.e. & LB & s.e. & $\overline{Y}_0^{N_{4}}$ & TE & UB & $\mathbb{E}[\overline{Y}_T^{N_{4}}]$ & s.e. \\
\hline\hline
%1&50.4479  &  0.9888 &  50.7803   & 0.3295  & 50.6038 &   0.0807   &50.9656  & 50.5404  &  0.0025\\
%3&   14.6337 &   0.0677 &  14.5727  &  0.0212 &  14.5527  &  0.0704  & 14.6359 &  14.5350 &   0.0022\\
10&    9.4387  &  0.0244 &   9.4069  &  0.0143 &   9.4331  &  0.0721  &  9.5063 &   9.4110 &   0.0021\\
30& 8.4419 &   0.0199   & 8.4161   & 0.0132  &  8.4445  &  0.0700 &   8.5146  &  8.4232   & 0.0020\\
 100&   8.1305  &  0.0186  &  8.1068   & 0.0129 &   8.1316  &  0.0678  &  8.1994  &  8.1113 &   0.0019\\
300&    8.0426  &  0.0183   & 8.0195  &  0.0128  &  8.0450  &  0.0678  &  8.1127 &   8.0249 &   0.0019\\
$1\mathord{,}000$&    8.0152  &  0.0182   & 7.9923  &  0.0128  &  8.0150  &  0.0678  &  8.0828 &   7.9949 &   0.0019\\
$10\mathord{,}000$&    8.0039  &  0.0181   & 7.9811  &  0.0127  &  8.0034  &  0.0678  &  8.0712 &   7.9833 &   0.0019\\
$\infty$&    8.0030  &  0.0181   & 7.9802  &  0.0127  &  8.0022  &  0.0678  &  8.0699 &   7.9821 &   0.0019\\
\end{tabular}
}
\caption{Bounds for $x_0=100$}
\label{tab:singleunicall-2}
%\end{table}
%\begin{table}[h!]
{\small
\vskip 0.1cm
\begin{tabular}{c||cc|cc|ccc|cc}
$1/\delta_{1}$ & LB without M & s.e. & LB & s.e. & $\overline{Y}_0^{N_{4}}$ & TE & UB & $\mathbb{E}[\overline{Y}_T^{N_{4}}]$ & s.e. \\
\hline\hline
%1&   63.2316&    1.1398&   63.5934&    0.4124&   63.3735  &  0.0879  & 63.7673  & 63.3101 &   0.0027\\
%3&   20.8274&    0.0782&   20.7538&    0.0282 &  20.7452  &  0.0733  & 20.8318  & 20.7226 &   0.0022\\
10&   14.7691&    0.0279&   14.7380&    0.0179&   14.7718  &  0.0682  & 14.8411  & 14.7460 &   0.0020\\
30 &  13.6662&    0.0228&   13.6414&    0.0166  & 13.6814   & 0.0669  & 13.7484  & 13.6570 &   0.0019\\
100 &  13.3221&    0.0214&   13.2994&    0.0163  & 13.3400   & 0.0658  & 13.4058  & 13.3164 &   0.0019\\
300  & 13.2257&    0.0211&   13.2037&    0.0162 &  13.2459   & 0.0658  & 13.3116  & 13.2225 &   0.0019\\
$1\mathord{,}000$  & 13.1924&    0.0209&   13.1706&    0.0162  & 13.2133   & 0.0657  & 13.2790  & 13.1900 &   0.0019\\
$10\mathord{,}000$  & 13.1791&    0.0209&   13.1574&    0.0162  & 13.2007   & 0.0657  & 13.2664  & 13.1775 &   0.0019\\
$\infty$&   13.1774&    0.0209&   13.1556&    0.0162  & 13.1994   & 0.0657  & 13.2650  & 13.1761 &   0.0019\\
\end{tabular}
}
\caption{Bounds for $x_0=110$}
\label{tab:singleunicall-3}
\end{center}
\end{table}

\begin{figure}
	\centering
		\includegraphics[scale=0.52]{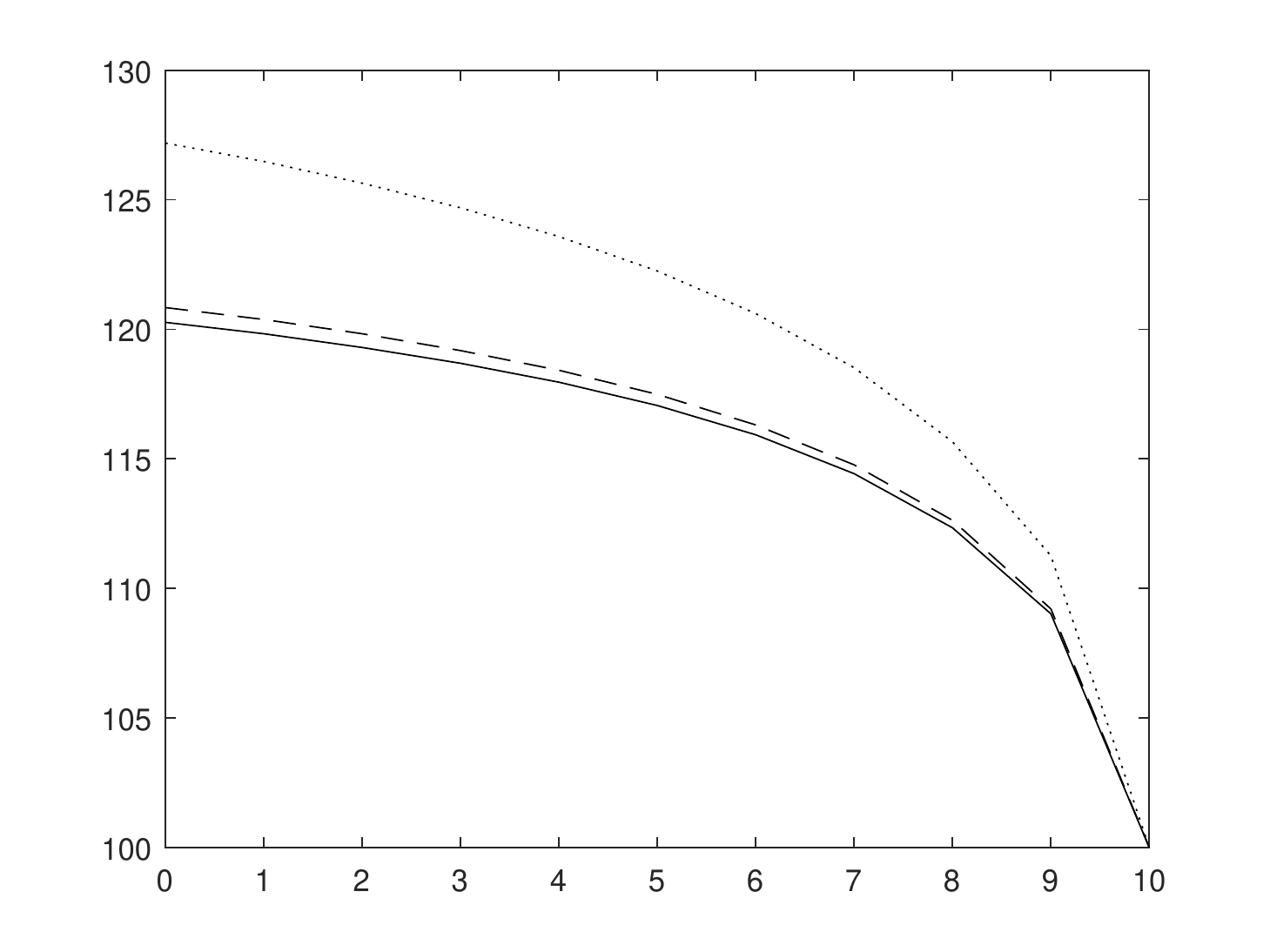}
        \includegraphics[scale=0.52]{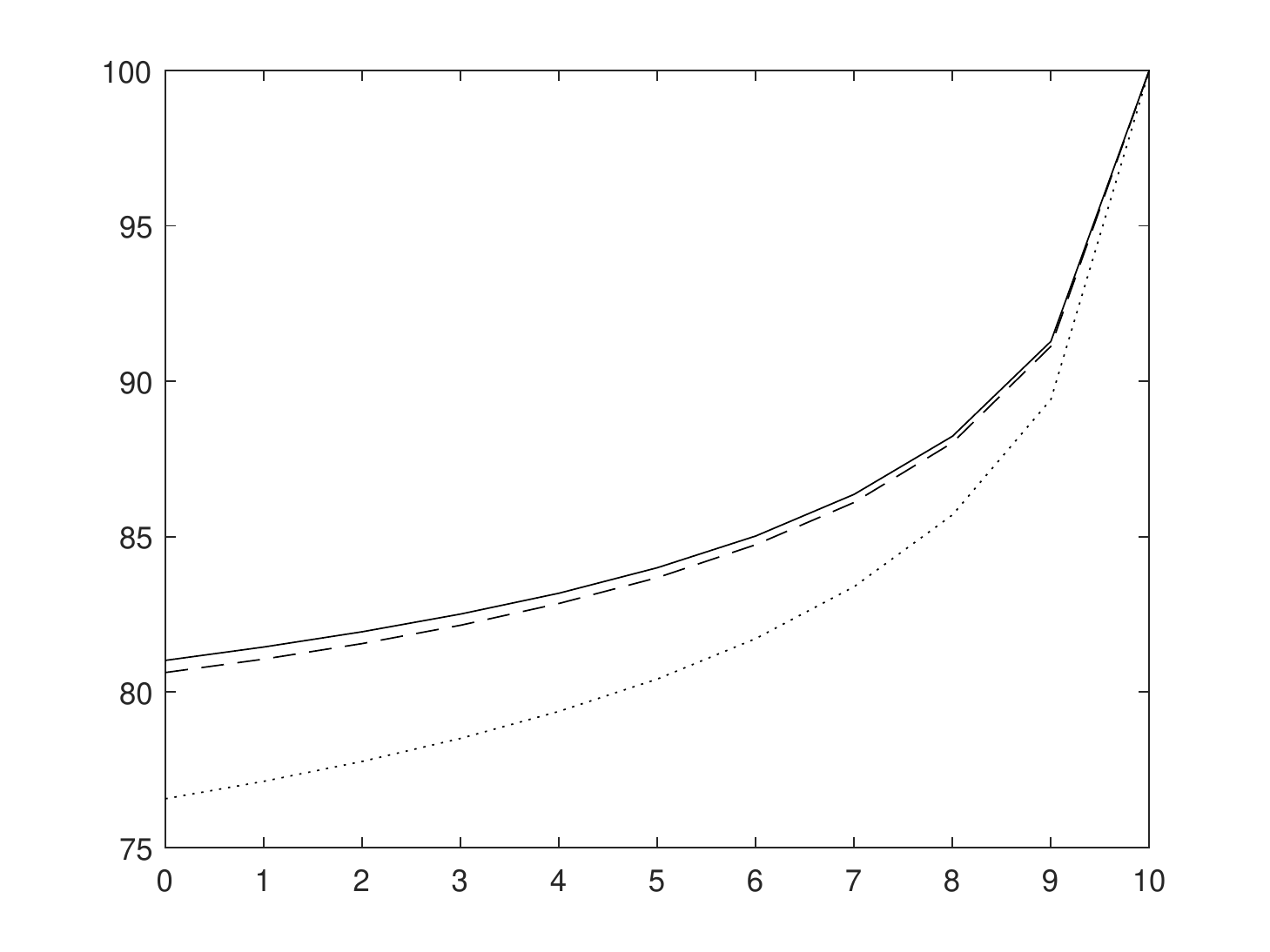}
	\caption{Exercise boundaries (left panel: call; right panel: put; solid: $1/\delta_1=\infty$;
dashed: $1/\delta_1=100$; and dotted: $1/\delta_1=10$).}
	\label{fig:singleuni}
\end{figure}

In %Tables~\ref{tab:singleunicall-4}--\ref{tab:singleunicall-6},
Figure~\ref{fig:singleuni}
we analyze exercise boundaries.
These boundaries are theoretically independent of the starting values of $x_{0}$.
In order to obtain accurate exercise boundaries already for early time points,
we need to make sure that the simulated trajectories are sufficiently spread out over the entire time span from $0$ to $T$
and are not concentrated around a fixed value of $x_0$ close to time $0$.
To this end, we start simulating trajectories at time $-1$
and set the drift of $X$ to zero over the interval $[-1,0]$.
Moreover, we double the number of basis functions by choosing a finer grid of quantiles, $0.005, 0.0015, \ldots, 0.995,$
and double the number of trajectories to $200\mathord{,}000$,
to account for the fact that we now have a wider space over which we approximate.

%Each column gives, for a specific exercise date,
In the left panel of Figure~\ref{fig:singleuni},
we plot the threshold value that the stock price has to exceed to make stopping optimal (i.e., the continuation value),
for different degrees of ambiguity ($\delta_{1}$),
as a function of the exercise dates.
%Exercise boundaries should, in principle, not depend on the value of $x_0$.
%However, due to the fact that all simulations start in $x_0$, which may be far from the exercise boundary,
%we cannot expect our estimated stopping decisions to be accurate at early time points.
%%Basically, all that matters is whether $x_0$ is above or below the exercise boundary.
%One can see this nicely from the three tables: from some point onwards (about after the third column),
%the impact of $x_0$ has disappeared and all three tables become nearly the same.
%%(which is actually a nice check of the stability of the method -- a priori we know very little about convergence of stopping times).
The solid line depicts exercise boundaries without ambiguity ($1/\delta_1=\infty$),
while the dashed line corresponds to $1/\delta_1=100$ and the dotted line corresponds to $1/\delta_1=10$.
In all cases, continuation values increase as we move away from the no-ambiguity case ($1/\delta_{1}=\infty$).
Thus, with more ambiguity, the decision-maker stops later.
Intuitively, for a call option's payoff function,
the possibility of a (prosperous) deviation from the reference model which accumulates over time
makes the option more valuable at later dates.
Hence, the decision-maker will stop later.

A similar pattern emerges when we replace the call option payoff function by a put option,
with strike price $K=100$ and $x_{0}=100$.
To make exercise decisions non-trivial, we set the dividend rate equal to zero in this case,
all else equal.
Table~\ref{tab:singleuniput-1} and Figure~\ref{fig:singleuni} (right panel)
%\ref{tab:singleuniput-2}
display the corresponding price bounds
and exercise boundaries.
For the put option, contrary to for the call option,
exercise becomes optimal
when the stock price falls \textit{below} the exercise boundaries in the right panel of Figure~\ref{fig:singleuni}.
Similar to the call option, with more ambiguity, the decision-maker stops later.
In this case, the possibility of an unfavorable deviation from the reference model that accumulates over time
makes the option more valuable at later dates.

\begin{table}[h!]
\begin{center}
{\small
\begin{tabular}{c||cc|cc|ccc|cc}
$1/\delta_1$ & LB without M & s.e. & LB & s.e. & $\overline{Y}_0^{N_{4}}$ & TE & UB & $\mathbb{E}[\overline{Y}_T^{N_{4}}]$ & s.e. \\
\hline\hline
%$     1 $  & $ 32.0447 $  & $  0.35965 $  & $ 32.1277 $  & $   0.20724 $  & $ 32.2991 $  & $ 0.071821 $  & $  32.621 $  & $ 32.2684 $  & $ 0.0022316 $  \\
%$     3 $  & $ 13.9494 $  & $ 0.042547 $  & $ 13.9606 $  & $  0.016318 $  & $ 13.9748 $  & $ 0.051037 $  & $  14.035 $  & $ 13.9649 $  & $ 0.0015754 $  \\
$    10 $  & $ 9.8647 $  & $ 0.0202 $  & $ 9.8810 $  & $  0.0105 $  & $  9.8767 $  & $ 0.0482 $  & $ 9.9256 $  & $ 9.8678 $  & $ 0.0015 $  \\
$    30 $  & $ 8.9584 $  & $ 0.0169 $  & $ 8.9749 $  & $  0.0101 $  & $ 8.9723 $  & $ 0.0480 $  & $ 9.0204 $  & $ 8.9622 $  & $ 0.0015 $  \\
$   100 $  & $ 8.6643 $  & $ 0.0160 $  & $ 8.6808 $  & $ 0.0100 $  & $ 8.6785 $  & $ 0.0478 $  & $ 8.7264 $  & $ 8.6684 $  & $ 0.0014 $  \\
$   300 $  & $ 8.5832 $  & $ 0.0157 $  & $ 8.5997 $  & $ 0.0100 $  & $ 8.5967 $  & $ 0.0476 $  & $ 8.6443 $  & $ 8.5865 $  & $ 0.0014 $  \\
$1\mathord{,}000$ & $ 8.5540 $  & $ 0.0156 $  & $ 8.5705 $  & $ 0.0100 $  & $ 8.5683 $  & $ 0.0475 $  & $ 8.6157 $  & $ 8.5581 $  & $ 0.0014 $  \\
$10\mathord{,}000$  & $ 8.5437 $  & $ 0.0156 $  & $ 8.5603 $  & $ 0.0100 $  & $ 8.5573 $  & $ 0.0475 $  & $ 8.6048 $  & $ 8.5471 $  & $ 0.0014 $  \\
$\infty $  & $ 8.5427 $  & $ 0.0156 $  & $ 8.5593 $  & $ 0.0100 $  & $ 8.5561 $  & $ 0.0475 $  & $ 8.6036 $  & $ 8.5459 $  & $ 0.0014 $  \\
\end{tabular}
}
\caption{Bounds for $x_0=100$ (put option)}
\label{tab:singleuniput-1}
\end{center}
\end{table}

%\begin{figure}[h!]
%	\centering
%		\includegraphics[scale=0.63]{exbdput.pdf}
%	\caption{Exercise boundaries (put option; solid: $1/\delta_1=\infty$;
%dashed: $1/\delta_1=100$; and dotted: $1/\delta_1=10$).}
%	\label{fig:singleuniput}
%\end{figure}

Next, we consider a two-dimensional version of this example.
We suppose that there are two risky assets $X^{1}$ and $X^{2}$,
which are assumed to be independent and identically distributed
with the same dynamics as the single dividend-paying stock in the univariate case.
The payoff function we consider is a max-call, that is,
stopping at time $t$ yields a reward of $(\max(X_t^{1},X_t^{2})-K)^+$.
We set $X_0^{1}=X_0^{2}=100$ and $K=100$,
and allow for eleven equidistant exercise opportunities including $0$ and $T$, as before.
All other problem parameters and specifications remain the same.

Regarding the numbers of trajectories and time steps in the different stages of the algorithm,
we maintain the same specifications as in the univariate case.
The function basis is constructed as follows.
We always use the same set of 441 basis functions for $Y$ and for $(Z,\tilde{Z})$.
These consist of the constant $1$, 20 univariate basis functions $\phi^{(1)}_i(x^{1})$, $i=1,\ldots,20$,
which only depend on $X^{1}$,
20 univariate basis functions $\phi^{(2)}_j(x^{2})$, $j=1,\ldots,20$,
which only depend on $X^{2}$,
and all products $\phi^{(1)}_i(x^{1})\phi^{(2)}_{j}(x^{2})$, $i,j=1,\ldots,20$.
The 20 univariate basis functions are constructed as in the univariate case but with a slightly coarser grid, i.e.,
for $d=1,2$, we choose $x^{d}$ and $(x^{d}-q_{it}^{(d)})^+$ where the $q_{it}^{(d)}$ are the $5,10,\ldots, 95$ percent quantiles of $X^{d}_t$
estimated from the trajectories in the initial run of least squares Monte Carlo.
This bivariate basis of linear splines is relatively large but fairly generic, i.e.,
it exploits additional knowledge about the problem far less than, e.g., the $2d$-implementations in Belomestny, Bender and Schoenmakers \cite{BBS09} or Kr\"atschmer \textit{et al.} \cite{KLLSS18},
which rely on prices of European max-call options.
In principle, one could increase efficiency by including a variable selection step in the first regression.

We observe from Table \ref{tab:singlebicall-1} that the gaps between the genuine lower bound and the approximate upper bound are fairly small,
corresponding to about only 0.4\% of the option value.
The presence of ambiguity amplifies in the multivariate setting and its impact is more pronounced than in the univariate case.
%The tracking error in the forward simulation is somewhat larger, yielding about 5\% of the option value as the difference between the two genuine bounds.

\begin{table}[h!]
\begin{center}
{\small
\begin{tabular}{c||cc|cc|c%cc
|cc}
$1/\delta_1$ & LB without M & s.e. & LB & s.e. & $\overline{Y}_0^{N_{4}}$ %& TE & UB
& $\mathbb{E}[\overline{Y}_T^{N_{4}}]$ & s.e. \\
\hline\hline
10 &  16.5266  &  0.0351 &  16.5252 &   0.0196 &  16.5902 %&  0.7242 &  17.3253
& 16.5705   & 0.0230\\
30 &  14.7575   & 0.0268  & 14.7513  &  0.0172  & 14.8054  %&  0.6772  & 15.4837
& 14.7992  &  0.0216\\
100 &  14.1997   & 0.0246  & 14.1921  &  0.0167  & 14.2418  %&  0.6613  & 14.9032
& 14.2365  &  0.0211\\
300  & 14.0434    &0.0240   &14.0348   & 0.0166   &14.0861   %&  0.6576   &14.7438
&14.0807   & 0.0210\\
$1\mathord{,}000$  & 13.9916 &   0.0238   &13.9830   & 0.0165   &14.0322   %&  0.6564   &14.6886
&14.0268   & 0.0210\\
$10\mathord{,}000$  & 13.9700 &   0.0238   &13.9618   & 0.0165   &14.0115   %&  0.6559   &14.6674
&14.0061   & 0.0210\\
$\infty$&   13.9679 &   0.0237 &  13.9597 &   0.0165 &  14.0092 &   %0.6559 &  14.6650 &
14.0038 &   0.0210\\
\end{tabular}
}
\caption{Bounds for $x_0^{1}=x_0^{2}=100$ (bivariate case)}
\label{tab:singlebicall-1}
\end{center}
\end{table}

\subsection{Multiple Stopping: Swing Option in a Two-Factor Jump-Diffusion Model}

Supported by the accuracy and stability of our pathwise duality approach for optimal single stopping,
we now proceed to multiple stopping.
In this subsection, we analyze a canonical multiple stopping problem,
that of swing option pricing in electricity markets,
in the presence of ambiguity.
For this purpose,
we consider a two-factor jump-diffusion model for the electricity log-price process,
which has been suggested by Hambley, Howison and Kluge \cite{HHK09} to be a more realistic extension
of the one-factor Gaussian model proposed by Lucia and Schwartz \cite{LS02}
and implemented e.g., by Bender, Schoenmakers and Zhang \cite{BSZ15}.

Specifically, we assume that the electricity price $X_t$ at time $t>0$ is given by $X_{t}=X_{0} \exp(f(t)+u_t+v_t)$,
where the two stochastic factors $u_t$ and $v_t$ are mutually independent
and $f(t)$ is a deterministic function of time that can be used to calibrate the model.
The factor $u$ is Gaussian and follows the SDE
\begin{equation*}
du_t =-\kappa_u u_t +\sigma_u dW_t,
\end{equation*}
with $u_0=0$, $\kappa_u, \sigma_u >0$, and $W_t$ a standard Brownian motion.
The jump component $v$ follows the SDE
\begin{equation*}
dv_t =-\kappa_v v_{t_-} + J dN_t,
\end{equation*}
where $v_0=0$, $N$ is a (non-compensated) Poisson process with arrival rate $\lambda_{\mathbb{P}}$,
and $\kappa_v$ and the (deterministic) jump size $J$ are positive constants.
In the special case $\kappa_u\equiv\kappa_v$,
this model reduces to a mean-reverting one-factor jump-diffusion model for the log-price process.
For $\kappa_v \gg \kappa_u$, the model combines a mean-reverting Gaussian component, like in the Lucia-Schwartz model,
with occasional highly transitory spikes;
see Hambley, Howison and Kluge \cite{HHK09} for further discussion.
Using their formula (3), the processes $u$ and $v$ can be simulated forward in time without discretization error.

In our illustration, we consider a swing option contract that gives the owner the right to purchase electricity at a strike price $K$,
and consider $L$ exercise rights, in the time interval $[0,T]$.
We assume a fixed number of equidistant exercise opportunities.
We set the parameters of the price process as $S_0=10$, $f\equiv 0$, $\kappa_u=10$, $\sigma_u=0.25$, $\kappa_v=50$, $\lambda_{\mathbb{P}}=1$ and $J=\{0,0.06\}$.
Furthermore, we set the contract parameters as $K=10$, $T=5$ and consider 21 equidistant exercise opportunities
(including one at time 0 and one at time $T$).
We allow for varying degrees of ambiguity towards the Gaussian and the jump components of the price process
and consider $g(t,z,\tilde{z})=\delta_1 |z|+\delta_2 |\tilde{z}|$, $\delta_{1},\delta_{2}>0$;
cf. Example~\ref{g-ex}.

The overall numerical implementation is very similar to the previous optimal single stopping example.
The three BSDEs---$\rho$-martingale, lower, and upper bounds---are solved with two sets of simulations.
We consider one set of simulations with $100\mathord{,}000$ trajectories and $1\mathord{,}000$ time steps
for the initial $\rho$-martingale BSDE
and a second set of simulations with $100\mathord{,}000$ trajectories and $1\mathord{,}000$ time steps
for the BSDEs associated with the lower and upper bounds.
We choose the same basis functions as in the single stopping example and thus obtain 52 basis functions for $Y$ and 52 for $Z$,
and now also 52 basis functions for $\tilde{Z}$.
Note that these functions depend only on $X$ but not on $u$ and $v$.
For the evaluation of the lower bound, we draw a new sample of $100\mathord{,}000$ trajectories and $1\mathord{,}000$ time steps.
For the evaluation of the upper bound, we again artificially create a finer time discretization,
by repeating each set of coefficients 100 times, and reduce the number of trajectories to $1\mathord{,}000$.

Tables~\ref{tab:multionecall-1}--\ref{tab:multionecall-2} and \ref{tab:multitwocall-1}--\ref{tab:multitwocall-6} (in the Appendix)
summarize lower and upper bounds for different numbers of exercise rights,
different values of $\delta_{1}$ and $\delta_{2}$,
and for different values of $J$ (i.e., without and with jump component).
Upon comparing the results for one exercise right to those for multiple exercise rights,
we readily observe that the impact of ambiguity is even more pronounced in the multiple stopping case.

\begin{table}[h!]
\begin{center}
{\small
\begin{tabular}{c||cccc c} %ccc|}
$L$ & 1&2&3&4&5\\%&6&7&8\\
\hline\hline
%    LB without M & $  0.95154 $  & $    1.7006 $  & $   2.3155 $  & $     2.828 $  & $    3.2539 $ \\% & $   3.6057 $  & $    3.8901 $  & $    4.1164 $  \\
    LB        & $  0.9526 $  & $     1.7020 $  & $   2.3178 $  & $     2.8290 $  & $    3.2539 $ \\% & $   3.6053 $  & $    3.8899 $  & $    4.1161 $  \\
    s.e.      & $ 0.0012 $  & $ 0.0018 $  & $ 0.0023 $  & $ 0.0027 $  & $ 0.0030 $ \\% & $ 0.003396 $  & $ 0.0037152 $  & $ 0.0040009 $  \\
    $\overline{Y}_0^{N_{4},L}$  & $  0.9599 $  & $    1.7138 $  & $   2.3319 $  & $    2.8453 $  & $    3.2717 $ \\% & $   3.6228 $  & $    3.9072 $  & $    4.1321 $  \\
    TE        & $ 0.0315 $  & $  0.0407 $  & $ 0.0484 $  & $  0.0526 $  & $  0.0568 $ \\% & $ 0.059198 $  & $  0.060106 $  & $  0.060587 $  \\
    UB        & $  0.9914 $  & $    1.7546 $  & $   2.3803 $  & $    2.8979 $  & $    3.3284 $ \\% & $    3.682 $  & $    3.9674 $  & $    4.1927 $  \\
\end{tabular}
}
\caption{Bounds for $\delta_1=0$ and $J=0$} %(i.e., without jumps)}
\label{tab:multionecall-1}
%\end{table}
%\begin{table}[h!]
{\small
\vskip 0.1cm
\begin{tabular}{c||cccc c} %ccc|}
$L$ & 1&2&3&4&5\\%&6&7&8\\
\hline\hline
%    LB without M & $   0.95275 $  & $    1.7026 $  & $    2.3202 $  & $    2.8344 $  & $    3.2618 $ \\% & $   3.6157 $  & $  3.9047 $  & $   4.1324 $  \\
    LB        & $   0.9796 $  & $    1.7572 $  & $     2.4040 $  & $    2.9293 $  & $    3.3946 $ \\% & $    3.794 $  & $   4.112 $  & $   4.3438 $  \\
    s.e.      & $ 0.0027 $  & $ 0.0048 $  & $ 0.0065 $  & $ 0.0079 $  & $ 0.0092 $ \\% & $ 0.010468 $  & $ 0.01161 $  & $ 0.012084 $  \\
    $\overline{Y}_0^{N_{4},L}$  & $    1.0129 $  & $    1.8174 $  & $    2.4851 $  & $     3.0470 $  & $    3.5205 $ \\% & $   3.9168 $  & $  4.2447 $  & $   4.5105 $  \\
    TE        & $  0.0331 $  & $  0.0426 $  & $  0.0469 $  & $  0.0520 $  & $  0.0562 $ \\% & $ 0.058795 $  & $ 0.06229 $  & $ 0.065304 $  \\
    UB        & $    1.0495 $  & $    1.8644 $  & $     2.5370 $  & $    3.1045 $  & $    3.5826 $ \\% & $   3.9818 $  & $  4.3135 $  & $   4.5826 $  \\
\end{tabular}
}
\caption{Bounds for $\delta_1=0.2$ and $J=0$} %(i.e., without jumps)}
\label{tab:multionecall-2}
\end{center}
\end{table}

%\section{Conclusion}
%
%\label{sec:con}
%
%[***TO FOLLOW.***]

%\bibliographystyle{plain}
%\bibliography{LSSSref}

%[***SUPPRESS REFERENCES NOT CITED.***]

%\newpage

%MOR:
\noindent\textbf{Acknowledgements.} We are very grateful to
Christian Bender and to seminar and conference participants at the University of Waterloo and the Lorentz Center for comments and suggestions.
This research was funded in part by
the Netherlands Organization for Scientific Research (NWO) under grants NWO-Vidi and NWO-Vici (Laeven)
and by the DFG Excellence Cluster Math+ Berlin, project AA4-2 (Schoenmakers).

\newpage

{\Huge
\begin{center}
SUPPLEMENTARY MATERIAL TO\\
``Robust Multiple Stopping --- A Pathwise Duality Approach''\\
(FOR ONLINE PUBLICATION)
\end{center}
}

\newpage

\appendix

\setcounter{equation}{0}

\section*{ONLINE APPENDIX}

\setcounter{equation}{0}

\section{Proofs and Auxiliary Results for Section \ref{sec:robstop}}

\subsection{An Auxiliary Lemma on Sensitivity (P2) and the Proofs of Lemmas \ref{lem:Doob} and \ref{L3}}
\label{sec:A1}

We state the following auxiliary lemma:
\begin{lemma}
\label{lem:P2a}
%For subadditive and sensitive $\rho$ we have the following result:
%\begin{enumerate}
%\item[(P2a)]
(P2a) If subadditivity (P1) applies, then sensitivity (P2) of $\rho$ implies
\begin{equation}
\left[  X\geq0\text{ \ \ and \ \ }\rho_{t}(X)\leq0\right]  \Longrightarrow
X=0,\text{ \ \ for all }X\in\mathfrak{X,}\text{ and }t\in\left\{
0,\ldots,T\right\}  . \label{se1}%
\end{equation}
%\end{enumerate}
\end{lemma}
\begin{Proof}
Let $\rho$ be subadditive (P1) and sensitive (P2).
Suppose $\rho_{i}\left(  Y\right)  \leq0$ and $Y\geq0.$
Then, $-Y\leq0$ and so, by subadditivity,
\begin{equation*}
0=\rho_{t}\left(  Y-Y\right)  \leq\rho_{t}\left(  Y\right)  +\rho_{t}\left(
-Y\right)  \leq\rho_{t}\left(  -Y\right).
\end{equation*}
Hence, by (P2), $-Y=0,$ i.e., $Y=0$ a.s.
\end{Proof}

\vskip 0.3cm

\begin{Proof}
\textbf{of Lemma \ref{lem:Doob}.}
For $i=T$ the statement is trivial.
Assuming that it holds for $0<i\leq T,$ we have%
\begin{align*}
\rho_{i-1}\left(  M_{\tau_{i-1}}\right)   &  =\rho_{i-1}\left(  1_{\tau
_{i-1}=i-1}M_{i-1}+1_{\tau_{i-1}>i-1}M_{\tau_{i-1}\vee i}\right) \\
\text{(by (C2) and (C4))}  &  =1_{\tau_{i-1}=i-1}M_{i-1}+1_{\tau_{i-1}%
>i-1}\rho_{i-1}\left(  M_{\tau_{i-1}\vee i}\right) \\
\text{(by (C3))}  &  =1_{\tau_{i-1}=i-1}M_{i-1}+1_{\tau_{i-1}>i-1}\rho
_{i-1}\circ\rho_{i}\left(  M_{\tau_{i-1}\vee i}\right) \\
\text{(by induction)}  &  =1_{\tau_{i-1}=i-1}M_{i-1}+1_{\tau_{i-1}>i-1}%
\rho_{i-1}\left(  M_{i}\right) \\
\text{(property of }\rho\text{-martingale)}  &  =M_{i-1}.
\end{align*}
\end{Proof}

\vskip 0.3cm

\begin{Proof}
\textbf{of Lemma \ref{L3}.}
(i) For an arbitrary $X\in\mathfrak{X},$ and an
arbitrary set $A\in\mathcal{B}\left(\mathbb{R}\right),$
we have for any $i\geq t$,
\begin{align*}
\left\{  \rho_{\tau}(X)\in A\right\}  \cap\left\{  \tau=i\right\}   &
=\left\{  \sum_{j=t}^{T}1_{\tau=j}\rho_{j}(X)\in A\right\}  \cap\left\{
\tau=i\right\} \\
&  =\left\{  \rho_{i}(X)\in A\right\}  \cap\left\{  \tau=i\right\}
\in\mathcal{F}_{i},
\end{align*}
hence $\rho_{\tau_{t}}(X)\in\mathcal{F}_{\tau}.$\newline\vskip 0.01cm
\noindent(ii) Induction: For $t=T$ the statements are trivial.
Suppose they are true for $0<t\leq T.$
Now let $\tau\geq t-1$ and define $\tau_{1}:=$ $\tau\vee t.$
Then, %obviously,%
\begin{equation}
\rho_{\tau}(X)=1_{\tau=t-1}\rho_{t-1}(X)+1_{\tau>t-1}\rho_{\tau_{1}}(X).
\label{ta}%
\end{equation}
Thus, for any $X\in\mathfrak{X}$,
\begin{align*}
\rho_{t-1}(X) =  &  1_{\tau=t-1}\rho_{t-1}(X)+1_{\tau>t-1}\rho_{t-1}(X)\\
\text{(by (C3) and induction)} =  &  1_{\tau=t-1}\rho_{t-1}(X)
+1_{\tau>t-1}\rho_{t-1}\circ\rho_{t}\circ\rho_{\tau_{1}}(X)\\
\text{(by (C2) and (\ref{ta}))} =  &  1_{\tau=t-1}\rho_{t-1}(X)+1_{\tau
>t-1}\rho_{t-1}\circ\rho_{t}
\circ\left(  \rho_{\tau}(X)-1_{\tau=t-1}\rho_{t-1}(X)\right) \\
\text{(by (C4))} =  &  1_{\tau=t-1}\rho_{t-1}(X)
+ 1_{\tau>t-1}\left(  \rho_{t-1}\circ\rho_{t}\circ\rho_{\tau}%
(X)-1_{\tau=t-1}\rho_{t-1}(X)\right) \\
\text{(by (C3) and (\ref{ta}))} =  &  1_{\tau=t-1}\rho_{\tau}(X)+1_{\tau
>t-1}\rho_{t-1}\circ\rho_{\tau}(X)\\
\text{(by i) and (C6))} =  &  \rho_{t-1}\circ\rho_{\tau}(X).
\end{align*}
\vskip 0.01cm
\noindent(iii) Let $X\in\mathcal{F}_{\tau}$ and $Y\in\mathcal{F}_{T}.$
Then, $1_{\tau=t-1}X\in\mathcal{F}_{t-1}.$
Indeed, for any $A\in\mathcal{B}\left(\mathbb{R}\right)$ one has%
\begin{align*}
\left\{  1_{\tau=t-1}X\in A\right\}  =  &  \left(  \left\{  1_{\tau=t-1}X\in
A\right\}  \cap\left\{  \tau=t-1\right\}  \right)
\cup\left(  \left\{  0\in A\right\}  \cap\left\{  \tau>t-1\right\}  \right)
\\
=  &  \left(  \left\{  X\in A\right\}  \cap\left\{  \tau=t-1\right\}  \right)
\cup\left(  \left\{  0\in A\right\}  \cap\left\{  \tau>t-1\right\}
\right)  \in\mathcal{F}_{t-1}.
\end{align*}
Furthermore, also $X\in\mathcal{F}_{\tau_{1}}$ since $\tau_{1}\geq\tau.$
Hence, we have by (\ref{ta})
\begin{align*}
\rho_{\tau}(X+Y) =  &  1_{\tau=t-1}\rho_{t-1}(X+Y)+1_{\tau>t-1}\rho_{\tau_{1}%
}(X+Y)\\
\text{((C2) and induction)} =  &  1_{\tau=t-1}\rho_{t-1}(1_{\tau=t-1}%
X+1_{\tau=t-1}Y)
+1_{\tau>t-1}\left(  X+\rho_{\tau_{1}}(Y)\right) \\
\text{((C4) and above argument)} =  &  1_{\tau=t-1}\left(  1_{\tau
=t-1}X+1_{\tau=t-1}\rho_{t-1}(Y)\right)
+1_{\tau>t-1}\left(  X+\rho_{\tau_{1}}(Y)\right) \\
=  &  X+1_{\tau=t-1}\rho_{t-1}(Y)+1_{\tau>t-1}\rho_{\tau_{1}}(Y)\\
\text{(by (\ref{ta}))} =  &  X+\rho_{\tau}(Y).
\end{align*}
\end{Proof}

\subsection{Auxiliary Results on the Robust Single Optimal Stopping Problem \eqref{eq:OST}}
\label{sec:singlestop}

As is well-known, for the robust optimal single stopping problem,
we may find an optimal stopping family
$\left(\tau_{t}^{\ast}\right)_{t\in\{0,\ldots,T\}}$
satisfying
\begin{equation*}
Y_{t}^{\ast}=\sup_{\tau\in\mathcal{T}_{t}}\rho_{t}(H_{\tau})=\rho_{t}(H_{\tau_{t}^{\ast}}),\qquad t\in\{0,\ldots,T\},
\end{equation*}
and, furthermore, the Bellman principle
\begin{equation}
Y_{t}^{\ast}=\max\left(H_{t},\rho_{t}\left(Y_{t+1}^{\ast}\right)\right),\qquad t\in\{0,\ldots,T-1\},
\label{eq:Bel}%
\end{equation}
is satisfied (see e.g., Kr{\"a}tschmer and Schoenmakers \cite{KS10}
and Kr{\"a}tschmer \textit{et al.} \cite{KLLSS18} for details).

%\subsubsection{Additive dual representation}
Let us briefly recall the already existing (non-pathwise) additive dual representation
for the optimal single stopping problem \eqref{eq:OST}
(cf. Kr{\"a}tschmer and Schoenmakers \cite{KS10}
and Kr{\"a}tschmer \textit{et al.} \cite{KLLSS18}),
but with a different proof adapted to the goals in this paper and exploited later.

\begin{proposition}
\label{prop:ADP}
Let $\rho$ be a DMU satisfying (C1)--(C4) and
let $M^{\ast}=M^{\ast\rho}\in\mathcal{M}_{0}^{\rho}$ be the unique $\rho$-martingale in
the $\rho$-Doob decomposition of $Y^{\ast}=\left(Y_{t}^{\ast}\right)_{0\leq t\leq T}$.
Then the optimal single stopping problem \eqref{eq:OST} has an additive dual representation
\begin{align}
Y_{t}^{\ast}  &  =\inf_{M\in\mathcal{M}_{0}^{\rho}}\rho_{t}\,\Big(\max
_{j\in\{t,\ldots,T\}}\left(  H_{j}+M_{T}-M_{j}\right)  \Big)\label{dual}\\
&  =\rho_{t}\,\Big(\max_{j\in\{t,\ldots,T\}}\left(  H_{j}+M_{T}^{\ast}%
-M_{j}^{\ast}\right)  \Big),\text{ \ \ }t\in\{0,\ldots,T\}.\nonumber
\end{align}
\end{proposition}

\begin{Proof}
\textbf{of Proposition \ref{prop:ADP}.}
For any $\rho$-martingale $M$ and any stopping time $\tau\geq t,$
we have by Lemmas \ref{lem:Doob} and \ref{L3} that
\begin{align*}
\rho_{t}\,\Big(\max_{j\in\{t,\ldots,T\}}\left(  H_{j}+M_{T}-M_{j}\right)
\Big)  &  \geq\rho_{t}\,\left(  H_{\tau}+M_{T}-M_{\tau}\right)
=\rho_{t}\,\circ\rho_{\tau}\left(  H_{\tau}+M_{T}-M_{\tau}\right) \\
&  =\rho_{t}\,\left(  H_{\tau}-M_{\tau}+\rho_{\tau}\left(  M_{T}\right)
\right)  =\rho_{t}\,\left(  H_{\tau}\right),
\end{align*}
which implies%
\begin{equation*}
Y_{t}^{\ast}\leq\inf_{M\in\mathcal{M}_{0}^{\rho}}\rho_{t}\,\Big(\max
_{j\in\{t,\ldots,T\}}\left(  H_{j}+M_{T}-M_{j}\right)  \Big).
\end{equation*}
On the other hand, for the $\rho$-Doob martingale $M^{\ast}$ it holds that
\begin{align*}
H_{j}+M_{t}^{\ast}-M_{j}^{\ast}  &  =H_{j}+\sum_{r=t}^{j-1}M_{r}^{\ast
}-M_{r+1}^{\ast}\\
\text{(by (\ref{eq:DD}))}  &  =H_{j}+\sum_{r=t}^{j-1}\rho_{r}\left(
Y_{r+1}^{\ast}\right)  -Y_{r+1}^{\ast}\\
\text{(Bellman)}  &  \leq H_{j}+\sum_{r=t}^{j-1}Y_{r}^{\ast}-Y_{r+1}^{\ast
}=Y_{t}^{\ast}+H_{j}-Y_{j}^{\ast}\leq Y_{t}^{\ast},
\end{align*}
whence%
\begin{align*}
\rho_{t}\,\Big(\max_{j\in\{t,\ldots,T\}}\left(  H_{j}+M_{T}^{\ast}-M_{j}%
^{\ast}\right)  \Big)  &  =\rho_{t}\,\Big(\max_{j\in\{t,\ldots,T\}}\left(
H_{j}+M_{t}^{\ast}-M_{j}^{\ast}\right)  +M_{T}^{\ast}-M_{t}^{\ast}\Big)\\
&  \leq Y_{t}^{\ast}+\rho_{t}\,\Big(M_{T}^{\ast}-M_{t}^{\ast}\Big)=Y_{t}%
^{\ast}.
\end{align*}
\end{Proof}

\setcounter{equation}{0}

\section{Proofs and Auxiliary Results for Section \ref{sec:pathdual}}

\subsection{Proofs and Auxiliary Results for Section \ref{sec:pathdualrep}}

We state the following lemma:
\begin{lemma}
\label{L2}
Suppose that $\rho$ satisfies (C1)--(C4) and (P1).
%, hence in particular $\rho$ is subadditive.
Then, for any adapted process $H$,
any $\rho$-martingale $M$,
and any stopping $\tau$, with $T\geq\tau\geq t$ a.s. it holds that
\begin{equation*}
\rho_{t}(H_{\tau})\leq\rho_{t}(H_{\tau}+M_{t}-M_{\tau}),\qquad 0\leq t\leq T.
\end{equation*}
\end{lemma}

\begin{Proof}
\textbf{of Lemma \ref{L2}.}
Using Lemma \ref{L3} and the proof of Proposition \ref{prop:ADP} one has
\begin{align*}
\rho_{t}(H_{\tau})  &  =\rho_{t}(H_{\tau}+M_{T}-M_{\tau})\\
&  =\rho_{t}(H_{\tau}+M_{t}-M_{\tau}+M_{T}-M_{t})\\
\text{(by (P1))}  &  \leq\rho_{t}\left(  H_{\tau}+M_{t}-M_{\tau}\right)
+\rho_{t}(M_{T}-M_{t})\\
\text{((C4) and }\rho\text{-mart. prop.)}  &  =\rho_{t}\left(H_{\tau}+M_{t}-M_{\tau}\right).
\end{align*}
\end{Proof}

\vskip 0.3cm

%\subsection{Proof of Theorem \ref{th:mulsub}}
\begin{Proof}
\textbf{of Theorem \ref{th:mulsub}.}
First we show that for any $t=0,\ldots,T$,
any sequence of stopping times with $t\leq\tau_{1}<\tau_{2}<\cdots<\tau_{L}$ a.s.,
and any set of $\rho$-martingales $M^{(1)},\ldots,M^{(L)}$,
we have
\begin{equation}
\rho_{t}\left(  \sum_{k=1}^{L}H_{\tau_{k}}\right)  \leq\rho_{t}\left(
\sum_{k=1}^{L}\left(H_{\tau_{k}}+M_{\tau_{k-1}}^{(k)}-M_{\tau_{k}}^{(k)}\right)\right).
\label{fi}%
\end{equation}
For $L=1$, this statement boils down to Lemma \ref{L2}.
Let us assume the statement is true for some $L\geq1$.
Take $0\leq t\leq T$ and $t\leq\tau_{1}<\tau_{2}<\cdots<\tau_{L+1}$ arbitrarily.
Observe that (with $\rho_{T+1}:=\rho_{T}$)
\begin{align}
&  \rho_{\tau_{1}}\left(  \sum_{k=2}^{L+1}\left(H_{\tau_{k}}+M_{\tau_{k-1}}%
^{(k)}-M_{\tau_{k}}^{(k)}\right)\right) \label{ind}\\
&  =\sum_{j=t}^{T}1_{\tau_{1}=j}\rho_{j}\circ\rho_{j+1}\left(  1_{j+1\leq
\tau_{2}<\cdots<\tau_{L+1}}\sum_{k=2}^{L+1}\left(H_{\tau_{k}}+M_{\tau_{k-1}}%
^{(k)}-M_{\tau_{k}}^{(k)}\right)\right) \nonumber\\
\text{(by induction)}  &  \geq\sum_{j=\tau_{1}}^{T}1_{\tau_{1}=j}\rho
_{j}\left(  1_{j<\tau_{2}<\cdots<\tau_{L+1}}\sum_{k=2}^{L+1}H_{\tau_{k}%
}\right)  =\rho_{\tau_{1}}\left(  \sum_{k=2}^{L+1}H_{\tau_{k}}\right)
.\nonumber
\end{align}
One may thus write, by Lemma \ref{L3},
\begin{align*}
&  \rho_{t}\left(  \sum_{k=1}^{L+1}\left(H_{\tau_{k}}+M_{\tau_{k-1}}^{(k)}%
-M_{\tau_{k}}^{(k)}\right)\right) \\
&  =\rho_{t}\left(  H_{\tau_{1}}+M_{t}^{(1)}-M_{\tau_{1}}^{(1)}+\rho_{\tau
_{1}}\left(  \sum_{k=2}^{L+1}\left(H_{\tau_{k}}+M_{\tau_{k-1}}^{(k)}-M_{\tau_{k}%
}^{(k)}\right)\right)  \right) \\
\text{(by Lemma \ref{L2})}  &  \geq\rho_{t}\left(  H_{\tau_{1}}+\rho_{\tau
_{1}}\left(  \sum_{k=2}^{L+1}\left(H_{\tau_{k}}+M_{\tau_{k-1}}^{(k)}-M_{\tau_{k}%
}^{(k)}\right)\right)  \right) \\
\text{(by (\ref{ind}))}  &  \geq\rho_{t}\left(  H_{\tau_{1}}+\rho_{\tau_{1}%
}\left(  \sum_{k=2}^{L+1}H_{\tau_{k}}\right)  \right) \\
\text{(Lemma \ref{L3})}  &  =\rho_{t}\left(  \sum_{k=1}^{L+1}H_{\tau_{k}%
}\right),
\end{align*}
which proves (\ref{fi}).
As a corollary, we obtain
\begin{equation}
Y_{t}^{\ast,L}\leq\rho_{t}\left(  \max_{t\leq j_{1}<j_{2}<\cdots<j_{L}}%
\sum_{k=1}^{L}\left(H_{j_{k}}+M_{j_{k-1}}^{(k)}-M_{j_{k}}^{(k)}\right)\right),
\label{fi1}%
\end{equation}
where we note that for any set $A$ of probability one has%
\begin{equation*}
1_{A}\rho_{t}\left(  X\right)  =\rho_{t}\left(  1_{A}X\right)  =\rho
_{t}\left(  X\right),
\end{equation*}
due to monotonicity (C1).
Since the $\rho$-martingales $M^{(k)}$ are arbitrary,
we thus arrive at
\begin{equation}
Y_{t}^{\ast,L}\leq\underset{M^{(1)},\ldots,M^{(L)}\in\mathcal{M}_{0}^{\rho}}{\inf}\rho_{t}\left(
\max_{t\leq j_{1}<j_{2}<\cdots<j_{L}}\sum_{k=1}^{L}\left(H_{j_{k}}+M_{j_{k-1}}%
^{(k)}-M_{j_{k}}^{(k)}\right)\right).
\label{fi2}%
\end{equation}
On the other hand, for the $\rho$-Doob martingales $M^{\ast,L-k+1}$, we may write
(with $j_{0}=t$)
\begin{align*}
&\sum_{k=1}^{L}\left(  H_{j_{k}}+M_{j_{k-1}}^{\ast,L-k+1}-M_{j_{k}}%
^{\ast,L-k+1}\right)  \\
& =\sum_{k=1}^{L}\left(  H_{j_{k}}+\sum_{r=j_{k-1}%
}^{j_{k}-1}\left(M_{r}^{\ast,L-k+1}-M_{r+1}^{\ast,L-k+1}\right)\right) \\
&  =\sum_{k=1}^{L}H_{j_{k}}+\sum_{k=1}^{L}\sum_{r=j_{k-1}}^{j_{k}-1}\left(\rho
_{r}\left(  Y_{r+1}^{\ast,L-k+1}\right)  -Y_{r+1}^{\ast,L-k+1}\right)\\
&  =\sum_{k=1}^{L}H_{j_{k}}+\sum_{k=1}^{L}\sum_{r=j_{k-1}}^{j_{k}-1}\left(
Y_{r}^{\ast,L-k+1}-Y_{r+1}^{\ast,L-k+1}\right) \\
&  \quad
+\sum_{k=1}^{L}\sum_{r=j_{k-1}}^{j_{k}-1}\left(  \rho_{r}\left(
Y_{r+1}^{\ast,L-k+1}\right)  -Y_{r}^{\ast,L-k+1}\right) \\
&  =\sum_{k=1}^{L}H_{j_{k}}+\sum_{k=1}^{L}\left(  Y_{j_{k-1}}^{\ast
,L-k+1}-Y_{j_{k}}^{\ast,L-k+1}\right) \\
&  \quad
+\sum_{k=1}^{L}\sum_{r=j_{k-1}}^{j_{k}-1}\left(  \rho_{r}\left(
Y_{r+1}^{\ast,L-k+1}\right)  -Y_{r}^{\ast,L-k+1}\right) \\
&  =Y_{j_{0}}^{\ast,L}+\underset{\leq0}{\underbrace{H_{j_{L}}-Y_{j_{L}}%
^{\ast,1}}}+\underset{\leq0}{\sum_{k=1}^{L-1}\underbrace{\left(  H_{j_{k}%
}+Y_{j_{k}}^{\ast,L-k}-Y_{j_{k}}^{\ast,L-k+1}\right)  }}\\
&  \quad+\sum_{k=1}^{L}\sum_{r=j_{k-1}}^{j_{k}-1}\underset{\leq0}%
{\underbrace{\left(  \rho_{r}\left(  Y_{r+1}^{\ast,L-k+1}\right)  -Y_{r}%
^{\ast,L-k+1}\right)  }}\text{ }\leq Y_{j_{0}}^{\ast,L}.
\end{align*}
That is,
\begin{equation*}
\max_{t\leq j_{1}<j_{2}<\cdots<j_{L}}\sum_{k=1}^{L}\left(H_{j_{k}}+M_{j_{k-1}}%
^{\ast,L-k+1}-M_{j_{k}}^{\ast,L-k+1}\right)\leq Y_{t}^{\ast,L},
\end{equation*}
while, due to (\ref{fi1}),
\begin{equation*}
\rho_{t}\left(  \max_{t\leq j_{1}<j_{2}<\cdots<j_{L}}\sum_{k=1}^{L}\left(H_{j_{k}%
}+M_{j_{k-1}}^{\ast,L-k+1}-M_{j_{k}}^{\ast,L-k+1}\right)\right)  \geq Y_{t}^{\ast
,L}.
\end{equation*}
Thus, by monotonicity (C1) and $\mathcal{F}_{t}$-invariance (C6) we obtain (ii),
and, by sensitivity (P2), we obtain (iii).
Finally, (ii) combined with (\ref{fi2}) yields (i).
\end{Proof}

\subsection{Proofs and Auxiliary Results for Section \ref{sec:suroptmart}}

\begin{Proof}
\textbf{of Theorem \ref{th:sure1}.}
Suppose that $\theta_{i}:=\max_{i\leq j\leq T}(H_{j}-M_{j}+M_{i})\in\mathcal{F}_{i}$
and define the stopping time
\begin{equation*}
\tau_{i}:=\inf\{j\geq i:H_{j}-M_{j}+M_{i}\geq\theta_{i}\}.
\end{equation*}
By the definition of $\theta_{i}$, clearly $i\leq\tau_{i}\leq T$ a.s.
Hence, we have on the one hand
\begin{equation*}
Y_{i}^{\ast}\geq\rho_{i}(H_{\tau_{i}})\geq\rho_{i}(M_{\tau_{i}}-M_{i}+\theta_{i})=\theta_{i},
\end{equation*}
by the fact that $-M_{i}+\theta_{i}\in\mathcal{F}_{i}$,
translation invariance (C4), and Lemma~\ref{lem:Doob}.
On the other hand, we have $\theta_{i}=\rho_{i}(\theta_{i})\geq Y_{i}^{\ast}$
due to Theorem~\ref{P1}, Eqn.~(\ref{as0}).
\end{Proof}

\vskip 0.3cm

\begin{Proof}
\textbf{of Lemma~\ref{th:sure1_cor}.}
By writing
\begin{equation}
\theta_{i+}=\underset{\in\mathcal{F}_{i+1}}{\underbrace{\max_{i+1\leq j\leq
T}(H_{j}-M_{j}+M_{i+1})}}+M_{i}-M_{i+1},
\label{iplus}
\end{equation}
and applying Theorem~\ref{th:sure1},
we have
\begin{equation}
\theta_{i+}+M_{i+1}-M_{i}=Y_{i+1}^{\ast}.
\label{eq:Corro}
\end{equation}
Then, (i) follows by applying $\rho_{i}$ on both sides, using
conditional translation invariance (C4) and the martingale property.
Next, (ii) is obvious from (\ref{eq:Corro}).
\end{Proof}

\vskip 0.3cm

\begin{Proof}
\textbf{of Proposition \ref{prop:ch}.}
It is sufficient to show that
\begin{equation}
\rho_{i}\left(  1_{\left\vert Y\right\vert \geq\epsilon}\right)  \leq
\frac{\rho_{i}\left(  Y^{2}\right)  }{\epsilon^{2}}.
\label{sen}%
\end{equation}
Indeed, one has by monotonicity and positive homogeneity,
\begin{equation*}
\rho_{i}\left(  Y^{2}\right)  =\rho_{i}\left(  Y^{2}1_{\left\vert Y\right\vert
\geq\epsilon}+Y^{2}1_{\left\vert Y\right\vert <\epsilon}\right)  \geq\rho
_{i}\left(  Y^{2}1_{\left\vert Y\right\vert \geq\epsilon}\right)  \geq\rho
_{i}\left(  \epsilon^{2}1_{\left\vert Y\right\vert \geq\epsilon}\right)
=\epsilon^{2}\rho_{i}\left(  1_{\left\vert Y\right\vert \geq\epsilon}\right).
\end{equation*}
\end{Proof}

\vskip 0.3cm

\begin{Proof}
\textbf{of Lemma \ref{sensv}.}
Indeed, $\mathrm{Var}_{\rho_{i}}\left(  X\right)  =\rho_{i}\left(  \left(
X-\rho_{i}\left(  X\right)  \right)  ^{2}\right)  =0$ implies,
by (\ref{se1}),
$X-\rho_{i}\left(  X\right)=0$,
hence $X\in\mathcal{F}_{i}$.
The reverse direction is trivial.
\end{Proof}

\vskip 0.3cm

\begin{Proof}
\textbf{of Theorem \ref{th:Th14}.}
Suppose that the assumptions of the theorem
are satisfied.
Fix an $i\in{\{0,\dots,T\}}$ and take an $\epsilon>0$.
Upon introducing an auxiliary time $\partial>T$ and setting $H_{\partial}=0$,
we next define the stopping time
$\tau_{i}^{(n)}=\inf\{j\geq i:H_{j}-M_{j}^{(n)}+M_{i}^{(n)}\geq\rho
_{i}(\theta_{i}^{(n)})-\epsilon\}\wedge\partial$.

Then, with $M_{\partial}^{(n)}:=M_{T}^{(n)}$, $n\geq1$,
\begin{align*}
Y_{i}^{\ast}  &  \geq\rho_{i}(H_{\tau_{i}^{(n)}})
%\\
%&
=\rho_{i}(H_{\tau_{i}^{(n)}}1_{\tau_{i}^{(n)}<\partial})\\
&  \geq\rho_{i}\Big((M_{\tau_{i}^{(n)}}^{(n)}-M_{i}^{(n)}+\rho_{i}(\theta
_{i}^{(n)})-\epsilon)1_{\tau_{i}^{(n)}<\partial}\Big)\\
&  \geq\rho_{i}\Big(M_{\tau_{i}^{(n)}}^{(n)}-M_{i}^{(n)}+\rho_{i}(\theta
^{(n)})-\epsilon\Big)%\\
%&\qquad
-\rho_{i}\Big((M_{T}^{(n)}-M_{i}^{(n)}+\rho_{i}(\theta_{i}^{(n)}%
)-\epsilon)1_{\tau_{i}^{(n)}=\partial}\Big)\\
&  =\rho_{i}(\theta_{i}^{(n)})-\epsilon-\rho_{i}\Big((M_{T}^{(n)}-M_{i}%
^{(n)}+\rho_{i}(\theta_{i}^{(n)})-\epsilon)1_{\tau_{i}^{(n)}=\partial
}\Big),\qquad \mathrm{almost\ surely},
\end{align*}
using subadditivity in the last inequality and translation invariance in the last equality.
Hence,
\begin{align*}
\rho_{i}\left(  \theta_{i}^{(n)}\right)   &  \leq Y_{i}^{\ast}+\epsilon
+\rho_{i}\left(  \left\vert M_{T}^{(n)}-M_{i}^{(n)}+\rho_{i}\left(  \theta
_{i}^{(n)}\right)  -\epsilon\right\vert 1_{\tau_{i}^{(n)}=\partial}\right) \\
&  =:Y_{i}^{\ast}+\epsilon+\rho_{i}\left(  \left\vert U_{i}^{(n)}\right\vert
1_{\tau_{i}^{(n)}=\partial}\right),\qquad \mathrm{almost\ surely}.
\end{align*}

By (\ref{ch}),
\begin{equation*}
\rho_{i}\left(  1_{\tau_{i}^{(n)}=\partial}\right)  =\rho_{i}\left(
1_{\left\vert \theta_{i}^{(n)}-\rho_{i}\left(  \theta_{i}^{(n)}\right)
\right\vert \geq\epsilon}\right)  \leq\frac{\mathrm{Var}_{\rho_{i}}\left(
\theta_{i}^{(n)}\right)  }{\epsilon^{2}}\overset{\text{P}}{\rightarrow}0,
\end{equation*}
and since moreover by monotonicity $0\leq\rho_{i}\left(  1_{\tau_{i}%
^{(n)}=\partial}\right)  \leq\rho_{i}\left(  1\right)  =1,$ it holds that%
\begin{equation}
\rho_{i}\left(  1_{\tau_{i}^{(n)}=\partial}\right)  \overset{L_{1}%
}{\rightarrow}0. \label{l1a}%
\end{equation}
Next, by subadditivity, monotonicity, and positive homogeneity, we have, for any
$K>0$,
\begin{gather*}
\mathbb{E}\rho_{i}\left(  \left\vert U_{i}^{(n)}\right\vert 1_{\tau_{i}%
^{(n)}=\partial}\right)  \leq\mathbb{E}\rho_{i}\left(  \left\vert U_{i}%
^{(n)}\right\vert 1_{\tau_{i}^{(n)}=\partial}1_{\left\vert U_{i}%
^{(n)}\right\vert \leq K}\right)  +\mathbb{E}\rho_{i}\left(  \left\vert
U_{i}^{(n)}\right\vert 1_{\tau_{i}^{(n)}=\partial}1_{\left\vert U_{i}%
^{(n)}\right\vert >K}\right) \\
\leq K\mathbb{E}\rho_{i}\left(  1_{\tau_{i}^{(n)}=\partial}\right)
+\mathbb{E}\rho_{i}\left(  \left\vert U_{i}^{(n)}\right\vert 1_{\left\vert
U_{i}^{(n)}\right\vert >K}\right)  .
\end{gather*}
Now Propositions~\ref{claim*} and~\ref{claim**} below imply that the family
$\left(U_{i}^{(n)}\right)_{n\geq1}$ is also uniformly integrable in the
sense of (\ref{gUI}), i.e., there exists $K_{1,\epsilon}$ large enough such
that%
\begin{equation*}
\sup_{n\geq1}\mathbb{E}\rho_{i}\left(  \left\vert U_{i}^{(n)}\right\vert
1_{\left\vert U_{i}^{(n)}\right\vert >K}\right)  <\epsilon,
\end{equation*}
hence%
\begin{equation*}
\mathbb{E}\rho_{i}\left(  \left\vert U_{i}^{(n)}\right\vert 1_{\tau_{i}%
^{(n)}=\partial}\right)  \leq K_{\epsilon}\underset{\rightarrow0\text{ by
(\ref{l1a})}}{\underbrace{\mathbb{E}\rho_{i}\left(  1_{\tau_{i}^{(n)}%
=\partial}\right)  }}+\epsilon\leq2\epsilon,
\end{equation*}
for $n>N_{K_{1,\epsilon},\epsilon}$.
Thus, since $\epsilon>0$ was arbitrary,%
\begin{equation*}
\overline{\lim}_{n\geq1}\mathbb{E}\rho_{i}\left(  \theta_{i}^{(n)}\right)
\leq\mathbb{E}Y_{i}^{\ast}+3\epsilon.
\end{equation*}

On the other hand, by monotonicity and the duality theorem for subadditive
functionals,
\begin{equation*}
\mathbb{E}\rho_{i}\left(  \theta_{i}^{(n)}\right)  \geq\mathbb{E}Y_{i}^{\ast},
\end{equation*}
so it follows that%
\begin{equation*}
\lim_{n\rightarrow\infty}\mathbb{E}\rho_{i}\left(  \theta_{i}^{(n)}\right)
=\mathbb{E}Y_{i}^{\ast}.
\end{equation*}
\end{Proof}

\begin{proposition}
\label{claim*} Suppose $\left(  A_{n}\right)  _{n\geq1},$ and $\left(
B_{n}\right)  _{n\geq1}$ satisfy (\ref{gUI}), i.e.,%
\begin{equation*}
\sup_{n\geq1}\mathbb{E}\rho_{i}\left(  \left\vert A_{n}\right\vert
1_{\left\vert A_{n}\right\vert >K_{\epsilon}}\right)  <\epsilon\text{ \ \ and
\ \ }\sup_{n\geq1}\mathbb{E}\rho_{i}\left(  \left\vert B_{n}\right\vert
1_{\left\vert B_{n}\right\vert >K_{\epsilon}}\right)  <\epsilon,\text{ }%
\end{equation*}
for $K_{\epsilon}$ large enough.
If $\rho_{i}$ is subadditive and positively homogeneous, then also $\left(  A_{n}+B_{n}\right)  _{n\geq1}$ satisfies
(\ref{gUI}).
\end{proposition}

\begin{Proof}
\textbf{of Proposition \ref{claim*}.}
By (P1),
\begin{align*}
\mathbb{E}\rho_{i}\left(  \left\vert A_{n}+B_{n}\right\vert 1_{\left\vert
A_{n}+B_{n}\right\vert >2K_{\epsilon}}\right)  &\leq
\mathbb{E}\rho_{i}\left(  \left(  \left\vert A_{n}\right\vert +\left\vert
B_{n}\right\vert \right)  1_{\left\vert A_{n}\right\vert +\left\vert
B_{n}\right\vert >2K_{\epsilon}}\right) \\
&\leq\mathbb{E}\rho_{i}\left(  2\left\vert A_{n}\right\vert 1_{\left\vert
A_{n}\right\vert >K_{\epsilon}}+2\left\vert B_{n}\right\vert 1_{\left\vert
B_{n}\right\vert >K_{\epsilon}}\right) \\
&\leq2\mathbb{E}\rho_{i}\left(  \left\vert A_{n}\right\vert 1_{\left\vert
A_{n}\right\vert >K_{\epsilon}}\right)  +2\mathbb{E}\rho_{i}\left(  \left\vert
B_{n}\right\vert 1_{\left\vert B_{n}\right\vert >K_{\epsilon}}\right)
<4\epsilon,
\end{align*}
hence $A_{n}+B_{n}$ satisfies (\ref{gUI}) also.
\end{Proof}

\begin{lemma}
\label{lem:inter}
Assume (P1) and (P3).
$\left(  A_{n}\right)_{n\geq1}$ (with w.l.o.g. $A_{n}\geq0$) satisfy
(\ref{gUI}) if and only if
\begin{itemize}
\item[(i)] $\sup_{n\geq1}\rho_{i}\left(  A_{n}\right)  <\infty$;
\item[(ii)] For every $\epsilon>0$ there exists $\delta>0$ such that for all
$B\in\mathcal{F}$ with $\rho_{i}\left(  1_{B}\right)  <\delta,$ it holds that
$\sup_{n\geq1}\rho_{i}\left(  A_{n}1_{B}\right)  <\epsilon.$
\end{itemize}
\end{lemma}

\begin{Proof}
\textbf{of Lemma \ref{lem:inter}.}
($\Longrightarrow$) Let $\left(  A_{n}\right)  _{n\geq1}$ (with w.l.o.g.
$A_{n}\geq0$) satisfy (\ref{gUI}).
Then, for any $n\geq1,$ by subadditivity,
monotonicity, and positive homogeneity,
\begin{equation*}
\rho_{i}\left(  A_{n}\right)  \leq\rho_{i}\left(  A_{n}1_{A_{n}\leq K}\right)
+\rho_{i}\left(  A_{n}1_{A_{n}>K}\right)  \leq K\rho_{i}\left(  1_{A_{n}\leq
K}\right)  +1\leq K+1,
\end{equation*}
for large enough $K.$\ So $\sup_{n\geq1}\rho_{i}\left(  A_{n}\right)  \leq
K+1, $ whence (i).
Now let $\epsilon>0$ and $K$ be so large that%
\begin{equation*}
\sup_{n\geq1}\rho_{i}\left(  A_{n}1_{A_{n}>K}\right)  <\epsilon/2.
\end{equation*}
For any $B\in\mathcal{F}$ with $\rho_{i}\left(  1_{B}\right)  <\epsilon
/(2K)=:\delta$ we then have%
\begin{align*}
\rho_{i}\left(  A_{n}1_{B}\right)   &  \leq\rho_{i}\left(  A_{n}1_{B}%
1_{A_{n}\leq K}\right)  +\rho_{i}\left(  A_{n}1_{B}1_{A_{n}>K}\right)
%\\
%&
\leq K\rho_{i}\left(  1_{B}\right)  +\rho_{i}\left(  A_{n}1_{A_{n}%
>K}\right)  <\epsilon.
\end{align*}

($\Longleftarrow$) Let $\left(  A_{n}\right)  _{n\geq1}$ satisfy (i) and (ii)
for $\epsilon>0$ and $\delta>0.$
For any $n\geq1$ we have%
\begin{equation*}
\rho_{i}\left(  A_{n}\right)  \geq\rho_{i}\left(  A_{n}1_{A_{n}>K}\right)
\geq K\rho_{i}\left(  1_{A_{n}>K}\right)  ,
\end{equation*}
so due to (i),%
\begin{equation*}
M:=\sup_{n\geq1}\rho_{i}\left(  A_{n}\right)  \geq K\sup_{n\geq1}\rho
_{i}\left(  1_{A_{n}>K}\right).
\end{equation*}
Hence,%
\begin{equation*}
\sup_{n\geq1}\rho_{i}\left(  1_{A_{n}>K}\right)  \leq\frac{M}{K}<\delta,
\end{equation*}
if $K>M/\delta$.
Thus, due to (ii), for all $n\geq1,$ and $K>M/\delta,$
\begin{equation*}
\rho_{i}\left(  A_{n}1_{A_{n}>K}\right)  <\epsilon.
\end{equation*}
\end{Proof}

\begin{proposition}
\label{claim**} Let $\rho_{i}$ be subadditive and positively homogeneous, and
let $\left(  A_{n}\right)  _{n\geq1}$ satisfy (\ref{gUI}).
Then $\left(\rho_{i}\left(  A_{n}\right)  \right)_{n\geq1}$ also satisfy (\ref{gUI}).
\end{proposition}

\begin{Proof}
\textbf{of Proposition \ref{claim**}.}
Due to Lemma~\ref{lem:inter}, (i) and (ii) apply
for $\left(  A_{n}\right)_{n\geq1}.$
Let $\epsilon>0$ and take $\delta>0$ such that (ii) holds for
$\left(A_{n}\right)_{n\geq1}.$
Observe that%
\begin{align*}
\rho_{i}\left(  A_{n}\right)   &  =\rho_{i}\left(  \rho_{i}\left(
A_{n}\right)  \right)  \geq\rho_{i}\left(  \rho_{i}\left(  A_{n}\right)
1_{\rho_{i}\left(  A_{n}\right)  >K}\right) %\\
%&
\geq K\rho_{i}\left(  1_{\rho_{i}\left(  A_{n}\right)  >K}\right).
\end{align*}
Hence,
\begin{equation*}
\sup_{n\geq1}\rho_{i}\left(  1_{\rho_{i}\left(  A_{n}\right)  >K}\right)
\leq\frac{1}{K}\sup_{n\geq1}\rho_{i}\left(  A_{n}\right)  =:\frac{M}{K}.
\end{equation*}
Take $K$ such that $M/K<\delta.$
Then, for all $n\geq1,$%
\begin{equation*}
\rho_{i}\left(  \rho_{i}\left(  A_{n}\right)  1_{\rho_{i}\left(  A_{n}\right)
>K}\right)  =\rho_{i}\left(  A_{n}1_{\rho_{i}\left(  A_{n}\right)  >K}\right)
<\epsilon,
\end{equation*}
since $\rho_{i}\left(  1_{\rho_{i}\left(  A_{n}\right)  >K}\right)  <\delta.$
That is, $\left(  \rho_{i}\left(  A_{n}\right)  \right)  _{n\geq1}$ satisfy
(\ref{gUI}).
\end{Proof}

\vskip 0.3cm

\begin{Proof}
\textbf{of Proposition \ref{prop:EM}.}
We have
\begin{gather*}
\sup_{n\geq1}\mathbb{E}\rho_{i}\left(  \left\vert M_{i}^{(n)}\right\vert
1_{\left\vert M_{i}^{(n)}\right\vert >K}\right)  =\sup_{n\geq1}\mathbb{E}%
\rho_{i}\left(  \frac{1}{\left\vert M_{i}^{(n)}\right\vert ^{\eta}}\left\vert
M_{i}^{(n)}\right\vert ^{1+\eta}1_{\left\vert M_{i}^{(n)}\right\vert
>K}\right) \\
\leq\frac{1}{K^{\eta}}\sup_{n\geq1}\mathbb{E}\rho_{i}\left(  \left\vert
M_{i}^{(n)}\right\vert ^{1+\eta}\right)  \rightarrow0\text{ \ for
\ }K\rightarrow\infty.
\end{gather*}
\end{Proof}

\vskip 0.3cm

\begin{Proof}
\textbf{of Theorem \ref{th:sure1_cor_mult}.}
For $L=1$, this follows from Lemma~\ref{th:sure1_cor}.
Now let us suppose that
\begin{equation*}
\Theta_{i+}^{q}\in\mathcal{F}_{i},\text{ \ \ for \ \ }q=1,\ldots,L+1,\qquad 0\leq i<T,
\end{equation*}
and that the theorem has been proved for $L\geq1$.
Then, by induction, we have $(i)$
and $(ii)$, and so, with $j_{0}^{\prime}=j_{1}$,
\begin{align*}
\Theta_{i+}^{L+1}  &  =\max_{i<j_{1}<j_{2}<\cdots<j_{L+1}}\sum
_{k=1}^{L+1}\left(H_{j_{k}}+M_{j_{k-1}}^{(L+2-k)}-M_{j_{k}}^{(L+2-k)}\right)\\
&  =\max_{i<j_{1}}\Bigg(  H_{j_{1}}+M_{i}^{(L+1)}-M_{j_{1}}^{(L+1)}%\right.
\\
%&\hspace{0.4cm}  \left.
&\quad +\max_{j_{1}<j_{2}<\cdots<j_{L+1}}\sum_{k=2}^{L+1}%
\left(H_{j_{k}}+M_{j_{k-1}}^{(L+1-k+1)}-M_{j_{k}}^{(L+1-k+1)}\right)\Bigg) \\
&  =\max_{i<j_{1}}\Bigg(  H_{j_{1}}+M_{i}^{(L+1)}-M_{j_{1}}^{(L+1)}%\right.
\\
%&\hspace{0.4cm}  \left.
&\quad +\max_{j_{1}<j_{1}^{\prime}<\cdots<j_{L}^{\prime}}%
\sum_{k=1}^{L}\left(H_{j_{k}^{\prime}}+M_{j_{k-1}^{\prime}}^{(L-k+1)}-M_{j_{k}%
^{\prime}}^{(L-k+1)}\right)\Bigg) \\
&  =\max_{i<j_{1}}\left(  H_{j_{1}}+\Theta_{j_{1}+}^{L}+M_{i}^{(L+1)}%
-M_{j_{1}}^{(L+1)}\right) \\
&  =\max_{i<j_{1}}\left(  H_{j_{1}}+\rho_{j_{1}}\left(  Y_{j_{1}+1}^{\ast,
L}\right)  +M_{i}^{(L+1)}-M_{j_{1}}^{(L+1)}\right).
\end{align*}
Next, since $\Theta_{i+}^{L+1}\in\mathcal{F}_{i}$,
Lemma~\ref{th:sure1_cor} implies
\begin{equation*}
\Theta_{i+}^{L+1}=\rho_{i}\left(  Y_{i+1}^{\ast, L+1}\right)  \text{ \ \ and
\ \ }M_{i+1}^{(L+1)}-M_{i}^{(L+1)}=Y_{i+1}^{\ast, L+1}-\rho_{i}\left(
Y_{i+1}^{\ast, L+1}\right).
\end{equation*}
\end{Proof}

\setcounter{equation}{0}

\section{Proofs of Section \ref{sec:MSAL}}

\begin{Proof}
\textbf{of Lemma \ref{var+}.}
Let us define $\mathfrak{m}^{\circ}:=\mathcal{Y}-\rho_{j}\left(  \mathcal{Y}%
\right),$ and write%
\begin{equation*}
\mathfrak{m}^{\circ}-\mathfrak{m}^{N}=\mathcal{Y}-\mathfrak{m}^{N}%
-\mathcal{C}^{N}+\mathcal{C}^{N}-\mathcal{C}+\mathcal{C}-\rho_{j}\left(
\mathcal{Y}\right).
\end{equation*}
Hence, due to (\ref{coy}),%
\begin{equation}
\mathfrak{m}^{\circ}-\mathfrak{m}^{N}\overset{L_{2}}{\rightarrow}%
\mathcal{C}-\rho_{j}\left(  \mathcal{Y}\right)  . \label{he}%
\end{equation}
Then also
\begin{equation*}
\rho_{j}\left(  \mathfrak{m}^{\circ}-\mathfrak{m}^{N}\right)  \overset{L_{2}%
}{\rightarrow}\mathcal{C}-\rho_{j}\left(  \mathcal{Y}\right),
\end{equation*}
since
\begin{align*}
&  \mathbb{E}\left[  \left\vert \rho_{j}\left(  \mathfrak{m}^{\circ
}-\mathfrak{m}^{N}\right)  -\mathcal{C}+\rho_{j}\left(  \mathcal{Y}\right)
\right\vert ^{2}\right] \\
\text{(by (C4))}  &  =\mathbb{E}\left[  \left\vert \rho_{j}\left(
\mathfrak{m}^{\circ}-\mathfrak{m}^{N}-\mathcal{C}+\rho_{j}\left(
\mathcal{Y}\right)  \right)  \right\vert ^{2}\right] \\
\text{(by (\ref{l1}) with }p=2\text{)}&\leq C_{2}\mathbb{E}\left[
\left\vert \mathfrak{m}^{\circ}-\mathfrak{m}^{N}-\mathcal{C}+\rho_{j}\left(
\mathcal{Y}\right)  \right\vert ^{2}\right]  \rightarrow0,
\end{align*}
because of (\ref{he}).
Due to subadditivity and $\rho_{j}\left(\mathfrak{m}^{\circ}\right)=0$,
we have
\begin{equation*}
\rho_{j}\left(  \mathfrak{m}^{\circ}-\mathfrak{m}^{N}\right)  \geq\rho_{j}\left(  \mathfrak{m}^{\circ}\right)  -\rho_{j}\left(  \mathfrak{m}%
^{N}\right)=0.
\end{equation*}
Thus, we must have $\mathcal{C}-\rho_{j}\left(\mathcal{Y}\right)\geq0$.
By the same reasoning,
\begin{align*}
\mathfrak{m}^{N}-\mathfrak{m}^{\circ}%&
\overset{L_{2}}{\rightarrow}\rho_{j}\left(  \mathcal{Y}\right)  -\mathcal{C}\text{ \ \ implies \ \ }%\\
\rho_{j}\left(  \mathfrak{m}^{N}-\mathfrak{m}^{\circ}\right)%&
\overset{L_{2}}{\rightarrow}\rho_{j}\left(  \mathcal{Y}\right)  -\mathcal{C},
\end{align*}
and now subadditivity and $\rho_{j}\left(\mathfrak{m}^{\circ}\right)=0$
implies $\rho_{j}\left(\mathfrak{m}^{N}-\mathfrak{m}^{\circ}\right)\geq0$.
Hence, we must also have that $\rho_{j}\left(  \mathcal{Y}\right)
-\mathcal{C}\geq0$.
Thus, $\rho_{j}\left(\mathcal{Y}\right)=\mathcal{C}$,
and then the other statement follows from (\ref{he}).
\end{Proof}

\vskip 0.3cm

\begin{Proof}
\textbf{of Theorem \ref{theoConv}.}
We will prove the theorem through an induction for $l=1,\ldots,L$.

For $l=1$ we do a second induction over $j=T,T-1,\ldots,0$.
Assume that (\ref{Mconv})--(\ref{Yconv}) hold for $j+1 \leq t$.
It follows from the Law of Large Numbers and the induction assumption for $l-1$ and $j+1$ that %$Y_r^{l,K,N},$
$\overline{M}_{j+1}^{l,K,N}$ and $\overline{c}_{j+1}^{l,K,N}$ converge a.s. to the projections of $M_{j+1}^{\ast,l}$ and $c_{j+1}^{\ast,l}$
on the spaces $\{\mathcal{E}_j^{(\beta_1,\ldots,\beta_{K+1})}\vert (\beta_1,\ldots,\beta_{K+1}) \in \mathbb{R}^{K^{\prime}}\}$
and $\{ \sum_{k=1}^{K^{\prime\prime}} \gamma_k \psi_k (X_j) \vert \gamma_k \in \mathbb{R}, k=1,\ldots, K^{\prime\prime}\}$.
Letting $K=\min(K^{\prime},K^{\prime\prime})$ tend to infinity and using that both spaces form a basis,
we can use Corollary~\ref{loc} to conclude that (\ref{Mconv})--(\ref{Yconv}) hold for $j$.
This completes the induction over $j$, and hence also the induction over $l$.

For simplicity, we drop the indexes $K,N,n$ in the sequel.
So we write $\overline{c}_{j}^{l} = \overline{c}_{j}^{l,K,N} (X_j^n)$.
We let $\overline{c}_{j}^{l}$ and $c_{j}^{\ast, l}$ be a set of approximate and
true continuation functions, respectively,
let
\begin{align*}
\overline{U}_{j}^{l}  %&
=f_{j}(X_{j})+\overline{c}_{j}^{l-1}(X_{j}),%\\
\qquad U_{j}^{\ast, l}  %&
=f_{j}(X_{j})+c_{j}^{\ast, l-1}(X_{j}),
\end{align*}
let $\overline{M}_{j}^{l}$ and $M_{j}^{\ast, l}$ be a set of approximate and true
$\rho$-Doob martingales, and let $\overline{Y}_{j}^{l}$ and $Y_{j}^{\ast, l}$ be a set of
approximate and true upper Snell envelopes.
Then consider
\begin{align*}
&\max_{j\leq r\leq T}\left(  \overline{U}_{r}^{l}-\overline{M}_{r}^{l}\right)
-Y_{j}^{\ast, l}  =\max_{j\leq r\leq T}\left(  \overline{U}_{r}^{l}%
-\overline{M}_{r}^{l}\right)  -\max_{j\leq r\leq T}\left(  U_{r}^{\ast,
l}-M_{r}^{\ast, l}\right)  \\
& =\max_{j\leq r\leq T}\left(  \overline{U}_{r}^{l}-\overline{M}_{r}%
^{l}\right)  -\max_{j\leq r\leq T}\left(  \overline{U}_{r}^{l}-M_{r}^{\ast,
l}\right)
\\
&\qquad
+\max_{j\leq r\leq T}\left(  \overline{U}_{r}^{l}-M_{r}^{\ast, l}\right)
-\max_{j\leq r\leq T}\left(  U_{r}^{\ast, l}-M_{r}^{\ast, l}\right)  \\
& \leq\max_{j\leq r\leq T}\left(  \overline{U}_{r}^{l}-\overline{M}_{r}%
^{l}-\left(  \overline{U}_{r}^{l}-M_{r}^{\ast, l}\right)  \right)
\\
&\qquad
+\max_{j\leq r\leq T}\left(  \overline{U}_{r}^{l}-M_{r}^{\ast, l}-\left(
U_{r}^{\ast, l}-M_{r}^{\ast, l}\right)  \right)  \\
& =\max_{j\leq r\leq T}\left(  M_{r}^{\ast, l}-\overline{M}_{r}^{l}\right)
+\max_{j\leq r\leq T}\left(  \overline{c}_{r}^{l-1}-c_{r}^{\ast, l-1}\right)
\\
&
\leq\max_{j\leq r\leq T}\left\vert M_{r}^{\ast, l}-\overline{M}_{r}%
^{l}\right\vert +\max_{j\leq r\leq T}\,\left\vert \overline{c}_{r}^{l-1}%
-c_{r}^{\ast, l-1}\right\vert .
\end{align*}
Similarly,
\begin{align*}
&Y_{j}^{\ast, l}-\max_{j\leq r\leq T}\left(  \overline{U}_{r}^{l}-\overline
{M}_{r}^{l}\right)    =\max_{j\leq r\leq T}\left(  U_{r}^{\ast, l}%
-M_{r}^{\ast, l}\right)  -\max_{j\leq r\leq T}\left(  \overline{U}_{r}%
^{l}-\overline{M}_{r}^{l}\right)  \\
& =\max_{j\leq r\leq T}\left(  U_{r}^{\ast, l}-M_{r}^{\ast, l}\right)
-\max_{j\leq r\leq T}\left(  \overline{U}_{r}^{l}-M_{r}^{\ast, l}\right)
\\
&\qquad
+\max_{j\leq r\leq T}\left(  \overline{U}_{r}^{l}-M_{r}^{\ast, l}\right)
-\max_{j\leq r\leq T}\left(  \overline{U}_{r}^{l}-\overline{M}_{r}^{l}\right)
\\
& =\max_{j\leq r\leq T}\left(  U_{r}^{\ast, l}-M_{r}^{\ast, l}-\max_{j\leq
r^{\prime}\leq T}\left(  \overline{U}_{r^{\prime}}^{l}-M_{r^{\prime}}^{\ast,
l}\right)  \right)
\\
&\qquad
+\max_{j\leq r\leq T}\left(  \overline{U}_{r}^{l}-M_{r}^{\ast, l}-\max_{j\leq
r^{\prime}\leq T}\left(  \overline{U}_{r^{\prime}}^{l}-\overline{M}%
_{r^{\prime}}^{l}\right)  \right)  \\
& \leq\max_{j\leq r\leq T}\left(  c_{r}^{\ast, l-1}-\overline{c}_{r}%
^{l-1}\right)  +\max_{j\leq r\leq T}\left(  \overline{M}_{r}^{l}-M_{r}^{\ast,
l}\right)
\\
&
\leq\max_{j\leq r\leq T}\left\vert M_{r}^{\ast, l}-\overline{M}_{r}%
^{l}\right\vert +\max_{j\leq r\leq T}\,\left\vert \overline{c}_{r}^{l-1}%
-c_{r}^{\ast, l-1}\right\vert,
\end{align*}
whence%
\begin{equation*}
\left\vert Y_{j}^{\ast, l}-\max_{j\leq r\leq T}\left(  \overline{U}_{r}%
^{l}-\overline{M}_{r}^{l}\right)  \right\vert \leq\max_{j\leq r\leq
T}\left\vert M_{r}^{\ast, l}-\overline{M}_{r}^{l}\right\vert +\max_{j\leq r\leq
T}\,\left\vert \overline{c}_{r}^{l-1}-c_{r}^{\ast, l-1}\right\vert .
\end{equation*}
That is, by monotonicity and subadditivity,
\begin{align}
\rho_{j}&\left(  \left\vert Y_{j}^{\ast, l}-\max_{j\leq r\leq T}\left(
\overline{U}_{r}^{l}-\overline{M}_{r}^{l}\right)  \right\vert \right)\nonumber\\
&\leq \rho_{j}\left(  \max_{j\leq r\leq T}\left\vert M_{r}^{\ast, l}-\overline
{M}_{r}^{l}\right\vert \right)
+\rho_{j}\left(  \max_{j\leq r\leq
T}\,\left\vert \overline{c}_{r}^{l-1}-c_{r}^{\ast, l-1}\right\vert \right).
\label{rhs}%
\end{align}
By the first part of the theorem, the right-hand side in \eqref{rhs} goes to zero.
Thus,
\begin{align*}
&\left\vert Y_{j}^{\ast, l}-\rho_{j}\left(  \max_{j\leq r\leq T}\left(
\overline{U}_{r}^{l}-\overline{M}_{r}^{l}\right)  \right)  \right\vert
= \left\vert \rho_{j}\left(  Y_{j}^{\ast, l}\right)  -\rho_{j}\left(  \max_{j\leq
r\leq T}\left(  \overline{U}_{r}^{l}-\overline{M}_{r}^{l}\right)  \right)
\right\vert  \\
&\leq %C_{\mathrm{Lipschitz}}
\rho_{j}\left(  \left\vert Y_{j}^{\ast,l}
-\max_{j\leq r\leq T}\left(  \overline{U}_{r}^{l}-\overline{M}_{r}^{l}\right) \right\vert \right) \\
&\leq %C_{\mathrm{Lipschitz}}
\rho_{j}\left(  \max_{j\leq r\leq T}\left\vert M_{r}^{\ast, l}-\overline{M}_{r}^{l}\right\vert \right)
+ %C_{\mathrm{Lipschitz}}
\rho_{j}\left(  \max_{j\leq r\leq T}\,\left\vert \overline{c}_{r}^{l-1}-c_{r}^{\ast, l-1}\right\vert \right)
\end{align*}
tends to zero as well.
(Here, the first inequality follows as, by monotoniciy and subadditivity,
$\rho(X)\leq\rho(Y+|X-Y|))\leq\rho(Y)+\rho(|X-Y|)$ yielding
$\rho(X)-\rho(Y)\leq\rho(|X-Y|)$, and switching the roles of $X$ and $Y$ then gives the desired inequality.)
\end{Proof}

\vskip 0.3cm

\begin{Proof}
\textbf{of Proposition \ref{prop:com}.}
We write
\begin{align*}
\Theta_{i}^{q}  &  =\max\left[  \max_{i<j_{2}<\cdots< j_{q}}\left(
f_{i}(X_{i})+\sum_{l=2}^{q}\left(  f_{j_{l}}(X_{j_{l}})-\overline{M}_{j_{l}%
}^{q-l+1}+\overline{M}_{j_{l-1}}^{q-l+1}\right)  \right),  \right. \\
&  \qquad\qquad \left.  \max_{i<j_{1}<j_{2}<\cdots< j_{q}}\sum_{l=1}^{q}\left(  f_{j_{l}%
}(X_{j_{l}})-\overline{M}_{j_{l}}^{q-l+1}+\overline{M}_{j_{l-1}}%
^{q-l+1}\right)  \right] \\
&  =\max\Bigg[  \max_{i+1\leq j_{2}<\cdots< j_{q}}\Bigg(  f_{i}(X_{i}%
)+\overline{M}_{i}^{q-1}-\overline{M}_{i+1}^{q-1}\\
&  \qquad\qquad
+\sum_{l=2}^{q}\left(
f_{j_{l}}(X_{j_{l}})-\overline{M}_{j_{l}}^{q-l+1}+\overline{M}_{j_{l-1}\vee
i+1}^{q-l+1}\right)  \Bigg)  , \\
&  \qquad\qquad   \max_{i+1\leq j_{1}<j_{2}<\cdots< j_{q}}\overline{M}_{i}%
^{q}-\overline{M}_{i+1}^{q}+\sum_{l=1}^{q}\left(  f_{j_{l}}(X_{j_{l}%
})-\overline{M}_{j_{l}}^{q-l+1}+\overline{M}_{j_{l-1}\vee i+1}^{q-l+1}\right)
\Bigg] \\
&  =\max\Bigg[  f_{i}(X_{i})+\overline{M}_{i}^{q-1}-\overline{M}_{i+1}%
^{q-1}\\
&  \qquad\qquad
+\max_{i+1\leq j_{2}<\cdots< j_{q}}\sum_{l=2}^{q}\left(  f_{j_{l}%
}(X_{j_{l}})-\overline{M}_{j_{l}}^{q-l+1}+\overline{M}_{j_{l-1}\vee
i+1}^{q-l+1}\right)  , \\
&\qquad\qquad   \overline{M}_{i}^{q}-\overline{M}_{i+1}^{q}+\max_{i+1\leq
j_{1}<j_{2}<\cdots< j_{q}}\sum_{l=1}^{q}\left(  f_{j_{l}}(X_{j_{l}}%
)-\overline{M}_{j_{l}}^{q-l+1}+\overline{M}_{j_{l-1}\vee i+1}^{q-l+1}\right)
\Bigg]  ,
\end{align*}
which is equal to (\ref{rec}).
\end{Proof}

\setcounter{equation}{0}

\section{Proofs of Section \ref{sec:BPcon}}

\begin{Proof}
\textbf{of Proposition \ref{PropConv}.}
The result follows by analogous arguments as those in the proof of Theorem 21 of Kr\"{a}tschmer \textit{et al.} \cite{KLLSS18}.
\end{Proof}

\vskip 0.3cm

\begin{Proof}
\textbf{of Theorem \ref{th:LB}.}
First, ${\mathbb{E}}\left[\widetilde{Y}_{0}^{\mathrm{low},L}\right]  \leq Y^{\ast,L}_{0}$
follows by \eqref{eq:genuinelowerbound}.
Second, the convergence statement follows by applying Proposition \ref{PropConv} and Theorem \ref{theoConv} three times.
\end{Proof}

\vskip 0.3cm

\begin{Proof}
\textbf{of Proposition \ref{prop:subtracting}.}
We write
\begin{align*}
\mathbb{E}_{\mathbb{Q}}\left[U-\overline{M}^{N}_T\right]&=\mathbb{E}_{\mathbb{Q}}\left[U\right]-\mathbb{E}_{\mathbb{Q}}\left[\overline{M}^{N}_T\right]\\
&=\mathbb{E}_{\mathbb{Q}}\left[U\right]-\mathbb{E}_{\mathbb{Q}}\bigg[\int_0^T \mathcal{Z}^{N}_s dW^{\mathbb{Q}}_s+\int_0^T \mathcal{\tilde{Z}}^{N}_s d\tilde{N}^{\mathbb{Q}}_s\\
&\hspace{1cm}+\int_0^T \bigg\{\mathcal{Z}^{N}_s q_s %^n
+\mathcal{\tilde{Z}}^{N}_s(\lambda_{s} %^{\mathbb{Q},N}
-\lambda_{\mathbb{P}})-g(s,\mathcal{Z}^{N}_s,\mathcal{\tilde{Z}}^{N}_s)\bigg\}ds\bigg]\\
&=\mathbb{E}_{\mathbb{Q}}\left[U\right]+0+0=\mathbb{E}_{\mathbb{Q}}\left[U\right],
\end{align*}
%[***SIGNS???***]
where we used in the one but last equality that the convex conjugate satisfies
\begin{equation*}
\sup_{q,\lambda}%^\mathbb{Q}}
\{zq +\tilde{z}(\lambda%^\mathbb{Q}
-\lambda_{\mathbb{P}})-g(t,z,\tilde{z})\}=0,
\end{equation*}
as $g$ is positively homogeneous.
Moreover, this equality is attained above in $(q_{s}%^N
,\lambda_{s}%^{\mathbb{Q},N}_s
-\lambda_{\mathbb{P}})\in \partial g(s,\mathcal{Z}^{N}_s,\mathcal{\tilde{Z}}^{N}_s)$.
\end{Proof}

\vskip 0.3cm

\begin{Proof}
\textbf{of Theorem \ref{th:UB}.}
First, ${\mathbb{E}}\left[\widetilde{Y}_{0}^{\mathrm{upp},L}\right]\geq Y^{\ast,L}_{0}$
follows by \eqref{inequality}.
Second, the two convergence statements follow by applying Proposition \ref{PropConv} and Theorem \ref{theoConv} three times
and two times, respectively.
\end{Proof}

\vskip 0.3cm

\begin{Proof}
\textbf{of Proposition \ref{prop:Lipschitz}.}
Fix $\alpha\in\mathbb{R}$.
Applying It\^o's generalized formula yields
\begin{align*}
&	e^{\alpha t}|\delta Y_t|^2 +\int_t^T e^{\alpha s}|\delta \mathcal{Z}^{N}_s|^2ds+\sum_{s\geq t:\,\,\delta Y
\mbox{jumps at }s} e^{\alpha s}(|\delta  Y_s|^2-|\delta  Y_{s-}|^2-2\delta Y_{s-}\delta\mathcal{\tilde{Z}}^{N}_s)\\
&\hspace{0.2cm}=e^{\alpha T} |\delta \xi|^2 +\int_t^T e^{\alpha s}\bigg\{2\delta Y_s (g(\mathcal{Z}^{N}_s,\mathcal{\tilde{Z}}^{N}_s)-g(Z'_s,\tilde{Z}'_s)) -\alpha |\delta Y_s|^2 \bigg\}ds\\
&\hspace{0.5cm}-2\int_t^T e^{\alpha s}\delta Y_s \delta \mathcal{Z}^{N}_s dW_s-2\int_t^T e^{\alpha s}\delta Y_s \delta \mathcal{\tilde{Z}}^{N}_s d\tilde{N}_s\\
&\hspace{0.2cm}\leq e^{\alpha T} |\delta \xi|^2 +\int_t^T e^{\alpha s}\bigg\{\mathcal{L}^2 |\delta Y_s|^2+ |\delta \mathcal{Z}^{N}_s|^2+|\delta \mathcal{\tilde{Z}}^{N}_s|^2  -\alpha |\delta Y_s|^2 \bigg\}ds\\
&\hspace{0.5cm}-2\int_t^T e^{\alpha s}\delta Y_s \delta \mathcal{Z}^{N}_s dW_s-2\int_t^T e^{\alpha s}\delta Y_s \delta\mathcal{\tilde{Z}}^{N}_s d\tilde{N}_s,
\end{align*}
using the Lipschitz continuity of $g$ in the equality, and that $2ab\leq \mathcal{L} a^2 +\frac{b^2}{\mathcal{L}}$
where $\mathcal{L}$ is the Lipschitz constant of $g$ in the inequality.
Choosing $\alpha=\mathcal{L}^2$ and observing that
\begin{equation*}
\sum_{s\geq t:\,\,\delta Y  \mbox{ jumps at }s} e^{\alpha s}(|\delta  Y_s|^2-|\delta  Y_{s-}|^2-2\delta Y_{s-}\delta \mathcal{\tilde{Z}}^{N}_s)
=\sum_{s\geq t:\,\,\delta Y \mbox{ jumps at }s} e^{\alpha s}|\delta \mathcal{\tilde{Z}}^{N}_s|^2,
\end{equation*}
(which is the quadratic variation of the jump part of $e^{\alpha s/2}Y_s$)
we obtain, for $t=0$,
\begin{align*}
&	|\delta Y_0|^2 +\int_0^T e^{\alpha s}|\delta \mathcal{Z}^{N}_s|^2ds+\sum_{s\geq 0:\,\,\delta Y  \mbox{ jumps at }s}  e^{\alpha s}|\delta \mathcal{\tilde{Z}}^{N}_s|^2\\	
&\leq e^{\mathcal{L}^2 T} |\delta \xi|^2 +\int_0^T e^{\alpha s}\bigg\{ |\delta \mathcal{Z}^{N}_s|^2+|\delta \mathcal{\tilde{Z}}^{N}_s|^2  \bigg\}ds\\
&\hspace{0.5cm}-2\int_0^T e^{\alpha s}\delta Y_s \delta \mathcal{Z}^{N}_s dW_s-2\int_0^T e^{\alpha s}\delta Y_s \delta \mathcal{\tilde{Z}}^{N}_s d\tilde{N}_s.
\end{align*}
Taking expectations on both sides and cancelling the $\delta \mathcal{Z}$ and $\delta \mathcal{\tilde{Z}}$ terms
corresponding to the quadratic variation
yields the proposition.
\end{Proof}

\setcounter{equation}{0}

\section{Additional Tables}

\begin{table}[H]
\begin{center}
{\small
\begin{tabular}{c||cccc c}
$L$ & 1&2&3&4&5\\
\hline\hline
%    LB without M & $   0.97114 $  & $    1.7365 $  & $    2.3677 $  & $    2.8937 $  & $    3.3332 $  \\
    LB        & $   0.9722 $  & $    1.7389 $  & $    2.3707 $  & $    2.8969 $  & $     3.3350 $  \\
    s.e.      & $ 0.0012 $  & $ 0.0018 $  & $ 0.0023 $  & $ 0.0027 $  & $ 0.0031 $  \\
    $\overline{Y}_0^{N_{4},L}$  & $   0.9869 $  & $    1.7603 $  & $    2.3953 $  & $    2.9235 $  & $    3.3635 $  \\
    TE        & $  0.0682 $  & $   0.1046 $  & $   0.1352 $  & $   0.1631 $  & $   0.1849 $  \\
    UB        & $     1.0550 $  & $    1.8649 $  & $    2.5306 $  & $    3.0866 $  & $    3.5484 $  \\
\end{tabular}
}
\caption{Bounds for $\delta_1=0$, $\delta_2=0$ and $J=0.06$} %(i.e., with jumps)} % with jumps ($\lambda=1$, $J=0.06$)}
\label{tab:multitwocall-1}
%\end{center}
%\end{table}

%\begin{table}[h!]
%{\small
%\begin{tabular}{c||cccc c}
%$L$ & 1&2&3&4&5\\
%\hline\hline
%    $\overline{Y}_0^{N_{4},L}$  & $    0.9955 $  & $    1.7772 $  & $      2.42 $  & $    2.9563 $  & $    3.4037 $  \\
%    LB without M & $   0.97489 $  & $    1.7418 $  & $     2.374 $  & $    2.9024 $  & $    3.3429 $  \\
%    UB        & $    1.0674 $  & $    1.8854 $  & $    2.5598 $  & $    3.1183 $  & $    3.5882 $  \\
%    TE        & $  0.071691 $  & $   0.10795 $  & $   0.13945 $  & $   0.16164 $  & $     0.184 $  \\
%    LB        & $   0.97867 $  & $    1.7496 $  & $    2.3835 $  & $    2.9149 $  & $    3.3653 $  \\
%    s.e.      & $ 0.0012694 $  & $ 0.0019137 $  & $ 0.0024465 $  & $ 0.0029376 $  & $ 0.0034334 $  \\
%\end{tabular}
%}
%\caption{Bounds for $\delta_1=\frac{1}{30}$, $\delta_2=0$ with jumps ($\lambda=1$, $J=0.06$)}
%\end{table}

%\begin{table}[H]
%\begin{center}
\vskip 0.1cm
{\small
\begin{tabular}{c||cccc c}
$L$ & 1&2&3&4&5\\
\hline\hline
%    LB without M & $   0.96991 $  & $    1.7352 $  & $   2.3679 $  & $    2.8945 $  & $    3.3351 $  \\
    LB        & $   0.9884 $  & $    1.7727 $  & $   2.4173 $  & $    2.9625 $  & $    3.4123 $  \\
    s.e.      & $ 0.0015 $  & $ 0.0025 $  & $ 0.0032 $  & $ 0.0040 $  & $ 0.0045 $  \\
   $\overline{Y}_0^{N_{4},L}$  & $     1.0120 $  & $    1.8091 $  & $   2.4672 $  & $    3.0184 $  & $    3.4806 $  \\
    TE        & $   0.0764 $  & $   0.1102 $  & $  0.1365 $  & $   0.1586 $  & $   0.1791 $  \\
    UB        & $    1.0903 $  & $    1.9221 $  & $   2.6071 $  & $     3.1810 $  & $    3.6642 $  \\
\end{tabular}
}
\caption{Bounds for $\delta_1=\tfrac{1}{10}$, $\delta_2=0$ and $J=0.06$} %(i.e., with jumps)} %with jumps ($\lambda=1$, $J=0.06$)}
\label{tab:multitwocall-2}
%\end{center}
%\end{table}

%\begin{table}[h!]
%\begin{center}
\vskip 0.1cm
{\small
\begin{tabular}{c||cccc c}
$L$ & 1&2&3&4&5\\
\hline\hline
%    LB without M & $   0.97287 $  & $     1.742 $  & $    2.3753 $  & $     2.904 $  & $    3.3441 $  \\
    LB        & $   0.9999 $  & $    1.8066 $  & $     2.4640 $  & $    3.0228 $  & $    3.4831 $  \\
    s.e.      & $ 0.0021 $  & $ 0.0037 $  & $ 0.0049 $  & $ 0.0061 $  & $ 0.0071 $  \\
    $\overline{Y}_0^{N_{4},L}$  & $     1.0370 $  & $    1.8587 $  & $    2.5407 $  & $     3.1150 $  & $    3.5996 $  \\
    TE        & $  0.0747 $  & $   0.1095 $  & $   0.1413 $  & $   0.1695 $  & $   0.1942 $  \\
    UB        & $    1.1196 $  & $    1.9797 $  & $    2.6969 $  & $    3.3023 $  & $    3.8142 $  \\
\end{tabular}
}
\caption{Bounds for $\delta_1=\tfrac{1}{5}$, $\delta_2=0$ and $J=0.06$} %(i.e., with jumps)} %with jumps ($\lambda=1$, $J=0.06$)}
\label{tab:multitwocall-3}
\end{center}
\end{table}

\begin{table}[H]
\begin{center}
\vskip 0.1cm
{\small
\begin{tabular}{c||cccc c}
$L$ & 1&2&3&4&5\\
\hline\hline
%    LB without M & $   0.97225 $  & $    1.7392 $  & $    2.3709 $  & $    2.8986 $  & $    3.3378 $  \\
    LB        & $   0.9776 $  & $    1.7484 $  & $    2.3814 $  & $    2.9123 $  & $    3.3529 $  \\
    s.e.      & $ 0.0020 $  & $ 0.0033 $  & $ 0.0044 $  & $ 0.0054 $  & $ 0.0062 $  \\
    $\overline{Y}_0^{N_{4},L}$  & $    1.0028 $  & $    1.7907 $  & $    2.4386 $  & $     2.9790 $  & $    3.4302 $  \\
    TE        & $  0.0657 $  & $   0.1003 $  & $   0.1301 $  & $    0.1569 $  & $   0.1807 $  \\
    UB        & $    1.0754 $  & $    1.9015 $  & $    2.5823 $  & $    3.1524 $  & $      3.6300 $  \\
\end{tabular}
}
\caption{Bounds for $\delta_1=0$, $\delta_2=\tfrac{1}{5}$ and $J=0.06$} %(i.e., with jumps)} %with jumps ($\lambda=1$, $J=0.06$)}
\label{tab:multitwocall-4}
%\end{center}
%\end{table}

%\begin{table}[h!]
%{\small
%\begin{tabular}{c||cccc c}
%$L$ & 1&2&3&4&5\\
%\hline\hline
%    $\overline{Y}_0^{N_{4},L}$  & $    1.0121 $  & $    1.8079 $  & $    2.4639 $  & $    3.0118 $  & $    3.4704 $  \\
%    LB without M & $   0.97356 $  & $    1.7417 $  & $    2.3724 $  & $       2.9 $  & $    3.3394 $  \\
%    UB        & $    1.0948 $  & $    1.9317 $  & $    2.6228 $  & $    3.1986 $  & $    3.6849 $  \\
%    TE        & $  0.074866 $  & $   0.11203 $  & $   0.14375 $  & $   0.16896 $  & $   0.19406 $  \\
%    LB        & $   0.97615 $  & $     1.752 $  & $    2.3844 $  & $    2.9173 $  & $    3.3604 $  \\
%    s.e.      & $ 0.0020518 $  & $ 0.0034508 $  & $ 0.0045698 $  & $ 0.0055345 $  & $ 0.0063346 $  \\
%\end{tabular}
%}
%\caption{Bounds for $\delta_1=\frac{1}{30}$, $\delta_2=\frac15$ with jumps ($\lambda=1$, $J=0.06$)}
%\end{table}

%\begin{table}[H]
%\begin{center}
\vskip 0.1cm
{\small
\begin{tabular}{c||cccc c}
$L$ & 1&2&3&4&5\\
\hline\hline
%    LB without M & $   0.97532 $  & $    1.7442 $  & $    2.3789 $  & $    2.9079 $  & $    3.3493 $  \\
    LB        & $   0.9929 $  & $    1.7811 $  & $    2.4296 $  & $    2.9688 $  & $    3.4244 $  \\
    s.e.      & $ 0.0023 $  & $ 0.0038 $  & $ 0.0052 $  & $ 0.0063 $  & $ 0.0072 $  \\
    $\overline{Y}_0^{N_{4},L}$  & $    1.0275 $  & $    1.8398 $  & $    2.5113 $  & $    3.0748 $  & $    3.5483 $  \\
    TE        & $  0.0740 $  & $   0.1106 $  & $   0.1395 $  & $   0.1654 $  & $   0.1875 $  \\
    UB        & $    1.1093 $  & $     1.9620 $  & $    2.6655 $  & $    3.2576 $  & $    3.7555 $  \\
\end{tabular}
}
\caption{Bounds for $\delta_1=\tfrac{1}{10}$, $\delta_2=\tfrac{1}{5}$ and $J=0.06$} %(i.e., with jumps)} %with jumps ($\lambda=1$, $J=0.06$)}
\label{tab:multitwocall-5}
%\end{center}
%\end{table}

%\begin{table}[h!]
%\begin{center}
\vskip 0.1cm
{\small
\begin{tabular}{c||cccc c}
$L$ & 1&2&3&4&5\\
\hline\hline
%    LB without M & $   0.97226 $  & $    1.7412 $  & $    2.3745 $  & $    2.9032 $  & $   3.3443 $  \\
    LB        & $    1.0026 $  & $    1.8031 $  & $    2.4755 $  & $    3.0091 $  & $   3.4702 $  \\
    s.e.      & $ 0.0028 $  & $ 0.0049 $  & $ 0.0067 $  & $ 0.0081 $  & $ 0.0094 $  \\
    $\overline{Y}_0^{N_{4},L}$  & $    1.0536 $  & $    1.8896 $  & $    2.5848 $  & $    3.1715 $  & $   3.6676 $  \\
    TE        & $  0.0749 $  & $   0.1097 $  & $   0.1420 $  & $   0.1699 $  & $  0.1947 $  \\
    UB        & $    1.1365 $  & $    2.0109 $  & $    2.7418 $  & $    3.3592 $  & $   3.8828 $  \\
\end{tabular}
}
\caption{Bounds for $\delta_1=\tfrac{1}{5}$, $\delta_2=\tfrac{1}{5}$ and $J=0.06$} %(i.e., with jumps)} %with jumps ($\lambda=1$, $J=0.06$)}
\label{tab:multitwocall-6}
\end{center}
\end{table}

\end{document}